\documentclass[reqno]{amsart}
\usepackage[margin=1.2in]{geometry}
\usepackage{amssymb}
\usepackage{mathrsfs}
\usepackage{float}
\usepackage{tikz}
\usetikzlibrary{arrows.meta}
\usepackage[colorlinks=true]{hyperref}

\hypersetup{
    pdftitle={Energy martingale solutions for stochastic mean-curvature flow},
    pdfauthor={Zhengyan Wu and Ao Zhang},
    pdfsubject={Stochastic mean-curvature flow of periodic graphs},
    pdfkeywords={stochastic mean curvature flow, Kawasaki--Ohta equation, thermal fluctuations, energy martingale solutions},
    linkcolor=blue,
    citecolor=red,
    filecolor=blue,
    urlcolor=cyan
}

\numberwithin{equation}{section}

\newtheorem{theorem}{Theorem}[section]
\newtheorem{lemma}[theorem]{Lemma}
\newtheorem{proposition}[theorem]{Proposition}
\theoremstyle{definition}
\newtheorem{definition}[theorem]{Definition}
\newtheorem{assumption}{Assumption}[section]
\theoremstyle{remark}
\newtheorem{remark}{Remark}
\theoremstyle{plain}

\newenvironment{nouppercase}{%
  \renewcommand{\uppercasenonmath}[1]{}}{}

\allowdisplaybreaks

\begin{document}

\title[Fluctuating Mean Curvature Flow]{\LARGE Mean Curvature Flow with Thermal Fluctuations in Dimensions Two and Higher}

\author[Zhengyan Wu]{\large Zhengyan Wu}
\address[Z. Wu]{Department of Mathematics, Technische Universit\"at M\"unchen, Boltzmannstr. 3, 85748 Garching, Germany}
\email{wuzh@cit.tum.de} 

\author[Ao Zhang]{\large Ao Zhang}
\address[A. Zhang]{School of Mathematics and Statistics, HNP-LAMA, Central South University, Changsha 410083, China}
\email{aozhang1993@csu.edu.cn}

\begin{abstract}
Surface tension drives interfaces towards smaller area, while thermal
motion produces shape fluctuations whose strength is linked to the
interfacial mobility through the fluctuation--dissipation relation.
These mechanisms lead to mean-curvature flow with thermal noise,
whose graph formulation is known as the Kawasaki--Ohta equation.
We establish global-in-time existence of energy martingale solutions
to the spatially regularized Kawasaki--Ohta equation for periodic
graphs in every dimension $d\geq2$, with initial data in
$H^1(\mathbb T^d)$. No restriction on the noise intensity is needed
for $d=2$, while a smallness condition is required for $d\geq3$.
The key new challenges arise from the
thermal noise, which has a $(1+|\nabla u|^2)^{1/4}$-type coefficient.
This coefficient gives rise to an irregular Stratonovich-to-It\^o
correction term of the form
$\nabla u\cdot D^2u\nabla u/(1+|\nabla u|^2)^{3/2}$.
Our main innovation is a compactness argument tailored to the structure
of the correction. By combining a factorization of this term with
$H^1(\mathbb T^d)$ and area estimates, we establish strong convergence
of the normalized gradients and reconstruct the weighted Hessian,
which allows us to define the correction at energy regularity and
identify its approximation limit.
\end{abstract}


\date{\today}

\begin{nouppercase}
\maketitle
\end{nouppercase}

%\setcounter{tocdepth}{1}
%\tableofcontents

%%%%%%%%%%%%%%%
\section{Introduction}

Surface tension and thermal fluctuations both shape the motion of
interfaces in phase transitions.  For an interface driven only by
isotropic surface tension, constant normal mobility leads to
mean-curvature flow.  The deterministic mean-curvature flow reduces interfacial area but leaves
out the random shape fluctuations caused by thermal motion.
One can include these fluctuations by adding noise whose strength is linked
to the interfacial mobility through the formal fluctuation--dissipation
relation.  For interfaces represented
as graphs, the Kawasaki and Ohta~\cite{KawasakiOhta1982} give such a description, also see \cite[Eq.~(3.15)]{KatsoulakisKho2001} for the
graph derivation.  Both the
curvature drift and thermal noise depend on the slope of the height
function.  The resulting equation is quasilinear, with degenerate
diffusion and gradient-dependent noise.  In this paper, we prove
global-in-time existence of energy martingale solutions to a spatially
regularized Kawasaki--Ohta equation for periodic graphs in every
dimension $d\geq2$.

Let the spatial domain be the flat torus $\mathbb T^d$,
where $d\geq2$.
We study the graph-valued Stratonovich equation
\begin{align}\label{SPDE-0}
 du={}&(1+|\nabla u|^2)^{1/2}\,\nabla\cdot
 \left(\frac{\nabla u}{(1+|\nabla u|^2)^{1/2}}\right)dt
 +\sqrt{\varepsilon}\,(1+|\nabla u|^2)^{1/4}\circ dW_K .
\end{align}
Here $u$ is the height of a periodic graph and $\varepsilon>0$ is the
noise intensity.  The formal time derivative of the Wiener process
$W_K$ is a spatial ultraviolet cutoff of space--time white noise.
We specify its finite-mode representation and covariance assumptions
in Subsection~\ref{subsec-intro-results}.

Equation~\eqref{SPDE-0} is a spatially regularized version of the
Kawasaki--Ohta model for thermally fluctuating interfaces
\cite{KatsoulakisKho2001,KawasakiOhta1982}.
The noise amplitude $(1+|\nabla u|^2)^{1/4}$ is the square root of the
graph mobility $(1+|\nabla u|^2)^{1/2}$, as prescribed by the formal
fluctuation--dissipation relation for local thermal forcing.
We provide a derivation of the model in Subsection~\ref{subsec-thermal-background}.
In fluctuating hydrodynamics, stochastic PDEs for thermal fluctuations
are often viewed as mesoscopic descriptions of coarse-grained
microscopic dynamics.  This motivates replacing space--time white
noise by noise that is white in time and correlated in space, with
a correlation length reflecting the spatial resolution.  Such regularized
equations are often written in Stratonovich form.  The
Stratonovich-to-It\^o correction can give stochastic coercivity through
cancellations in entropy-dissipation estimates; see
\cite[Section~5.1]{FehrmanGess2024}.  Under suitable structural
conditions on the multiplicative noise, the Stratonovich or
Klimontovich interpretation can also restore time-reversibility; see
Ayala, Dirr, Pavliotis, and Zimmer~\cite{AyalaDirrPavliotisZimmer2025}.\footnote{For certain one-dimensional Langevin models with state-dependent noise, the Klimontovich,
rather than Stratonovich, convention yields reversibility with respect
to the intended Gibbs measure; see
\cite[Section~1.1]{AyalaDirrPavliotisZimmer2025}.
For the same smooth noise coefficients, the two It\^o corrections
differ by a factor of two.  This affects reversibility and may change
constants and admissible noise ranges, while generally leaving the
basic analytical framework unchanged.  We use the Stratonovich
convention for the analysis here and do not claim equivalence of the
two dynamics or reversibility of the present model.}
Macroscopic fluctuation limits can then be studied by sending the
noise intensity and spatial correlation length to zero under a
suitable joint scaling.  For related models, this recovers fluctuation
statistics predicted by microscopic theories; see
\cite{FehrmanGess2023} and \cite[Section~1.2]{FehrmanGess2024}.
The first step is a solution theory for the spatially regularized
equations, which we establish in this paper.

The difficulties in \eqref{SPDE-0} depend on the spatial dimension
and noise amplitude.  In one dimension, formal differentiation with
$\rho=\partial_xu$ yields the scalar conservative SPDE
\begin{equation*}
 d\rho=\partial_{xx}(\arctan\rho)\,dt
 +\sqrt{\varepsilon}\,\partial_x
 \left((1+\rho^2)^{1/4}\circ dW_K\right).
\end{equation*}
Fehrman and Gess \cite[Section~1.2, Example~1.7]{FehrmanGess2024}
study nonlinear fluctuating-hydrodynamics equations of this type and
prove well-posedness under the assumptions stated there.\footnote{The
cited example concerns nonnegative
$L^3(\mathbb T)$ initial slopes.  An arbitrary periodic height has a signed,
mean-zero derivative, so the result does not directly give well-posedness
for arbitrary periodic heights.}
For $d\geq2$, the slope $\nabla u$ is vector valued, therefore the scalar reduction is no longer available.

For stochastic mean-curvature equations of graphs in higher dimensions,
existence was proved by
Hofmanov\'a, R\"oger, and von Renesse~\cite{HofmanovaRogerVonRenesse2017}
and Dabrock, Hofmanov\'a, and R\"oger~\cite{DabrockHofmanovaRoger2021}.
Subsection~\ref{subsec-thermal-background} explains the different
physical origins of their noise and ours.  Replacing their amplitude
$(1+|\nabla u|^2)^{1/2}$ by the thermal amplitude
$(1+|\nabla u|^2)^{1/4}$ changes the Stratonovich-to-It\^o correction
and creates new difficulties at the
energy regularity.   The It\^o
formulation derived in Section~\ref{sec-ito-formulation} shows that the correction reads 
\begin{equation*}
 \frac{\varepsilon\Lambda_K}{8}\,
 \frac{\nabla u\cdot D^2u\,\nabla u}{(1+|\nabla u|^2)^{3/2}}.
\end{equation*}
Here the noise covariance constant $\Lambda_K$ is specified in
Subsection~\ref{subsec-intro-results}.  Unlike the correction terms in
\cite{DabrockHofmanovaRoger2021,HofmanovaRogerVonRenesse2017},
this term cannot be written as a constant linear combination of
the Laplacian and the original mean-curvature drift. Defining the correction requires control of
the weighted Hessian, which the weighted mean-curvature bound from
the area estimate does not provide.  We develop new arguments to
define this term for weak solutions and identify it in the approximation
limit.  Their treatment is a central contribution of
this paper, and the main technical ideas are explained in
Subsection~\ref{subsec-intro-techniques}.

\subsection{Main results}\label{subsec-intro-results}

We first specify the noise and fix notation.  Let $K\in\mathbb N$ be fixed,
let $f_1,\ldots,f_K\in C^\infty(\mathbb T^d;\mathbb R)$, and let
$\beta^1,\ldots,\beta^K$ be independent real Brownian motions.  Write
$\beta:=(\beta^k)_{k=1,\ldots,K}$ for the standard $K$-dimensional
Brownian motion and set
\begin{equation*}
 W_K(t,x):=\sum_{k=1}^K f_k(x)\beta_t^k.
\end{equation*}
Retaining finitely many modes of a real Fourier basis gives the usual
spatial ultraviolet cutoff of space--time white noise.  We also allow
more general smooth modes satisfying the constant
pointwise variance condition
\begin{equation*}
 \sum_{k=1}^K f_k(x)^2\equiv\Lambda_K,
 \qquad \Lambda_K\geq0.
\end{equation*}
For every $p\in\mathbb R^d$, define
\begin{equation}
 H(p):=(1+|p|^2)^{1/2}.
\end{equation}
We also write
\begin{equation*}
 \kappa(u):=\nabla\cdot\left(\frac{\nabla u}{H(\nabla u)}\right),
 \qquad
 \mathcal A(u):=\int_{\mathbb T^d}H(\nabla u)\,dx,
 \qquad
 \lambda:=\frac{\varepsilon\Lambda_K}{8}.
\end{equation*}
Here $\kappa(u)$ is the mean curvature with the chosen upward normal,
and $\mathcal A(u)$ is the graph area.

We prove the following existence result.

\begin{theorem}\label{thm-target-existence}
Under the above assumptions on the noise, fix $T\in(0,\infty)$.
Assume $\lambda<(d+1)/(2(d-1))$ if $d\geq3$, with no restriction on
$\lambda$ if $d=2$.
Let $\Lambda$ be a Borel probability measure on
$H^1(\mathbb T^d)$ such that
\begin{equation*}
 \int_{H^1(\mathbb T^d)}\|z\|_{H^1(\mathbb T^d)}^2\,\Lambda(dz)<\infty.
\end{equation*}
Then the It\^o formulation~\eqref{SPDE-Ito} of~\eqref{SPDE-0} admits
energy martingale solutions on $[0,T]$ with initial law $\Lambda$, in
the sense of Definition~\ref{def-energy-martingale-solution}, satisfying
\begin{equation}\label{theorem-area-bound}
\begin{aligned}
 \mathbb E\!\left[\operatorname*{ess\,sup}_{0<t<T}
 \mathcal A(u(t))^2\right]
 +\mathbb E\!\left[\left(\int_0^T\!\int_{\mathbb T^d}
 H(\nabla u)|\kappa(u)|^2\,dx\,dt\right)^2\right]\leq C_T\left(1+\int_{H^1(\mathbb T^d)}
 \|z\|_{H^1(\mathbb T^d)}^2\,\Lambda(dz)\right),
\end{aligned}
\end{equation}
and
\begin{equation}\label{theorem-weighted-hessian-bound}
\begin{aligned}
 \operatorname*{ess\,sup}_{0<t<T}\mathbb E
 \|\nabla u(t)\|_{L^2(\mathbb T^d)}^2
 +\mathbb E\int_0^T\!\int_{\mathbb T^d}
 \left|\frac{D^2u}{H(\nabla u)}\right|^2\,dx\,dt\leq C_T\left(1+\int_{H^1(\mathbb T^d)}
 \|z\|_{H^1(\mathbb T^d)}^2\,\Lambda(dz)\right).
\end{aligned}
\end{equation}
The constant $C_T$ may depend on $d$, $T$, $\varepsilon$,
$K$, and $\lambda$ when $d\geq3$,
but not on the initial law $\Lambda$.
\end{theorem}

\begin{remark}
The solutions are weak in both the analytic and probabilistic sense of
Definition~\ref{def-energy-martingale-solution}. 
\end{remark}

\subsection{Background of the model and the thermal noise}
\label{subsec-thermal-background}

A typical example behind \eqref{SPDE-0} is a moving boundary between
magnetic domains with positive and negative spins.  In a phase-field
description, the coarse-grained order parameter $\phi(t,X)$ is close to
$+1$ or $-1$ in the two bulk phases of a material domain $\mathcal D$ and
varies rapidly across a thin transition layer.  The evolution of local spin flips
can be modelled by
nonconserved Allen--Cahn or time-dependent Ginzburg--Landau dynamics.
As the layer thickness tends to zero, the sharp-interface description
replaces the diffuse layer by a moving hypersurface $\Gamma_t$.

For isotropic surface tension and in the absence of a bulk bias, the
deterministic Allen--Cahn sharp-interface limit yields mean-curvature
motion, which reduces the interfacial area.  At finite temperature,
thermal fluctuations perturb this motion, as described by the
Kawasaki--Ohta model~\cite{KawasakiOhta1982}.
Katsoulakis and Kho~\cite[Eq.~(3.15)]{KatsoulakisKho2001}
formally derive the corresponding stochastic graph equation from
stochastic Ginzburg--Landau dynamics.  Equation~\eqref{SPDE-0} is its
spatially regularized version.  We outline the geometry and
fluctuation--dissipation relation that lead to the formal height
equation, and explain how the local thermal forcing differs from
the spatially coherent normal forcing in
\cite{DabrockHofmanovaRoger2021,HofmanovaRogerVonRenesse2017}.

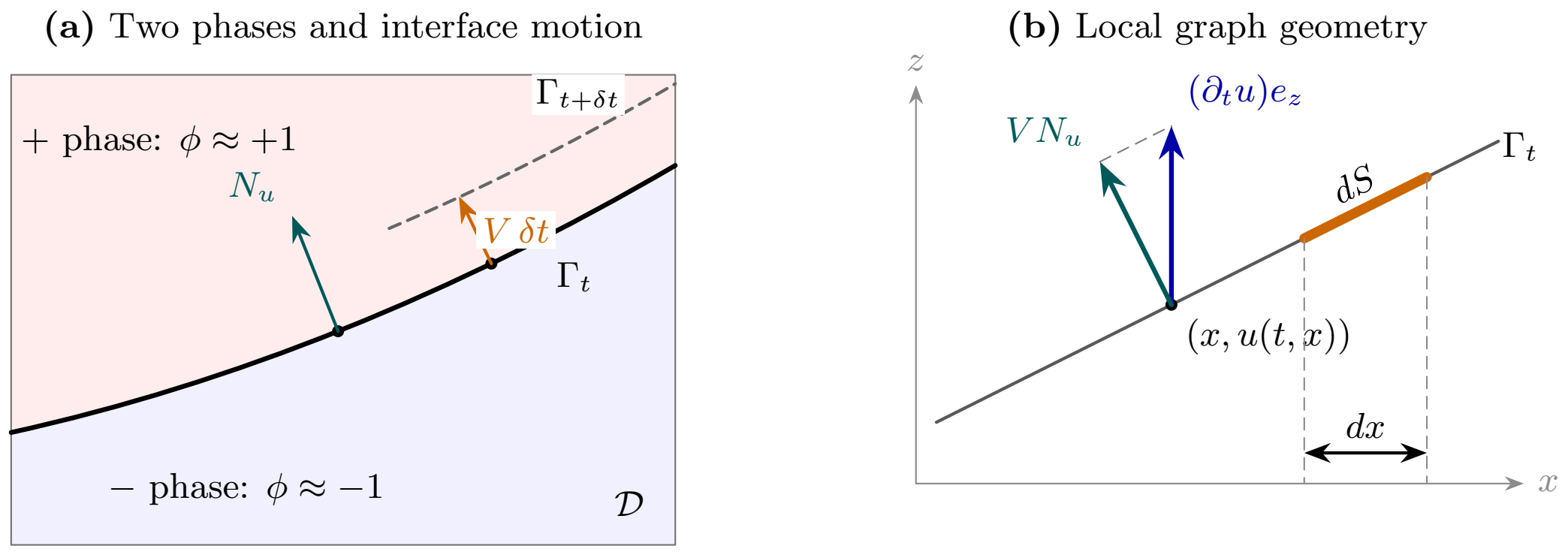
\begin{figure}[H]
\centering
\begin{tikzpicture}[x=0.9cm,y=0.9cm,font=\small,
    >=Stealth,line cap=round,line join=round]
  % A two-dimensional section through the two-phase material.
  \begin{scope}
    \fill[blue!6] (0,0) rectangle (6.5,4.6);
    \fill[red!7]
      (0,1.1) .. controls (2.1667,1.5767) and (4.3333,2.4482) ..
      (6.5,3.713) -- (6.5,4.6) -- (0,4.6) -- cycle;
    \draw[black!55] (0,0) rectangle (6.5,4.6);
    \draw[very thick]
      (0,1.1) .. controls (2.1667,1.5767) and (4.3333,2.4482) ..
      (6.5,3.713);
    \node at (1.45,3.95) {$+$ phase: $\phi\approx+1$};
    \node at (2.3,0.55) {$-$ phase: $\phi\approx-1$};
    \node at (6.05,0.4) {$\mathcal D$};
    \node[below right] at (5.2,2.9) {$\Gamma_t$};
    \fill (3.2,2.0907) circle (1.5pt);
    \draw[->,thick,teal!70!black]
      (3.2,2.0907) -- (2.75,3.22) node[above left] {$N_u$};
    \draw[densely dashed,thick,black!60]
      plot[domain=3.7:6.5,samples=30]
      (\x,{1.9+0.22*\x+0.028*\x*\x});
    \node[above,fill=white,inner sep=1pt] at (5.55,4.2)
      {$\Gamma_{t+\delta t}$};
    \fill (4.7,2.7525) circle (1.5pt);
    \draw[->,thick,orange!80!black]
      (4.7,2.7525) -- (4.382,3.411)
      node[midway,right,fill=white,inner sep=1.5pt] {$V\,\delta t$};
    \node at (3.25,5.05) {\textbf{(a)} Two phases and interface motion};
  \end{scope}
  % Local tangent geometry and the projection of the vertical velocity.
  \begin{scope}[xshift=7.7cm]
    \draw[->,black!45] (0.3,0.6) -- (6.25,0.6) node[right] {$x$};
    \draw[->,black!45] (0.3,0.6) -- (0.3,4.5) node[above] {$z$};
    \draw[thick,black!65] (0.5,1.2) -- (6,3.95);
    \node[right] at (5.9,3.9) {$\Gamma_t$};
    \fill (2.8,2.35) circle (1.5pt);
    \node[below right] at (2.8,2.35) {$(x,u(t,x))$};
    \draw[->,very thick,blue!65!black]
      (2.8,2.35) -- (2.8,4.1) node[above right] {$(\partial_tu)e_z$};
    \draw[->,very thick,teal!70!black]
      (2.8,2.35) -- (2.1,3.75) node[above left] {$V N_u$};
    \draw[densely dashed,black!50] (2.1,3.75) -- (2.8,4.1);
    \draw[line width=2.5pt,orange!80!black] (4.1,3) -- (5.3,3.6);
    \node[above,rotate=26.565] at (4.7,3.3) {$dS$};
    \draw[densely dashed,black!45] (4.1,3) -- (4.1,0.6);
    \draw[densely dashed,black!45] (5.3,3.6) -- (5.3,0.6);
    \draw[<->,thick] (4.1,0.9) -- (5.3,0.9)
      node[midway,above] {$dx$};
    \node at (3.25,5.05) {\textbf{(b)} Local graph geometry};
  \end{scope}
  \node at (7.55,-0.65)
    {$dS=H(\nabla u)\,dx,\qquad
      V=((\partial_tu)e_z)\cdot N_u=\dfrac{\partial_tu}{H(\nabla u)}$};
\end{tikzpicture}
\caption{Two-phase interface motion (a) and local graph geometry (b).}
\label{fig-interface-geometry}
\end{figure}

We write the interface as a graph with height
$u(t,x)\in\mathbb R$ over $x\in\mathbb T^d$,
$\Gamma_t=\Gamma_{u(t)}=\{(x,u(t,x)):x\in\mathbb T^d\}$.
Here $u$ gives the interface location, while $\phi$ describes the
phases in the material.  The formal calculations below are for smooth
graphs, with the geometry shown in Figure~\ref{fig-interface-geometry}.
We denote the upward unit normal by $N_u$, the
surface-area element by $dS$, the corresponding base-volume
element by $dx$, and the vertical unit vector by $e_z$. 

At fixed $x$, the point $(x,u(t,x))$ has vertical velocity
$(\partial_tu)e_z$, whose normal component is the signed interface
velocity $V$.  Over a short time $\delta t$, the normal displacement
is $V\delta t$ to first order, and $V>0$ means motion along $N_u$.
In graph coordinates, we have
\[
 N_u=\frac{(-\nabla u,1)}{H(\nabla u)},
 \qquad dS=H(\nabla u)\,dx,
 \qquad V=\frac{\partial_tu}{H(\nabla u)}.
\]
The factor $H=(1+|\nabla u|^2)^{1/2}$ both normalizes the perpendicular
vector $(-\nabla u,1)$ and gives the area Jacobian of the graph
parametrization.  The normal component of the vertical velocity is
$V=((\partial_tu)e_z)\cdot N_u=\partial_tu/H$, so the corresponding
height evolution is $\partial_tu=HV$.

Let $\gamma>0$ denote the surface tension and $\mu>0$ the normal
mobility.  The surface energy and its first variation in height
coordinates are
\begin{equation*}
 \mathcal E(u)=\gamma\mathcal A(u),
 \qquad \frac{\delta\mathcal E}{\delta u}=-\gamma\kappa(u).
\end{equation*}
For constant normal mobility, the height mobility operator is
$M(u)\psi:=\mu H(\nabla u)\psi$, and
\begin{equation*}
 \partial_tu=-M(u)\frac{\delta\mathcal E}{\delta u}
 =\mu\gamma H(\nabla u)\kappa(u),
 \qquad
 \frac{d}{dt}\mathcal E(u)
 =-\mu\gamma^2\int_{\Gamma_u}\kappa(u)^2\,dS.
\end{equation*}

Write $\Theta=k_B\vartheta$ for the thermal energy, with Boltzmann
constant $k_B$ and absolute temperature $\vartheta$.
Local thermal forcing with constant normal mobility has the formal law
\begin{equation}\label{intro-normal-thermal-law}
 V=\mu\gamma\kappa+\sqrt{2\mu\Theta}\,\xi_\Gamma,
\end{equation}
where $\xi_\Gamma$ is white noise with respect to surface area and
time, $dS\,dt$.  To express this noise in graph coordinates, first
fix a smooth graph $u$ and let $\xi_x$ be white noise with respect to
$dx\,dt$.  Since $dS=H(\nabla u)\,dx$, the formal noise densities satisfy
\begin{equation*}
 \xi_\Gamma(t,(x,u(x)))
 \overset{\mathrm{law}}{=}
 H(\nabla u(x))^{-1/2}\xi_x(t,x).
\end{equation*}
Using $\partial_tu=HV$ in \eqref{intro-normal-thermal-law} now gives
\begin{equation}\label{intro-height-thermal-law}
 \partial_tu=\mu\gamma H(\nabla u)\kappa(u)
 +\sqrt{2\mu\Theta}\,H(\nabla u)^{1/2}\xi_x.
\end{equation}
Normalizing $\mu\gamma=1$ and
writing $\varepsilon=2\mu\Theta$ gives the local-noise counterpart of
\eqref{SPDE-0}; see also \cite[Eq.~(3.15)]{KatsoulakisKho2001}.

In contrast, a single Brownian forcing acting equally at all points
in the normal direction gives
\begin{equation*}
 V\,dt=\kappa\,dt+\circ d\beta_t,
 \qquad
 du=H(\nabla u)\kappa(u)\,dt+H(\nabla u)\circ d\beta_t.
\end{equation*}
This is the forcing considered in
\cite{DabrockHofmanovaRoger2021,HofmanovaRogerVonRenesse2017}.
This normal forcing is spatially homogeneous, since the same Brownian
motion acts at every point of the interface. In contrast, the thermal
noise in \eqref{intro-height-thermal-law} describes local fluctuations
normalized by surface area.
The different noise amplitudes also lead to different
Stratonovich-to-It\^o corrections. In the following, we will explain the new challenges and our key ideas. 

\subsection{Key ideas and technical comments}\label{subsec-intro-techniques}

We first introduce a viscous approximation equation, see \eqref{viscous-equation}. Existence of weak solutions to this viscous approximation model can be shown by the standard Galerkin method for fixed $\nu>0$, based on the
$H^1(\mathbb T^d)$ energy estimate. The key challenge is to identify the nonlinear correction term
in the vanishing-viscosity limit $\nu\to0$.   

To pass to zero viscosity, we combine the $H^1(\mathbb T^d)$ energy estimate
\begin{equation}\label{intro-viscous-dirichlet-estimate}
\begin{aligned}
 \sup_{0\leq t\leq T}\mathbb E
 \|\nabla u^\nu(t)\|_{L^2(\mathbb T^d)}^2
 +\mathbb E\int_0^T\!\int_{\mathbb T^d}
 \frac{|D^2u^\nu|^2}{H(\nabla u^\nu)^2}\,dx\,dt+\nu\mathbb E\int_0^T
 \|\Delta u^\nu\|_{L^2(\mathbb T^d)}^2\,dt
 \leq C(T,u_0),
\end{aligned}
\end{equation}
with the area estimate
\begin{equation}\label{intro-viscous-area-estimate}
\begin{aligned}
 \mathbb E\left[\sup_{0\leq t\leq T}\mathcal A(u^\nu(t))^2\right]
 +\mathbb E\left[\left(\int_0^T\!\int_{\mathbb T^d}
 H(\nabla u^\nu)|\kappa(u^\nu)|^2\,dx\,dt\right)^2\right]\leq C(T,u_0),
\end{aligned}
\end{equation}
with constants independent of $0<\nu\leq1$.
The first follows from the Galerkin bounds by lower semicontinuity; see
Proposition~\ref{lem-fixed-nu-energy-limit} and
Theorem~\ref{thm-viscous-equation}.
The second follows from It\^o's formula for the graph area and a
martingale estimate; see
Proposition~\ref{prop-viscous-area-identity} and
Lemma~\ref{lem-second-moment-area}.
Here we remark that their derivation differs for $d=2$ and $d\geq3$.  In dimension two,
identities for $2\times2$ matrices and vanishing periodic weighted
Hessian-determinant integrals cancel the Stratonovich-to-It\^o
correction against the Hessian-dependent quadratic-variation terms
after spatial integration, therefore no smallness condition on $\lambda$ is needed. However, for $d\geq3$, quadratic Hessian integrals remain.
Weighted Newton identities bound them by the deterministic dissipation,
allowing absorption under \eqref{dimension-noise-assumption}.
The identities and coercivity estimates are proved in
Section~\ref{sec-geometric-estimates}.

At the energy regularity, we still need to define and identify the
Stratonovich-to-It\^o correction.  We use the three-factor decomposition
\begin{equation*}
 \frac{\nabla u\cdot D^2u\,\nabla u}{H(\nabla u)^3}
 =\frac{\nabla u}{H(\nabla u)}\cdot
 \frac{D^2u}{H(\nabla u)}
 \frac{\nabla u}{H(\nabla u)}.
\end{equation*}
The two outer factors coincide and have norm at most one. Noting that weak convergence of the gradients alone does not give strong
convergence of the normalized gradients, a special structure of $H(p)=(1+|p|^2)^{1/2}$ plays an important role through the quantitative monotonicity inequality
\begin{equation*}
 \left|\frac{p}{H(p)}-\frac{q}{H(q)}\right|^2
 \leq
 \left(\frac{p}{H(p)}-\frac{q}{H(q)}\right)\cdot(p-q),
 \qquad p,q\in\mathbb R^d.
\end{equation*}
Tightness and Skorokhod representation give a sequence $u^{\nu}$ and limit
$u$ on a new probability space.  We apply this inequality to their
gradients and integrate by parts.  The curvature bound
\eqref{intro-viscous-area-estimate}, strong convergence of $ u^{\nu}$
in $L^2((0,T)\times\mathbb T^d)$, and weak convergence of its gradients
make the integrated right-hand side tend to zero.  Hence, almost surely,
\begin{equation*}
 \frac{\nabla u^{\nu}}{H(\nabla u^{\nu})}
 \longrightarrow
 \frac{\nabla u}{H(\nabla u)}
 \quad\text{strongly in }L^2((0,T)\times\mathbb T^d;\mathbb R^d).
\end{equation*}
See
Proposition~\ref{lem-strong-normal-limit} for more details.

For the middle factor, the uniform bound
\eqref{intro-viscous-dirichlet-estimate} gives a weak limit $ Y_D$,
\begin{equation*}
 \frac{D^2 u^{\nu}}{H(\nabla u^{\nu})}
 \rightharpoonup Y_D
 \quad\text{weakly in }
 L^2((0,T)\times\mathbb T^d;\mathbb R^{d\times d}),
 \quad\text{almost surely}.
\end{equation*}
To identify $ Y_D$ as the weighted Hessian of $ u$, we use the
Sobolev chain rule for the approximations,
\begin{equation*}
 \nabla\left(
 \frac{\nabla u^{\nu}}{H(\nabla u^{\nu})}
 \right)
 =
 \left(I-\frac{\nabla u^{\nu}\otimes
 \nabla u^{\nu}}{H(\nabla u^{\nu})^2}\right)
 \frac{D^2 u^{\nu}}{H(\nabla u^{\nu})}.
\end{equation*}
We observe that the first matrix-valued factor on the right hand side is bounded, and the weighted-Hessian estimate give
uniform control of $\frac{\nabla u^{\nu}}{H(\nabla u^{\nu})}$ in
$L^2(0,T;H^1(\mathbb T^d;\mathbb R^d))$.
The aforementioned strong convergence of $\frac{\nabla u^{\nu}}{H(\nabla u^{\nu})}$ identifies the weak derivative limit as
$\nabla(\nabla u/H(\nabla u))$.  The matrix coefficients also converge
strongly and are uniformly bounded, so for every
$\Phi\in L^2((0,T)\times\mathbb T^d;\mathbb R^{d\times d})$,
\begin{equation*}
 \left(I-\frac{\nabla u^{\nu}\otimes
 \nabla u^{\nu}}{H(\nabla u^{\nu})^2}\right)\Phi
 \longrightarrow
 \left(I-\frac{\nabla u\otimes\nabla u}
 {H(\nabla u)^2}\right)\Phi
 \quad\text{strongly in }
 L^2((0,T)\times\mathbb T^d;\mathbb R^{d\times d}),
\end{equation*}
almost surely.  Pairing the weak and strong convergences, we then obtain
\begin{align*}
 \int_0^T\!\int_{\mathbb T^d}
 \nabla\left(\frac{\nabla u}{H(\nabla u)}\right)
 :\Phi\,dx\,dt&=\lim_{\nu\to0}\int_0^T\!\int_{\mathbb T^d}
 \nabla\left(
 \frac{\nabla u^{\nu}}{H(\nabla u^{\nu})}
 \right):\Phi\,dx\,dt\\
 &=\lim_{\nu\to0}\int_0^T\!\int_{\mathbb T^d}
 \left[\left(I-\frac{\nabla u^{\nu}\otimes
 \nabla u^{\nu}}{H(\nabla u^{\nu})^2}\right)
 \frac{D^2 u^{\nu}}{H(\nabla u^{\nu})}
 \right]:\Phi\,dx\,dt\\
 &=\lim_{\nu\to0}\int_0^T\!\int_{\mathbb T^d}
 \frac{D^2 u^{\nu}}{H(\nabla u^{\nu})}:
 \left[\left(I-\frac{\nabla u^{\nu}\otimes
 \nabla u^{\nu}}{H(\nabla u^{\nu})^2}\right)
 \Phi\right]dx\,dt\\
 &=\int_0^T\!\int_{\mathbb T^d} Y_D:
 \left[\left(I-\frac{\nabla u\otimes\nabla u}
 {H(\nabla u)^2}\right)\Phi\right]dx\,dt\\
 &=\int_0^T\!\int_{\mathbb T^d}
 \left[\left(I-\frac{\nabla u\otimes\nabla u}
 {H(\nabla u)^2}\right) Y_D\right]:
 \Phi\,dx\,dt.
\end{align*}
Since $\Phi$ is arbitrary,
\begin{equation*}
 \nabla\left(\frac{\nabla u}{H(\nabla u)}\right)
 =
 \left(I-\frac{\nabla u\otimes\nabla u}
 {H(\nabla u)^2}\right) Y_D.
\end{equation*}
Since $\nabla u$ is finite almost everywhere, the matrix is invertible
with inverse $I+\nabla u\otimes\nabla u$.  We thus define and identify
the weighted Hessian by
\begin{equation*}
 \frac{D^2 u}{H(\nabla u)}
 :=
 \left(I-\frac{\nabla u\otimes\nabla u}
 {H(\nabla u)^2}\right)^{-1}
 \nabla\left(\frac{\nabla u}{H(\nabla u)}\right)
 = Y_D.
\end{equation*}
This agrees with the classical expression for smooth functions and with
Definition~\ref{def-energy-martingale-solution}.
We emphasize that unweighted $H^2(\mathbb T^d)$ regularity of the limit $u$ is not needed, and the inverse need not be uniformly bounded as well, but the identified field is in
$L^2( \Omega\times(0,T)\times\mathbb T^d;
\mathbb R^{d\times d})$ by the uniform weighted-Hessian estimate.

This identification resolves the main technical difficulty and is the
main new contribution of this paper. Existence then follows by
relatively standard martingale arguments and passage to the limit.

%%%%%%%%
\subsection{Comments on the literature}

We collect related existing literature in this part. 

\textbf{Physical background and sharp-interface limits.}
Mean-curvature flow describes local interface relaxation under isotropic
surface tension.  The Allen--Cahn model~\cite{AllenCahn1979} uses a
diffuse phase field with a thin transition layer to describe this motion.
Evans, Soner, and Souganidis~\cite{EvansSonerSouganidis1992} prove
convergence to generalized mean-curvature motion beyond the formation
of singularities, while Ilmanen~\cite{Ilmanen1993} obtains convergence
to Brakke's mean-curvature flow.  These sharp-interface limits justify
the deterministic drift in \eqref{SPDE-0}.  Kawasaki and
Ohta~\cite{KawasakiOhta1982} include thermal fluctuations in the
interface motion.  Katsoulakis and Kho~\cite{KatsoulakisKho2001}
formally derive the corresponding stochastic graph equation from a
phase-field model.  This motivates the noise amplitude discussed in
Subsection~\ref{subsec-thermal-background}.

For stochastic sharp-interface limits, the results depend on the noise
structure and scaling.  Funaki~\cite{Funaki1999} derives randomly
perturbed curvature motion for planar interfaces from a stochastic
reaction--diffusion equation.  Weber~\cite{Weber2010Dynamics} extends
this short-time result to arbitrary ambient dimensions for spatially
homogeneous, temporally regularized additive noise whose integral
converges to Brownian motion.  For Allen--Cahn equations with spatially
correlated stochastic transport noise, R\"oger and
Weber~\cite{RogerWeber2013} prove uniform diffuse-area estimates and
tightness, with phase-indicator functions as sharp-interface limits.
For one-dimensional interface fluctuations and equilibrium
sharp-interface limits, see \cite{Funaki1995,Weber2010Invariant},
respectively.

\textbf{Hamilton--Jacobi equations and viscosity solutions.}
In the level-set formulation, prescribed normal motion gives a
first-order Hamilton--Jacobi equation, while curvature dependence
introduces a degenerate second-order term.  Viscosity solutions give
comparison, existence, and stability for these equations; see
Crandall, Ishii, and Lions~\cite{CrandallIshiiLions1992}.
Evans and Spruck~\cite{EvansSpruck1991} and Chen, Giga, and
Goto~\cite{ChenGigaGoto1991} develop this theory for mean-curvature
motion.  For suitable entire graphs, Ecker and
Huisken~\cite{EckerHuisken1989} prove global smooth evolution and
gradient estimates.  Our equation has the same curvature operator,
and, for smooth driving signals, its gradient-dependent forcing is
a first-order Hamiltonian term.  Lions and
Souganidis~\cite{LionsSouganidis1998a,LionsSouganidis1998b} initiated
the pathwise viscosity theory for fully nonlinear equations with rough
time dependence, including gradient-dependent stochastic forcing and
applications to front propagation.

\textbf{Existence and uniqueness for stochastic interface equations.}
Dirr, Luckhaus, and Novaga~\cite{DirrLuckhausNovaga2001} study a
forced planar curvature flow with fattening in the deterministic
setting.  In the limit of vanishing stochastic forcing, the law
selects the upper and lower barrier evolutions.  For spatially
constant normal noise, Souganidis and
Yip~\cite{SouganidisYip2004} prove non-fattening results for particular
planar curves.  These results concern geometric selection and
uniqueness under spatially coherent normal forcing, which differs
from the local thermal forcing in \eqref{SPDE-0}.  This geometric
uniqueness is distinct from uniqueness for weak solutions of a graph SPDE.

In one graph dimension, Es-Sarhir and
von Renesse~\cite{EsSarhirVonRenesse2012} prove variational
existence and uniqueness for stochastic curve-shortening equations
under suitable assumptions on the forcing.  For spatially homogeneous
normal noise, Gess and R\"ockner~\cite{GessRockner2017} prove
existence and uniqueness in the sense of stochastic variational
inequalities when the noise intensity is at most $\sqrt2$ in their
normalization.  Dareiotis and Gess~\cite[Section~7]{DareiotisGess2020} treat
spatially inhomogeneous normal noise by differentiating the graph
equation.  The slope solves a scalar conservative stochastic
diffusion equation with nonlinearity $\arctan$, for which they prove
existence and uniqueness of entropy solutions under the stated
coefficient and moment assumptions.  The reduction is related to
the one-dimensional observation above, though their height-noise
amplitude is $(1+|\partial_xu|^2)^{1/2}$ rather than the thermal
amplitude $(1+|\partial_xu|^2)^{1/4}$ used here.

In higher graph dimensions, Hofmanov\'a, R\"oger, and
von Renesse~\cite{HofmanovaRogerVonRenesse2017} construct weak
martingale solutions in dimension two using energy estimates and
stochastic compactness.  Dabrock, Hofmanov\'a, and
R\"oger~\cite[Theorem~3.5]{DabrockHofmanovaRoger2021} obtain
martingale solutions that are strong in the PDE sense in arbitrary
graph dimension.  They assume a deterministic uniform bound on the
initial gradient and prove that it is preserved by the evolution.

Gassiat, Gess, Lions, and
Souganidis~\cite{GassiatGessLionsSouganidis2024} study the long-time
behavior of pathwise solutions to stochastic Hamilton--Jacobi
equations, including stochastic mean-curvature equations.  They
prove decay of spatial oscillations and, in particular settings,
accelerated relaxation under spatially homogeneous noise.  Although
their appendix discusses the formal thermal-noise graph equation
underlying our model, the mean-curvature decay results concern
homogeneous normal noise rather than the coloured thermal forcing
studied here.
For long-time behavior under additive forcing and regularization by
nonlinear noise, see \cite{DirrSouganidis2005,GassiatGess2019}.
For geometric constructions and numerical approximations of
stochastic curvature flows, see
\cite{FengLiProhl2014,Yip1998}.

\subsection{Structure of the paper}

The structure of this paper is organized as follows.  In
Section~\ref{sec-ito-formulation}, we introduce the It\^o equation
and define the solution concept.
Section~\ref{sec-geometric-estimates} is devoted to the geometric and
algebraic identities used in the energy estimates.  In
Section~\ref{sec-viscous-construction}, we construct solutions of the
fixed-viscosity equation by Galerkin approximations. We also establish path regularity,
and pass to the limit in the Dirichlet inequality.  In
Section~\ref{sec-vanishing-viscosity}, we derive estimates uniform
in viscosity, prove tightness, and identify
the geometric coordinates and martingale covariations to complete the
proof of the main theorem.

%%%%%%%%%%%%%%%%%
\section{The It\^o formulation and solution concept}
\label{sec-ito-formulation}

We fix integers $d\geq2$, $T\in(0,\infty)$, $\varepsilon>0$, and $K\in\mathbb N^+$ throughout the paper. The negative Sobolev index used for path spaces is fixed as
\begin{equation*}
 m:=3+\frac d2, 
\end{equation*}
so that the embeddings $H^m(\mathbb T^d)\hookrightarrow C^2(\mathbb T^d)$ and
$L^1(\mathbb T^d)\hookrightarrow H^{-m}(\mathbb T^d)$ hold continuously.
We impose the following assumptions on the covariance and dimension
to obtain the cancellations and coercivity in the energy estimates.

\begin{assumption}
\label{constant-noise-covariance}
The real-valued modes $f_1,\ldots,f_K\in C^\infty(\mathbb T^d)$ satisfy
\begin{equation}
F_1(x):=\sum_{k=1}^Kf_k(x)^2\equiv\Lambda_K,\qquad x\in\mathbb T^d,
\end{equation}
for a finite constant $\Lambda_K\geq0$.
We set $\lambda:=\varepsilon\Lambda_K/8$ and assume
\begin{equation}\label{dimension-noise-assumption}
 \begin{cases}
 d=2,\quad \lambda\geq0\text{ arbitrary},\\
 d\geq3,\quad 0\leq\lambda<\dfrac{d+1}{2(d-1)}.
 \end{cases}
\end{equation}
Throughout the paper, constants in inequalities may change from line to line.
Unless otherwise stated, they may depend on the dimension, the torus,
$T$, $\varepsilon$, the fixed modes, and on $(1-2\lambda(d-1)/(d+1))^{-1}$ when $d\geq3$.
We write $a\lesssim b$ when $a\leq Cb$ for a constant $C$ with this dependence.
\end{assumption}

We rewrite the equation \eqref{SPDE-0} into the following It\^o form
\begin{equation}\label{SPDE-Ito}
 du=\left[H(\nabla u)\nabla\cdot
 \left(\frac{\nabla u}{H(\nabla u)}\right)
 +\frac{\varepsilon\Lambda_K}{8}
 \frac{\nabla u\cdot D^2u\,\nabla u}{H(\nabla u)^3}\right]dt+\sqrt{\varepsilon}\,H(\nabla u)^{1/2}
 \sum_{k=1}^Kf_k\,d\beta_t^k.
\end{equation}
Equivalently, 
\begin{equation}\label{SPDE-Ito-matrix}
 du=\left[
 \left(I-\frac{\nabla u\otimes\nabla u}{H(\nabla u)^2}\right):D^2u
 +\frac{\varepsilon\Lambda_K}{8}
 \frac{\nabla u\otimes\nabla u}{H(\nabla u)^3}:D^2u
\right]dt+\sqrt{\varepsilon}\,H(\nabla u)^{1/2}
 \sum_{k=1}^Kf_k\,d\beta_t^k,
\end{equation}
where $A:B:=\sum_{i,j=1}^d A_{ij}B_{ij}$ is the Frobenius contraction.
For any matrix $A$, $|A|:=(A:A)^{1/2}$ denotes the Frobenius norm. We define all vector
norms throughout the paper by the Euclidean norm.
We regard vectors as columns and use the Jacobian convention
$(\nabla z)_{ij}=\partial_jz_i$ for vector fields $z$.

We next provide a calculation concerning the Stratonovich-to-It\^o correction terms. The following calculation is performed for smooth functions under
Assumption~\ref{constant-noise-covariance}, which implies
$\sum_{k=1}^K f_k\nabla f_k=\tfrac12\nabla F_1=0$. We write
\begin{equation}\label{galerkin-noise-definition}
 G_k(u):=f_kH(\nabla u)^{1/2}.
\end{equation}
For smooth $u$ and a smooth variation $\eta$, we have 
\begin{align*}
 DG_k(u)[\eta]
 &=\frac{f_k}{2H(\nabla u)^{3/2}}\nabla u\cdot\nabla\eta,\\
 \frac{\varepsilon}{2}\sum_{k=1}^K DG_k(u)[G_k(u)]
 &=\frac{\varepsilon}{4H(\nabla u)}
   \nabla u\cdot\sum_{k=1}^K f_k\nabla f_k
   +\frac{\varepsilon F_1}{8}
   \frac{\nabla u\cdot D^2u\,\nabla u}{H(\nabla u)^3}=\lambda\frac{\nabla u\cdot D^2u\,\nabla u}{H(\nabla u)^3}.
\end{align*}

By the Stratonovich-to-It\^o conversion formula, the Stratonovich
equation~\eqref{SPDE-0} formally becomes
\begin{equation*}
 du=\left[H(\nabla u)\nabla\cdot
 \left(\frac{\nabla u}{H(\nabla u)}\right)
 +\frac{\varepsilon}{2}\sum_{k=1}^K DG_k(u)[G_k(u)]\right]dt
 +\sqrt{\varepsilon}\sum_{k=1}^K G_k(u)\,d\beta_t^k.
\end{equation*}
Substituting the above correction and the definition of $G_k$, we obtain
\eqref{SPDE-Ito}.

We interpret the second-order terms in \eqref{SPDE-Ito} and
\eqref{SPDE-Ito-matrix} through the weighted Hessian in the following
definition, and we remark that no unweighted $H^2(\mathbb T^d)$ regularity is needed.

\begin{definition}[Energy martingale solution]
\label{def-energy-martingale-solution}
Let $\Lambda$ be a Borel probability measure on $H^1(\mathbb T^d)$ with
$\int_{H^1(\mathbb T^d)}\|z\|_{H^1(\mathbb T^d)}^2\,\Lambda(dz)<\infty$.
An energy martingale solution to \eqref{SPDE-Ito} with initial law
$\Lambda$ consists of a stochastic basis
$(\Omega,\mathcal F,(\mathcal F_t)_{t\in[0,T]},\mathbb P)$ satisfying
the usual conditions, independent real $(\mathcal F_t)$-Brownian motions
$\beta^1,\ldots,\beta^K$, and an
$(\mathcal F_t)$-adapted, $H^{-m}(\mathbb T^d)$-valued process $u$,
where $m=3+d/2$, such that
\begin{align*}
 &u\in L^2(\Omega\times(0,T);H^1(\mathbb T^d)),\\
 &\frac{\nabla u}{H(\nabla u)}
 \in L^2(\Omega\times(0,T);H^1(\mathbb T^d;\mathbb R^d)),\\
 &\frac{D^2u}{H(\nabla u)},\quad H(\nabla u)^{1/2}\kappa(u)
 \in L^2(\Omega\times(0,T)\times\mathbb T^d).
\end{align*}
Here the weighted Hessian and the It\^o correction are defined by
\begin{equation}\label{weak-weighted-hessian-definition}
 \frac{D^2u}{H(\nabla u)}
 :=\left(I-\frac{\nabla u\otimes\nabla u}{H(\nabla u)^2}\right)^{-1}
 \nabla\left(\frac{\nabla u}{H(\nabla u)}\right)
\end{equation}
and
\begin{equation*}
 \frac{\nabla u\cdot D^2u\,\nabla u}{H(\nabla u)^3}
 :=\frac{\nabla u}{H(\nabla u)}\cdot
 \frac{D^2u}{H(\nabla u)}
 \frac{\nabla u}{H(\nabla u)},
\end{equation*}
respectively.
The process $u$ has almost surely continuous $H^{-m}(\mathbb T^d)$-valued paths,
and $u(0)$ is $\mathcal F_0$-measurable with values in $H^1(\mathbb T^d)$
and law $\Lambda$.
On a common probability-one event, the following identity holds for every
$\varphi\in C^\infty(\mathbb T^d)$ and $t\in[0,T]$:
\begin{equation}\label{weak-martingale-formulation}
  \begin{aligned}
 \langle u(t),\varphi\rangle
 ={}&\langle u(0),\varphi\rangle
 +\int_0^t\int_{\mathbb T^d}
 H(\nabla u)\kappa(u)\varphi\,dx\,ds+\frac{\varepsilon\Lambda_K}{8}\int_0^t\int_{\mathbb T^d}
 \frac{\nabla u\cdot D^2u\,\nabla u}
 {H(\nabla u)^3}\varphi\,dx\,ds\\
 &+\sqrt{\varepsilon}\sum_{k=1}^K\int_0^t
 \left\langle f_kH(\nabla u)^{1/2},\varphi\right\rangle d\beta_s^k.
\end{aligned}
\end{equation}
\end{definition}

\begin{remark}
In Definition~\ref{def-energy-martingale-solution},
$\kappa(u)=\operatorname{div}(\nabla u/H(\nabla u))$ is understood
distributionally.  For every $p\in\mathbb R^d$, the matrix
$I-p\otimes p/H(p)^2$ is positive definite, with inverse
\begin{equation*}
 \left(I-\frac{p\otimes p}{H(p)^2}\right)^{-1}=I+p\otimes p.
\end{equation*}
Since $\nabla u$ is finite almost everywhere,
\eqref{weak-weighted-hessian-definition} is well defined almost everywhere
and equivalently reads
\begin{equation*}
 \frac{D^2u}{H(\nabla u)}
 :=\left(I+\nabla u\otimes\nabla u\right)
 \nabla\left(\frac{\nabla u}{H(\nabla u)}\right).
\end{equation*}
Noting that the inverse need not be uniformly bounded, and the $L^2$ condition
on the weighted Hessian is a separate requirement, not a consequence
of the normalized gradient belonging to $H^1(\mathbb{T}^d)$.
For every $u\in H^2(\mathbb T^d)$, it follows from the Sobolev chain rule that
\begin{equation*}
 \nabla\left(\frac{\nabla u}{H(\nabla u)}\right)
 =\left(I-\frac{\nabla u\otimes\nabla u}{H(\nabla u)^2}\right)
 \frac{D^2u}{H(\nabla u)}.
\end{equation*}
Consequently, the weighted Hessian and It\^o correction in
Definition~\ref{def-energy-martingale-solution} agree with their
classical expressions for smooth $u$ and their almost-everywhere
expressions for $u\in H^2(\mathbb T^d)$.
\end{remark}

\begin{remark}
Using $H(\nabla u)\leq1+|\nabla u|$ and the assumed $H^1(\mathbb{T}^d)$ regularity, we have
$\mathbb E\int_0^T\mathcal A(u(t))\,dt<\infty$.
The mean-curvature drift is integrable because
\begin{align*}
 \mathbb E\int_0^T\!\int_{\mathbb T^d}
 |H(\nabla u)\kappa(u)|\,dx\,dt\leq
 \left(\mathbb E\int_0^T\mathcal A(u(t))\,dt\right)^{1/2}
 \left(\mathbb E\int_0^T\!\int_{\mathbb T^d}
 H(\nabla u)|\kappa(u)|^2\,dx\,dt\right)^{1/2}<\infty.
\end{align*}
The It\^o correction satisfies
\begin{equation*}
 \left|\frac{\nabla u\cdot D^2u\,\nabla u}{H(\nabla u)^3}\right|
 \leq\frac{|D^2u|}{H(\nabla u)}
 \in L^2(\Omega\times(0,T)\times\mathbb T^d).
\end{equation*}
These terms are therefore integrable and well-defined in \eqref{weak-martingale-formulation}.
\end{remark}

\begin{remark}
We take predictable versions of the $H^1(\mathbb{T}^d)$-representative of $u$,
the normalized gradient, the weighted fields, and $G_k(u)$ in their
respective Sobolev or spatial $L^2$ spaces, up to $d\mathbb P\,dt$-null sets.
By Assumption~\ref{constant-noise-covariance}, we have
\begin{equation*}
 \sum_{k=1}^K\|G_k(u)\|_{L^2(\mathbb T^d)}^2
 =\Lambda_K\mathcal A(u).
\end{equation*}
Hence, for every smooth $\varphi$,
\begin{equation*}
 \mathbb E\int_0^T\sum_{k=1}^K
 |\langle G_k(u),\varphi\rangle|^2\,dt
 \leq\Lambda_K\|\varphi\|_{L^2(\mathbb T^d)}^2
 \mathbb E\int_0^T\mathcal A(u(t))\,dt<\infty,
\end{equation*}
so the stochastic integrals are well defined.
\end{remark}

%%%%%%%%%%%%%%%%%
\section{Geometric identities and spatial energy estimates}
\label{sec-geometric-estimates}

In this section, we present elementary chain rule calculations for matrix-valued processes and derive estimates for the terms appearing in the Dirichlet energy
$\frac12\|\nabla u\|_{L^2(\mathbb T^d)}^2$
and graph area
$\mathcal A(u)$
estimates. 

Under Assumption~\ref{constant-noise-covariance}, for fixed $\nu\geq0$,
we consider the viscous approximation equation 
\begin{equation*}
\begin{aligned}
 du={}&\left[\nu\Delta u
 +H(\nabla u)\nabla\cdot\left(\frac{\nabla u}{H(\nabla u)}\right)
 +\lambda\frac{\nabla u\cdot D^2u\,\nabla u}{H(\nabla u)^3}
\right]dt+\sqrt{\varepsilon}\sum_{k=1}^K
 f_kH(\nabla u)^{1/2}\,d\beta_t^k.
\end{aligned}
\end{equation*}

We provide the following informal estimates. 

\noindent\textit{i)~Formal $H^1(\mathbb T^d)$ estimate.}
It\^o's formula for the Dirichlet energy and periodic integration by
parts give
\begin{align*}
 d\left(\frac12\|\nabla u\|_{L^2(\mathbb T^d)}^2\right)
 =&\left[-\nu\|\Delta u\|_{L^2(\mathbb T^d)}^2
 -\int_{\mathbb T^d}\Delta u\,H(\nabla u)\kappa(u)\,dx\right.\\
 &\left.-\lambda\int_{\mathbb T^d}\Delta u\,
 \frac{\nabla u\cdot D^2u\,\nabla u}{H(\nabla u)^3}\,dx
 +\frac{\varepsilon}{2}\sum_{k=1}^K
 \|\nabla G_k(u)\|_{L^2(\mathbb T^d)}^2\right]dt\\
 &-\sqrt\varepsilon\sum_{k=1}^K
 \left(\int_{\mathbb T^d}\Delta u\,G_k(u)\,dx\right)d\beta_t^k.
\end{align*}
Since
\begin{equation*}
 \nabla G_k(u)=H(\nabla u)^{1/2}\nabla f_k
 +\frac{f_k}{2H(\nabla u)^{3/2}}D^2u\,\nabla u,
 \qquad k=1,\ldots,K,
\end{equation*}
the identities $\sum_{k=1}^Kf_k^2=\Lambda_K$ and
$\sum_{k=1}^Kf_k\nabla f_k=0$ imply
\begin{equation*}
 \frac{\varepsilon}{2}\sum_{k=1}^K\|\nabla G_k(u)\|_{L^2(\mathbb T^d)}^2
 =\frac{\varepsilon}{2}\int_{\mathbb T^d}
 H(\nabla u)\sum_{k=1}^K|\nabla f_k|^2\,dx
 +\lambda\int_{\mathbb T^d}
 \frac{|D^2u\,\nabla u|^2}{H(\nabla u)^3}\,dx.
\end{equation*}
Consequently,
\begin{align*}
 &d\left(\frac12\|\nabla u\|_{L^2(\mathbb T^d)}^2\right)
 +\nu\|\Delta u\|_{L^2(\mathbb T^d)}^2\,dt\\
 =&\left[-\int_{\mathbb T^d}\Delta u\,H(\nabla u)\kappa(u)\,dx
 +\lambda\int_{\mathbb T^d}
 \frac{|D^2u\,\nabla u|^2
 -(\Delta u)(\nabla u\cdot D^2u\,\nabla u)}
 {H(\nabla u)^3}\,dx\right.\\
 &\left.
 +\frac{\varepsilon}{2}\int_{\mathbb T^d}
 H(\nabla u)\sum_{k=1}^K|\nabla f_k|^2\,dx\right]dt
 -\sqrt\varepsilon\sum_{k=1}^K
 \left(\int_{\mathbb T^d}\Delta u\,G_k(u)\,dx\right)d\beta_t^k.
\end{align*}
We estimate these terms as follows.
\begin{enumerate}
 \item The curvature contribution
 $-\int_{\mathbb T^d}\Delta u\,H(\nabla u)\kappa(u)\,dx$
 is rewritten using Lemma~\ref{lem-dirichlet-geometric-identity}.
 \item The integral multiplied by $\lambda$ is identified and estimated
 in Lemma~\ref{lem-dirichlet-noise-identity}.  Combining it with the
 curvature contribution, Proposition~\ref{lem-spatial-dirichlet-coercivity}
 gives the weighted-Hessian dissipation, with no restriction on
 $\lambda$ in dimension two and with $\lambda<(d+1)/(2(d-1))$ for $d\geq3$.
 \item The remaining noise drift is bounded by $C\mathcal A(u)$. After localization, the stochastic
 term has zero expectation.  Applying Gronwall's inequality, we then obtain the
 expected $H^1(\mathbb T^d)$ estimate.  We justify this for the Galerkin solutions in
 Lemmas~\ref{lem-projected-dirichlet-energy} and
 \ref{lem-stopped-L2-estimate}.
\end{enumerate}

\medskip
\noindent\textit{ii)~Formal area estimate.}
Using $DH(\nabla u)=\nabla u/H(\nabla u)$ and
$\operatorname{div}(DH(\nabla u))=\kappa(u)$, It\^o's formula yields
\begin{align*}
 d\mathcal A(u)
 ={}&-\left[\nu\int_{\mathbb T^d}\kappa(u)\Delta u\,dx
 +\int_{\mathbb T^d}H(\nabla u)|\kappa(u)|^2\,dx\right]dt\\
 &+\left[-\lambda\int_{\mathbb T^d}\kappa(u)
 \frac{\nabla u\cdot D^2u\,\nabla u}{H(\nabla u)^3}\,dx
 +\frac{\varepsilon}{2}\sum_{k=1}^K\int_{\mathbb T^d}
 \nabla G_k(u)\cdot D^2H(\nabla u)\nabla G_k(u)\,dx\right]dt\\
 &-\sqrt\varepsilon\sum_{k=1}^K
 \left(\int_{\mathbb T^d}\kappa(u)G_k(u)\,dx\right)d\beta_t^k.
\end{align*}
We estimate these terms as follows.
\begin{enumerate}
 \item Lemma~\ref{lem-viscous-area-dissipation} shows that
 $\nu\int_{\mathbb T^d}\kappa(u)\Delta u\,dx\geq0$, so viscosity
 contributes additional dissipation.
 \item Lemma~\ref{lem-H2-area-correction} rewrites the two terms in the
 second bracket as
 \begin{equation*}
  \frac{\varepsilon}{2}\int_{\mathbb T^d}
  \left(\operatorname{tr}C_K-
  \frac{\nabla u\cdot C_K\nabla u}{H(\nabla u)^2}\right)dx
  +\lambda\int_{\mathbb T^d}\mathcal J_d(u)\,dx,
 \end{equation*}
 where $C_K$ and $\mathcal J_d$ are defined below in
 \eqref{noise-gradient-covariance} and \eqref{dimension-area-remainder}.
 Lemma~\ref{lem-area-spatial-coercivity} bounds the covariance term
 independently of $u$ and absorbs the $\mathcal J_d$ term into
 $\int_{\mathbb T^d}H(\nabla u)|\kappa(u)|^2\,dx$.
 \item For the stochastic term, the Cauchy--Schwarz inequality gives
 \begin{equation*}
  \varepsilon\sum_{k=1}^K
  \left|\int_{\mathbb T^d}\kappa(u)G_k(u)\,dx\right|^2
  \leq\varepsilon\left(\sum_{k=1}^K\|f_k\|_{L^2(\mathbb T^d)}^2\right)
  \int_{\mathbb T^d}H(\nabla u)|\kappa(u)|^2\,dx.
 \end{equation*}
 Taking the time supremum and applying the
 Burkholder--Davis--Gundy and Young inequalities gives the expected
 area supremum and weighted-curvature bound.  We carry out the argument
 in Proposition~\ref{prop-viscous-area-identity}, where we also justify
 It\^o's formula at $H^2(\mathbb T^d)$ regularity.
\end{enumerate}

%%%%%%%%%%%%%%%%
\subsection{Matrix identities and dimension-dependent Newton formulas}

For a sufficiently regular function $u:\mathbb T^d\to\mathbb R$, let
\begin{equation*}
 \Gamma_u:=\{(x,u(x)):x\in\mathbb T^d\}
 \subset\mathbb T^d\times\mathbb R
\end{equation*}
denote its graph.  We use the following notation for the
induced metric $g=(g_{ij})_{i,j=1}^d$, its inverse
$g^{-1}=(g^{ij})_{i,j=1}^d$, and the second fundamental form
$\mathcal B=(\mathcal B_{ij})_{i,j=1}^d$ of $\Gamma_u$ (with respect to the
upward unit normal):
\begin{equation*}
 g_{ij}=\delta_{ij}+\partial_i u\,\partial_j u,\qquad
 g^{ij}=\delta_{ij}-\frac{\partial_i u\,\partial_j u}{H(\nabla u)^2},
 \qquad
 \mathcal B_{ij}=\frac{\partial_{ij}u}{H(\nabla u)},\qquad
 i,j=1,\ldots,d.
\end{equation*}
We use these objects as notation for the estimates below.  We have
\begin{equation*}
 \kappa=\sum_{i=1}^d\sum_{j=1}^d g^{ij}\mathcal B_{ij},
 \qquad
 |\mathcal B|_g^2
=\sum_{i=1}^d\sum_{j=1}^d\sum_{k=1}^d\sum_{\ell=1}^d
 g^{ik}g^{j\ell}\mathcal B_{ij}\mathcal B_{k\ell}.
\end{equation*}
The following identities distinguish the two-dimensional case.

\begin{lemma}
Let $d\geq1$ and $u\in H^2(\mathbb T^d)$.  The following identities hold
almost everywhere:
\begin{align}
 \kappa(u)&=\frac{\Delta u}{H(\nabla u)}
 -\frac{\nabla u\cdot D^2u\,\nabla u}{H(\nabla u)^3},
 \label{matrix-curvature-formula}\\
 |\mathcal B(u)|_g^2
 &=\frac{|D^2u|^2}{H(\nabla u)^2}
 -2\frac{|D^2u\,\nabla u|^2}{H(\nabla u)^4}
 +\frac{|\nabla u\cdot D^2u\,\nabla u|^2}{H(\nabla u)^6}.
 \label{matrix-second-form-formula}
\end{align}
In particular,
\begin{equation}\label{matrix-curvature-bounds}
 |\kappa(u)|^2\leq d|\mathcal B(u)|_g^2,
 \qquad
 |\mathcal B(u)|_g^2\leq\frac{|D^2u|^2}{H(\nabla u)^2}.
\end{equation}
When $d=2$, we define the cofactor of the Hessian by 
\begin{equation*}
 \operatorname{cof}D^2u
 :=\begin{pmatrix}
 \partial_{22}u&-\partial_{12}u\\
 -\partial_{12}u&\partial_{11}u
 \end{pmatrix}.
\end{equation*}
The following pointwise identities also hold:
\begin{equation}\label{matrix-two-dimensional-identities}
\begin{aligned}
 \operatorname{cof}D^2u&=(\Delta u)I-D^2u,\qquad
 (D^2u)^2-(\Delta u)D^2u+(\det D^2u)I=0,\\
 (\Delta u)^2-|D^2u|^2&=2\det D^2u,\qquad
 |D^2u\,\nabla u|^2-(\Delta u)(\nabla u\cdot D^2u\,\nabla u)
 =-|\nabla u|^2\det D^2u,\\
 |D^2u(\nabla u)^\perp|^2
 &=|\nabla u|^2|D^2u|^2-|D^2u\,\nabla u|^2.
\end{aligned}
\end{equation}
Here
\begin{equation*}
 (\nabla u)^\perp:=(-\partial_2u,\partial_1u)
\end{equation*}
denotes the counterclockwise rotation of $\nabla u$ by $\pi/2$.
For every $d\geq1$, periodic integration by parts yields
\begin{equation}\label{periodic-hessian-laplacian-equality}
 \|D^2u\|_{L^2(\mathbb T^d)}^2=\|\Delta u\|_{L^2(\mathbb T^d)}^2.
\end{equation}
\end{lemma}

\begin{proof}
Since $H(\nabla u)\geq1$, the coefficient matrix $g^{-1}$ has operator
norm at most one, and the Sobolev chain rule implies
$\kappa(u)\in L^2(\mathbb T^d)$, then all quantities in the statement are well defined.  Moreover, every weighted quadratic
Hessian density displayed above is bounded in absolute value by
$C_d|D^2u|^2\in L^1(\mathbb T^d)$, while the two norms in
\eqref{periodic-hessian-laplacian-equality} are finite because
$u\in H^2(\mathbb T^d)$.

Since $H(\nabla u)^2=1+|\nabla u|^2$, direct multiplication implies
\begin{equation*}
 (I+\nabla u\otimes\nabla u)
 \left(I-\frac{\nabla u\otimes\nabla u}{H(\nabla u)^2}\right)=I+\left(1-\frac{1+|\nabla u|^2}{H(\nabla u)^2}\right)
 \nabla u\otimes\nabla u=I.
\end{equation*}
Then we have 
\begin{equation*}
 g^{-1}=(I+\nabla u\otimes\nabla u)^{-1}
 =I-\frac{\nabla u\otimes\nabla u}{H(\nabla u)^2}.
\end{equation*}
By the definition of $\kappa$, 
\begin{equation*}
 \kappa(u)
=g^{-1}:\frac{D^2u}{H(\nabla u)}=\frac{I:D^2u}{H(\nabla u)}
 -\frac{(\nabla u\otimes\nabla u):D^2u}{H(\nabla u)^3}=\frac{\Delta u}{H(\nabla u)}
 -\frac{\nabla u\cdot D^2u\,\nabla u}{H(\nabla u)^3},
\end{equation*}
which is \eqref{matrix-curvature-formula}.  Expanding both inverse-metric
factors yields
\begin{align*}
 H(\nabla u)^2|\mathcal B(u)|_g^2
 &=\sum_{i=1}^d\sum_{j=1}^d\sum_{k=1}^d\sum_{\ell=1}^d
 \left(\delta_{ik}-\frac{\partial_i u\,\partial_k u}{H(\nabla u)^2}\right)
 \left(\delta_{j\ell}-\frac{\partial_j u\,\partial_\ell u}{H(\nabla u)^2}\right)
 \partial_{ij}u\,\partial_{k\ell}u\\
 &=|D^2u|^2-\frac{2|D^2u\,\nabla u|^2}{H(\nabla u)^2}
 +\frac{(\nabla u\cdot D^2u\,\nabla u)^2}{H(\nabla u)^4},
\end{align*}
which proves \eqref{matrix-second-form-formula}.  The matrix $g^{-1}$
is symmetric, with $d-1$ eigenvalues equal to $1$ and the remaining
eigenvalue equal to $H(\nabla u)^{-2}>0$ almost everywhere.
It is therefore positive definite almost everywhere.
Its symmetric positive square root $g^{-1/2}$ therefore exists, and
symmetry and cyclicity of the trace give
\begin{align*}
 \kappa(u)
 =\operatorname{tr}\left(
 \frac{g^{-1/2}D^2u\,g^{-1/2}}{H(\nabla u)}\right),\quad |\mathcal B(u)|_g^2
 =\frac{\operatorname{tr}(g^{-1}D^2u\,g^{-1}D^2u)}
 {H(\nabla u)^2}
 =\left|\frac{g^{-1/2}D^2u\,g^{-1/2}}
 {H(\nabla u)}\right|^2.
\end{align*}
Consequently, the Cauchy--Schwarz inequality and $|I|^2=d$ yield
\begin{align*}
 |\kappa(u)|^2
 =\left|I:\frac{g^{-1/2}D^2u\,g^{-1/2}}{H(\nabla u)}\right|^2\leq d\left|\frac{g^{-1/2}D^2u\,g^{-1/2}}
 {H(\nabla u)}\right|^2
 =d|\mathcal B(u)|_g^2.
\end{align*}
These are almost-everywhere matrix identities and require no smooth
principal-curvature representation.  Since
$(\nabla u\cdot D^2u\,\nabla u)^2
\leq|\nabla u|^2|D^2u\,\nabla u|^2$ and
$2H(\nabla u)^2-|\nabla u|^2=2+|\nabla u|^2\geq0$,
the expanded formula gives the second.

For $d=2$, the cofactor definition yields
$\operatorname{cof}D^2u=(\Delta u)I-D^2u$.
Multiplying this identity by $D^2u$ and using
$(\operatorname{cof}D^2u)D^2u=(\det D^2u)I$ proves the
Cayley--Hamilton identity in \eqref{matrix-two-dimensional-identities}.  Taking its
trace and its quadratic form against $\nabla u$ gives the next two
identities.  Finally, $\nabla u$ and $(\nabla u)^\perp$ are orthogonal
vectors of equal length, therefore 
$|D^2u\,\nabla u|^2+|D^2u(\nabla u)^\perp|^2
=|\nabla u|^2|D^2u|^2$.

For smooth periodic $u$ and every $i,j=1,\ldots,d$, applying the integrations by
parts formula twice, we can see that 
$\int_{\mathbb T^d}\partial_{ii}u\,\partial_{jj}u\,dx
=\int_{\mathbb T^d}|\partial_{ij}u|^2\,dx$.  Summing explicitly over
$i,j=1,\ldots,d$ proves
\eqref{periodic-hessian-laplacian-equality}, and density extends it to
$H^2(\mathbb T^d)$.  The pointwise formulas also hold at $H^2(\mathbb T^d)$ regularity by
the Sobolev chain rule.
\end{proof}

We use the next stability result to extend quadratic Hessian identities
from smooth graphs to $H^2(\mathbb T^d)$ graphs.

\begin{lemma}\label{lem-quadratic-hessian-continuity}
Let $d\geq1$.  Suppose $u_r\to u$ strongly in $H^2(\mathbb T^d)$ and
$a:\mathbb R^d\to\mathbb R$ is bounded and continuous.  Then, for
any $i,j,k,\ell=1,\ldots,d$,
\begin{equation}\label{quadratic-hessian-convergence}
 a(\nabla u_r)\,\partial_{ij}u_r\,\partial_{k\ell}u_r
 \longrightarrow
 a(\nabla u)\,\partial_{ij}u\,\partial_{k\ell}u
 \quad\text{in }L^1(\mathbb T^d).
\end{equation}
For every bounded continuous matrix map $A$, one also has
$A(\nabla u_r):D^2u_r\to A(\nabla u):D^2u$ in $L^2(\mathbb T^d)$.
\end{lemma}

\begin{proof}
Boundedness of $a$ and H\"older's inequality show that each product
in \eqref{quadratic-hessian-convergence} belongs to $L^1(\mathbb T^d)$,
whereas boundedness of $A$ and $D^2u_r,D^2u\in L^2(\mathbb T^d)$ show that
the two matrix contractions belong to $L^2(\mathbb T^d)$.

The Hessians converge strongly in $L^2(\mathbb T^d)$, and the gradients converge
strongly in $L^2(\mathbb T^d)$.  For the scalar assertion, write
\begin{align*}
 a(\nabla u_r)\partial_{ij}u_r\partial_{k\ell}u_r
 -a(\nabla u)\partial_{ij}u\partial_{k\ell}u&=a(\nabla u_r)(\partial_{ij}u_r-\partial_{ij}u)
   \partial_{k\ell}u_r
 +a(\nabla u_r)\partial_{ij}u
   (\partial_{k\ell}u_r-\partial_{k\ell}u)\\
 &\quad+[a(\nabla u_r)-a(\nabla u)]
   \partial_{ij}u\partial_{k\ell}u.
\end{align*}
The first two terms tend to zero in $L^1(\mathbb T^d)$ by the Cauchy--Schwarz inequality, the
boundedness of $a$, and the strong $L^2(\mathbb T^d)$ convergence of the Hessians.
For the last term, take an arbitrary subsequence.  It has a further
subsequence along which $\nabla u_r\to\nabla u$ almost everywhere.
Continuity of $a$ implies pointwise convergence to zero, while
\begin{equation*}
 |a(\nabla u_r)-a(\nabla u)|
 |\partial_{ij}u\,\partial_{k\ell}u|
 \leq2\|a\|_{L^\infty(\mathbb R^d)}
 |\partial_{ij}u\,\partial_{k\ell}u|\in L^1(\mathbb T^d).
\end{equation*}
By the dominated convergence theorem, we obtain $L^1(\mathbb T^d)$ convergence along this
further subsequence.  Since every subsequence has a further subsequence
converging to the same limit, the last term converges to zero in $L^1(\mathbb T^d)$
along the original sequence.  This proves
\eqref{quadratic-hessian-convergence}.

For a bounded continuous matrix map $A$, write
\begin{align*}
 A(\nabla u_r):D^2u_r-A(\nabla u):D^2u=A(\nabla u_r):(D^2u_r-D^2u)
 +[A(\nabla u_r)-A(\nabla u)]:D^2u.
\end{align*}
The first term tends to zero in $L^2(\mathbb T^d)$.  For the second, take an arbitrary
subsequence and then a further subsequence such that
$\nabla u_r\to\nabla u$ almost everywhere.  Since $A$ is bounded, we can see that 
\begin{equation*}
 |[A(\nabla u_r)-A(\nabla u)]:D^2u|^2
 \leq4C_A^2|D^2u|^2\in L^1(\mathbb T^d).
\end{equation*}
Continuity of $A$ and the dominated convergence theorem give convergence
to zero in $L^2(\mathbb T^d)$ along this further subsequence.  The same subsequence
criterion yields convergence of the original sequence and completes the
proof.
\end{proof}

The following identity gives the exact cancellations in the Dirichlet energy estimate and the area estimate in dimension two.

\begin{lemma}\label{lem-weighted-hessian-determinant}
For every $u\in H^2(\mathbb T^2)$ and
$\Phi\in C^1(\mathbb R^2;\mathbb R)\cap
L^\infty(\mathbb R^2;\mathbb R)$,
\begin{equation}\label{weighted-hessian-determinant-zero}
 \int_{\mathbb T^2}\Phi(\nabla u)\det D^2u\,dx=0.
\end{equation}
The same identity holds for $u\in C^2(\mathbb T^2)$ and arbitrary
$\Phi\in C^1(\mathbb R^2;\mathbb R)$.
\end{lemma}

\begin{proof}
In two dimensions,
$|\det D^2u|\leq |D^2u|^2/2\in L^1(\mathbb T^2)$, and $\Phi$ is bounded
in the $H^2(\mathbb T^2)$ case.  For $u\in C^2(\mathbb T^2)$ and
$\Phi\in C^1(\mathbb R^2;\mathbb R)$, the integrand
$\Phi(\nabla u)\det D^2u$ is continuous on $\mathbb T^2$ and hence integrable.
Thus the integral in
\eqref{weighted-hessian-determinant-zero} is finite.

For $u\in C^3(\mathbb T^2)$, define the vector field
$\mathcal F:\mathbb R^2\to\mathbb R^2$ by
\begin{equation*}
 \mathcal F(q):=\left(
 \int_0^{q_1}\Phi\bigl((s,q_2)\bigr)\,ds,0\right),
 \qquad q=(q_1,q_2)\in\mathbb R^2.
\end{equation*}
Then $\operatorname{div}_q\mathcal F(q)=\Phi(q)$.
Combining  
$\operatorname{div}(\operatorname{cof}D^2u)=0$ with 
$(\operatorname{cof}D^2u)D^2u=(\det D^2u)I$, we have 
\begin{equation*}
 \operatorname{div}\bigl[(\operatorname{cof}D^2u)\mathcal F(\nabla u)\bigr]
 =\Phi(\nabla u)\det D^2u.
\end{equation*}
Consequently, the torus integral vanishes.  For $u\in C^2(\mathbb T^2)$, periodic
mollifications converge in $C^2(\mathbb T^2)$, and the range of their gradients
remains bounded, so the identity passes to the limit for any $C^1(\mathbb R^2;\mathbb R)$
weight.  For $u\in H^2(\mathbb T^2)$, use periodic mollifications converging
strongly in $H^2(\mathbb T^2)$.  Since $\Phi$ is bounded,
Lemma~\ref{lem-quadratic-hessian-continuity} applies to both terms of
$\det D^2u=\partial_{11}u\,\partial_{22}u-(\partial_{12}u)^2$;
the integrands converge in $L^1(\mathbb T^2)$, proving the claim.
\end{proof}

The scalar-curvature divergence identity holds in every dimension.
In dimension two, the determinant identity gives an additional weighted
cancellation.

\begin{proposition}\label{lem-weighted-newton-identities}
Let $d\geq1$, $u\in C^3(\mathbb T^d)$, and
$\chi\in C^1([1,\infty))$.  Then
\begin{equation}\label{newton-divergence}
 \kappa(u)^2-|\mathcal B(u)|_g^2
 =\operatorname{div}\left(
 \frac{[(\Delta u)I-D^2u]\nabla u}{H(\nabla u)^2}\right),
\end{equation}
and
\begin{equation}\label{weighted-newton-identity}
 \int_{\mathbb T^d}\chi(H(\nabla u))
 \bigl(\kappa(u)^2-|\mathcal B(u)|_g^2\bigr)\,dx=\int_{\mathbb T^d}
 \frac{\chi'(H(\nabla u))}{H(\nabla u)^3}
 \bigl(|D^2u\,\nabla u|^2
 -\Delta u\,(\nabla u\cdot D^2u\,\nabla u)\bigr)\,dx.
\end{equation}
If $d=2$, one has in addition
\begin{equation}\label{two-dimensional-curvature-determinant}
 \kappa(u)^2-|\mathcal B(u)|_g^2
 =\frac{2\det D^2u}{H(\nabla u)^4},
\end{equation}
and, for every such $\chi$,
\begin{equation}\label{weighted-newton-two-dimensional}
 \int_{\mathbb T^2}\chi(H(\nabla u))
 \bigl(\kappa(u)^2-|\mathcal B(u)|_g^2\bigr)\,dx=2\int_{\mathbb T^2}
 \frac{\chi(H(\nabla u))}{H(\nabla u)^4}
 \det D^2u\,dx=0.
\end{equation}
For $u\in H^2(\mathbb T^d)$, the divergence identity
\eqref{newton-divergence} holds in distributions, and \eqref{two-dimensional-curvature-determinant} holds almost everywhere when $d=2$.
At this regularity, the integrated identity \eqref{weighted-newton-identity}
and, when $d=2$, \eqref{weighted-newton-two-dimensional} hold for
$\chi(r)=r$ and $\chi(r)=\log r$.
\end{proposition}

\begin{proof}
Commuting third derivatives, we obtain
\begin{equation*}
 \operatorname{div}[(\Delta u)I-D^2u]=0.
\end{equation*}
Consequently, it follows from
\eqref{matrix-curvature-formula}--\eqref{matrix-second-form-formula} that 
\begin{align*}
 \operatorname{div}\left(
 \frac{[(\Delta u)I-D^2u]\nabla u}{H(\nabla u)^2}\right)
 &=\frac{(\Delta u)^2-|D^2u|^2}{H(\nabla u)^2}
 +\frac{2\bigl(|D^2u\,\nabla u|^2
 -(\Delta u)(\nabla u\cdot D^2u\,\nabla u)\bigr)}
 {H(\nabla u)^4}\\
 &=\kappa(u)^2-|\mathcal B(u)|_g^2,
\end{align*}
Multiplying \eqref{newton-divergence} by $\chi(H(\nabla u))$,
integrating by parts, and using
$\nabla H(\nabla u)=D^2u\,\nabla u/H(\nabla u)$, we obtain
\begin{align*}
 \int_{\mathbb T^d}\chi(H(\nabla u))
 (\kappa(u)^2-|\mathcal B(u)|_g^2)\,dx
 &=-\int_{\mathbb T^d}
 \frac{\chi'(H(\nabla u))}{H(\nabla u)^3}
 (D^2u\,\nabla u)\cdot[(\Delta u)I-D^2u]\nabla u\,dx\\
 &=\int_{\mathbb T^d}
 \frac{\chi'(H(\nabla u))}{H(\nabla u)^3}
 \bigl(|D^2u\,\nabla u|^2
 -(\Delta u)(\nabla u\cdot D^2u\,\nabla u)\bigr)\,dx.
\end{align*}
For $d=2$, \eqref{matrix-two-dimensional-identities} reduces the
preceding pointwise expansion to
\begin{align*}
 \kappa(u)^2-|\mathcal B(u)|_g^2
 &=\frac{2\det D^2u}{H(\nabla u)^2}
 -\frac{2|\nabla u|^2\det D^2u}{H(\nabla u)^4}
 =\frac{2\det D^2u}{H(\nabla u)^4}.
\end{align*}
Lemma~\ref{lem-weighted-hessian-determinant} with
$\Phi(q)=\chi(H(q))/H(q)^4$ proves
\eqref{weighted-newton-two-dimensional}.

For $u\in H^2(\mathbb T^d)$, take smooth periodic mollifications converging in
$H^2(\mathbb T^d)$.  The vector field in \eqref{newton-divergence} converges in
$L^2(\mathbb T^d)$, since $q\mapsto q/H(q)^2$ is bounded and continuous. Furthermore, the curvature difference converges in $L^1(\mathbb T^d)$ by
Lemma~\ref{lem-quadratic-hessian-continuity}.  This proves the
distributional identity.  For $\chi(r)=r$ or $\chi(r)=\log r$,
the coefficients $\chi(H(q))/H(q)^2$ and
$\chi'(H(q))|q|^2/H(q)^3$ are bounded.  The same lemma therefore
passes all integrated identities to the limit, including the
two-dimensional zero.
\end{proof}

\begin{remark}
The divergence identity \eqref{newton-divergence} is the graph
scalar-curvature formula in \cite[Lemma~10]{Lam2010}, here written
on the periodic domain.  In dimension two, the cancellation
\eqref{weighted-newton-two-dimensional} with $\chi(r)=r$ is used in
\cite[proof of Proposition~5.1]{HofmanovaRogerVonRenesse2017}, where
it is explained by Gauss--Bonnet and also proved directly using a
divergence identity.  Here we use periodic integration by parts and
Lemma~\ref{lem-weighted-hessian-determinant} to obtain the weighted
identities without invoking Gauss--Bonnet.  This also gives the
cancellations for both $\chi(r)=r$ and $\chi(r)=\log r$ at
$H^2(\mathbb T^2)$ regularity, as needed below.
\end{remark}

\subsection{Spatial terms in the Dirichlet estimate}

We estimate the curvature part and noise correction separately, then
combine them to obtain coercivity.  In dimension two, the cancellation
is exact.  In dimensions $d\geq3$, coercivity follows under the small-noise
condition.

\begin{lemma}\label{lem-dirichlet-geometric-identity}
For $d\geq1$ and $u\in H^2(\mathbb T^d)$,
\begin{equation}\label{newton-dirichlet-identity}
  \begin{aligned}
 \int_{\mathbb T^d}\Delta u\,H(\nabla u)\kappa(u)\,dx=\int_{\mathbb T^d}\left[
 \frac{H(\nabla u)^2}{2}
 \bigl(\kappa(u)^2+|\mathcal B(u)|_g^2\bigr)
 +\frac{|D^2u\,\nabla u|^2}{H(\nabla u)^2}
 -\frac{|\nabla u\cdot D^2u\,\nabla u|^2}{H(\nabla u)^4}
 \right]dx.
  \end{aligned}
\end{equation}
For $d=2$, this also has the exact form
\begin{equation}\label{dirichlet-two-geometric-identity}
 \int_{\mathbb T^2}\Delta u\,H(\nabla u)\kappa(u)\,dx
 =\int_{\mathbb T^2}
 \frac{|D^2u|^2+|D^2u(\nabla u)^\perp|^2}
 {H(\nabla u)^2}\,dx.
\end{equation}
\end{lemma}

\begin{proof}
Since 
\begin{equation*}
 |H(\nabla u)\kappa(u)|\leq\sqrt d\,|D^2u|,
 \qquad
 H(\nabla u)^2|\mathcal B(u)|_g^2\leq |D^2u|^2,
\end{equation*}
then the left-hand side of \eqref{newton-dirichlet-identity} is
integrable by the Cauchy--Schwarz inequality.  The remaining densities are integrable
because
\begin{equation*}
 \frac{|D^2u\,\nabla u|^2}{H(\nabla u)^2}
 +\frac{|\nabla u\cdot D^2u\,\nabla u|^2}{H(\nabla u)^4}
 \leq 2|D^2u|^2,
\end{equation*}
and when $d=2$,
$|D^2u(\nabla u)^\perp|^2/H(\nabla u)^2\leq|D^2u|^2$.

By using equations \eqref{matrix-curvature-formula} and
\eqref{matrix-second-form-formula}, we have the pointwise identity
\begin{align*}
 &\frac{H(\nabla u)^2}{2}
 (\kappa(u)^2+|\mathcal B(u)|_g^2)
 +\frac{|D^2u\,\nabla u|^2}{H(\nabla u)^2}
 -\frac{(\nabla u\cdot D^2u\,\nabla u)^2}{H(\nabla u)^4}\\
 =&\frac{(\Delta u)^2+|D^2u|^2}{2}
 -\frac{(\Delta u)(\nabla u\cdot D^2u\,\nabla u)}{H(\nabla u)^2}.
\end{align*}
Noting that \eqref{periodic-hessian-laplacian-equality} turns
the first summand into $\int_{\mathbb T^d}(\Delta u)^2dx$, while
\begin{equation*}
 H(\nabla u)\kappa(u)
 =\Delta u-\frac{\nabla u\cdot D^2u\,\nabla u}{H(\nabla u)^2}.
\end{equation*}
This proves
\eqref{newton-dirichlet-identity} directly at $H^2(\mathbb T^d)$ regularity.

For $d=2$, it follows from \eqref{matrix-two-dimensional-identities} that
\begin{align*}
 \Delta u\,H(\nabla u)\kappa(u)
 &=|D^2u|^2+2\det D^2u
 -\frac{|D^2u\,\nabla u|^2
 +|\nabla u|^2\det D^2u}{H(\nabla u)^2}\\
 &=\frac{|D^2u|^2+|D^2u(\nabla u)^\perp|^2}
 {H(\nabla u)^2}
 +\left(1+\frac1{H(\nabla u)^2}\right)\det D^2u.
\end{align*}
The final term has zero torus integral directly at $H^2(\mathbb T^d)$ regularity
by Lemma~\ref{lem-weighted-hessian-determinant}, applied with the
bounded weight $\Phi(q)=1+H(q)^{-2}$. This completes the proof. 
\end{proof}

Under Assumption~\ref{constant-noise-covariance}, the noise term can
be written as a Newton curvature difference, whose integral vanishes
in dimension two.

\begin{lemma}\label{lem-dirichlet-noise-identity}
Under Assumption~\ref{constant-noise-covariance}, let
$u\in H^2(\mathbb T^d)$.  Then
\begin{equation}\label{dirichlet-noise-identity}
 \int_{\mathbb T^d}
 \frac{|D^2u\,\nabla u|^2
 -\Delta u\,(\nabla u\cdot D^2u\,\nabla u)}
 {H(\nabla u)^3}\,dx=\int_{\mathbb T^d}H(\nabla u)
 \bigl(\kappa(u)^2-|\mathcal B(u)|_g^2\bigr)\,dx.
\end{equation}
If $d=2$, then both sides equal zero.  If $d\geq3$, then 
\begin{equation}\label{dirichlet-noise-general-bound}
 \int_{\mathbb T^d}H(\nabla u)
 (\kappa(u)^2-|\mathcal B(u)|_g^2)\,dx
 \leq\frac{d-1}{d+1}\int_{\mathbb T^d}H(\nabla u)
 (\kappa(u)^2+|\mathcal B(u)|_g^2)\,dx.
\end{equation}
\end{lemma}

\begin{proof}
For $u\in H^2(\mathbb T^d)$, we have
\begin{equation*}
 \frac{|D^2u\,\nabla u|^2}{H(\nabla u)^3}
 +\frac{|\Delta u|\,|\nabla u\cdot D^2u\,\nabla u|}
 {H(\nabla u)^3}
 \leq (1+\sqrt d)|D^2u|^2.
\end{equation*}
Furthermore,
$|\kappa(u)|\leq\sqrt d\,|D^2u|/H(\nabla u)$ and
$|\mathcal B(u)|_g^2\leq|D^2u|^2/H(\nabla u)^2$, so both
$H(\nabla u)\kappa(u)^2$ and
$H(\nabla u)|\mathcal B(u)|_g^2$ belong to $L^1(\mathbb T^d)$.
All integrals above are therefore well defined.

Use Proposition~\ref{lem-weighted-newton-identities} with $\chi(r)=r$.
For $d=2$, \eqref{weighted-newton-two-dimensional} gives the
integral zero.
For $d\geq3$, \eqref{matrix-curvature-bounds} implies
$\kappa(u)^2\leq d|\mathcal B(u)|_g^2$, and hence
\begin{equation*}
 \kappa^2-|\mathcal B|_g^2
 \leq\frac{d-1}{d+1}(\kappa^2+|\mathcal B|_g^2). 
\end{equation*}
This proves \eqref{dirichlet-noise-general-bound}.
\end{proof}

Combining the two identities yields the dissipation used in the
Galerkin estimate.

\begin{proposition}\label{lem-spatial-dirichlet-coercivity}
Under Assumption~\ref{constant-noise-covariance}, let
$u\in H^2(\mathbb T^d)$, and set
\begin{equation}\label{dirichlet-coercivity-constant}
 c_{d,\lambda}:=
 \begin{cases}
 1,&d=2,\\
 \dfrac12-\dfrac{d-1}{d+1}\lambda,&d\geq3.
 \end{cases}
 \qquad c_{d,\lambda}>0.
\end{equation}
If $d=2$, then
\begin{equation}\label{spatial-dirichlet-two-identity}
\begin{aligned}
 &-\int_{\mathbb T^2}\Delta u\,H(\nabla u)\kappa(u)\,dx
 +\lambda\int_{\mathbb T^2}
 \frac{|D^2u\,\nabla u|^2
 -\Delta u\,(\nabla u\cdot D^2u\,\nabla u)}
 {H(\nabla u)^3}\,dx\\
 =&-\int_{\mathbb T^2}
 \frac{|D^2u|^2+|D^2u(\nabla u)^\perp|^2}
 {H(\nabla u)^2}\,dx.
\end{aligned}
\end{equation}
If $d\geq3$, then
\begin{equation}\label{spatial-dirichlet-general-bound}
  \begin{aligned}
 &-\int_{\mathbb T^d}\Delta u\,H(\nabla u)\kappa(u)\,dx
 +\lambda\int_{\mathbb T^d}
 \frac{|D^2u\,\nabla u|^2
 -\Delta u\,(\nabla u\cdot D^2u\,\nabla u)}
 {H(\nabla u)^3}\,dx\\
 \leq&-c_{d,\lambda}
 \int_{\mathbb T^d}
 \left(H(\nabla u)^2\kappa(u)^2
 +\frac{|D^2u|^2}{H(\nabla u)^2}\right)dx.
\end{aligned}
\end{equation}
\end{proposition}

\begin{proof}
The two-dimensional assertion is
\eqref{dirichlet-two-geometric-identity} together with the zero in
Lemma~\ref{lem-dirichlet-noise-identity}.
For $d\geq3$, the Cauchy--Schwarz inequality and
$H(\nabla u)^2=1+|\nabla u|^2$ give
\begin{align*}
 &\frac{|D^2u\,\nabla u|^2}{H(\nabla u)^2}
 -\frac{(\nabla u\cdot D^2u\,\nabla u)^2}{H(\nabla u)^4}
 \geq\frac{|D^2u\,\nabla u|^2}{H(\nabla u)^4}\geq0.
\end{align*}
Moreover, \eqref{matrix-second-form-formula} yields
\begin{align*}
 &H(\nabla u)^2|\mathcal B(u)|_g^2
 +\frac{|D^2u\,\nabla u|^2}{H(\nabla u)^2}
 -\frac{(\nabla u\cdot D^2u\,\nabla u)^2}{H(\nabla u)^4}=|D^2u|^2-\frac{|D^2u\,\nabla u|^2}{H(\nabla u)^2}
 \geq\frac{|D^2u|^2}{H(\nabla u)^2}.
\end{align*}
By Lemmas~\ref{lem-dirichlet-geometric-identity} and
\ref{lem-dirichlet-noise-identity}, the left side of
\eqref{spatial-dirichlet-general-bound} is the negative integral of
the first expression below.  By \eqref{dirichlet-noise-general-bound}
and its pointwise proof,
$\kappa(u)^2-|\mathcal B(u)|_g^2
\leq\frac{d-1}{d+1}(\kappa(u)^2+|\mathcal B(u)|_g^2)$. Since
$H(\nabla u)\geq1$ and $0<c_{d,\lambda}\leq1/2$, the preceding two
inequalities give, pointwise,
\begin{align*}
 &\frac{H(\nabla u)^2}{2}
 (\kappa(u)^2+|\mathcal B(u)|_g^2)
 -\lambda H(\nabla u)(\kappa(u)^2-|\mathcal B(u)|_g^2)+\frac{|D^2u\,\nabla u|^2}{H(\nabla u)^2}
 -\frac{(\nabla u\cdot D^2u\,\nabla u)^2}{H(\nabla u)^4}\\
 \geq& c_{d,\lambda}H(\nabla u)^2
 (\kappa(u)^2+|\mathcal B(u)|_g^2)+\frac{|D^2u\,\nabla u|^2}{H(\nabla u)^2}
 -\frac{(\nabla u\cdot D^2u\,\nabla u)^2}{H(\nabla u)^4}\\
 \geq& c_{d,\lambda}\Biggl[
 H(\nabla u)^2\kappa(u)^2
 +H(\nabla u)^2|\mathcal B(u)|_g^2+\frac{|D^2u\,\nabla u|^2}{H(\nabla u)^2}
 -\frac{(\nabla u\cdot D^2u\,\nabla u)^2}{H(\nabla u)^4}
 \Biggr]\\
 \geq& c_{d,\lambda}\left(
 H(\nabla u)^2\kappa(u)^2
 +\frac{|D^2u|^2}{H(\nabla u)^2}\right).
\end{align*}
The first inequality uses the pointwise noise bound above and
$H(\nabla u)\geq1$, and the second retains a fraction
$c_{d,\lambda}\leq1/2$ of the nonnegative Hessian difference.
Integrating and changing the sign proves
\eqref{spatial-dirichlet-general-bound}.
\end{proof}

\subsection{Spatial terms in the graph-area estimate}
\label{subsec-area-spatial}

Under Assumption~\ref{constant-noise-covariance}, we have
\begin{equation}\label{noise-gradient-covariance}
 F_2:=\sum_{k=1}^K f_k\nabla f_k=\tfrac12\nabla F_1=0,
 \qquad
 C_K:=\sum_{k=1}^K\nabla f_k\otimes\nabla f_k.
\end{equation}
For every $u\in H^2(\mathbb T^d)$ define
\begin{equation}\label{dimension-area-remainder}
 \mathcal J_d(u):=
 \begin{cases}
 0,&d=2,\\
 \log H(\nabla u)
 \bigl(\kappa(u)^2-|\mathcal B(u)|_g^2\bigr),&d\geq3.
 \end{cases}
\end{equation}
This remainder is integrable by
\eqref{matrix-curvature-bounds} and
$\sup_{r\geq1}\log r/r^2<\infty$.
The first lemma gives the dissipative sign of the viscous term in
every dimension.

\begin{lemma}\label{lem-viscous-area-dissipation}
Under Assumption~\ref{constant-noise-covariance}, let
$u\in H^2(\mathbb T^d)$.  Then
\begin{equation}\label{viscosity-area-compatibility}
  \begin{aligned}
 \int_{\mathbb T^d}\kappa(u)\Delta u\,dx
 &=\sum_{j=1}^d\int_{\mathbb T^d}
 D^2H(\nabla u)
 [\partial_j\nabla u,\partial_j\nabla u]\,dx\\
 &=\int_{\mathbb T^d}\left(
 \frac{|D^2u|^2}{H(\nabla u)}
 -\frac{|D^2u\,\nabla u|^2}{H(\nabla u)^3}
 \right)dx\geq\int_{\mathbb T^d}
 \frac{|D^2u|^2}{H(\nabla u)^3}\,dx\geq0.
\end{aligned}
\end{equation}
\end{lemma}

\begin{proof}
Since $\kappa(u),\Delta u\in L^2(\mathbb T^d)$, their product is integrable.
Moreover, $0\leq D^2H(q)\leq I$, and each of the remaining quadratic
densities is bounded by $|D^2u|^2\in L^1(\mathbb T^d)$.
Thus all integrals in \eqref{viscosity-area-compatibility} are finite.

For smooth $u$, two periodic integrations by parts give
\begin{align*}
 \int_{\mathbb T^d}\kappa(u)\Delta u\,dx
 &=-\int_{\mathbb T^d}\frac{\nabla u}{H(\nabla u)}
 \cdot\nabla\Delta u\,dx=\sum_{j=1}^d\int_{\mathbb T^d}
 \partial_j\!\left(\frac{\nabla u}{H(\nabla u)}\right)
 \cdot\partial_j\nabla u\,dx\\
 &=\sum_{j=1}^d\int_{\mathbb T^d} D^2H(\nabla u)
 [\partial_j\nabla u,\partial_j\nabla u]\,dx.
\end{align*}
The Sobolev chain rule and
$D^2H(q)=I/H(q)-q\otimes q/H(q)^3$ give the second equality.
For $u\in H^2(\mathbb T^d)$, periodic mollification, the $L^2(\mathbb T^d)$ assertion of
Lemma~\ref{lem-quadratic-hessian-continuity}, and its $L^1(\mathbb T^d)$
quadratic-form assertion justify the same identity.
Finally, $|D^2u\,\nabla u|^2
\leq|D^2u|^2|\nabla u|^2$ and
$H(\nabla u)^2-|\nabla u|^2=1$ give the lower bound.
\end{proof}

For the area It\^o correction, we use the logarithmic Newton weight
in place of the linear weight.  The remainder vanishes in dimension
two but may be nonzero in higher dimensions.  Recall that $G_k$,
$k=1,...,K$, are defined by \eqref{galerkin-noise-definition}, and
$C_K$ is defined by \eqref{noise-gradient-covariance}.  The
noise correction and quadratic variation terms in the area estimate
satisfy the following identity. 

\begin{lemma}\label{lem-H2-area-correction}
Under Assumption~\ref{constant-noise-covariance}, every
$u\in H^2(\mathbb T^d)$ satisfies
\begin{equation}\label{simplified-area-remainder}
  \begin{aligned}
 &\frac12\sum_{k=1}^K\int_{\mathbb T^d}
 \nabla G_k(u)\cdot D^2H(\nabla u)\nabla G_k(u)\,dx-\frac{\Lambda_K}{8}\int_{\mathbb T^d}\kappa(u)
 \frac{\nabla u\cdot D^2u\,\nabla u}{H(\nabla u)^3}\,dx\\
 =&\frac12\int_{\mathbb T^d}\left(
 \operatorname{tr}C_K-
 \frac{\nabla u\cdot C_K\nabla u}{H(\nabla u)^2}\right)dx
 +\frac{\Lambda_K}{8}\int_{\mathbb T^d}\mathcal J_d(u)\,dx.
\end{aligned}
\end{equation}
\end{lemma}

\begin{proof}
To verify the integrability, first note that
$H(\nabla u)\in L^1(\mathbb T^d)$ and
\begin{equation*}
 \nabla G_k(u)=H(\nabla u)^{1/2}\nabla f_k
 +\tfrac12 f_kH(\nabla u)^{-3/2}D^2u\,\nabla u\in L^2(\mathbb T^d).
\end{equation*}
Since $D^2H$ is bounded, the quadratic-variation density is integrable.
Also, $\kappa(u)\in L^2(\mathbb T^d)$ and
\begin{equation*}
 \left|\frac{\nabla u\cdot D^2u\,\nabla u}
 {H(\nabla u)^3}\right|\leq |D^2u|,
\end{equation*}
so the correction product is in $L^1(\mathbb T^d)$.  Finally, the covariance density
is bounded, while $\mathcal J_d(u)\in L^1(\mathbb T^d)$ follows from
\eqref{matrix-curvature-bounds} and
$\sup_{r\geq1}\log(r)/r^2<\infty$.

The Sobolev chain rule, valid for $u\in H^2(\mathbb T^d)$, yields
\begin{equation}\label{noise-spatial-gradient-formula}
 \nabla G_k(u)=H(\nabla u)^{1/2}\nabla f_k
 +\tfrac12 f_kH(\nabla u)^{-3/2}D^2u\,\nabla u.
\end{equation}
Both summands belong to $L^2(\mathbb T^d)$, since
$H(\nabla u)\in L^1(\mathbb T^d)$ and
$|\nabla u|/H(\nabla u)^{3/2}\leq1$.
Using \eqref{noise-gradient-covariance} and
$D^2H(q)=H(q)^{-1}(I-q\otimes q/H(q)^2)$, direct expansion gives,
\begin{align*}
 \frac12\sum_{k=1}^K\nabla G_k(u)\cdot
 D^2H(\nabla u)\nabla G_k(u)
 &=\frac12\left(\operatorname{tr}C_K-
 \frac{\nabla u\cdot C_K\nabla u}{H(\nabla u)^2}\right)\\
 &+\frac{\Lambda_K}{8}
 \left(\frac{|D^2u\,\nabla u|^2}{H(\nabla u)^4}
 -\frac{(\nabla u\cdot D^2u\,\nabla u)^2}
 {H(\nabla u)^6}\right).
\end{align*}
Equation \eqref{matrix-curvature-formula} then yields
\begin{align*}
 &\frac{|D^2u\,\nabla u|^2}{H(\nabla u)^4}
 -\frac{(\nabla u\cdot D^2u\,\nabla u)^2}{H(\nabla u)^6}
 -\kappa(u)\frac{\nabla u\cdot D^2u\,\nabla u}
 {H(\nabla u)^3}=\frac{|D^2u\,\nabla u|^2
 -(\Delta u)(\nabla u\cdot D^2u\,\nabla u)}
 {H(\nabla u)^4}.
\end{align*}
For $d=2$, it has zero integral by
Proposition~\ref{lem-weighted-newton-identities} with $\chi(r)=\log r$.
For $d\geq3$, the same proposition with this weight identifies its integral as
$\int\log H(\nabla u)(\kappa(u)^2-|\mathcal B(u)|_g^2)$.
These $H^2(\mathbb T^d)$ identities prove \eqref{simplified-area-remainder}.
The displayed bounds, smoothness of the modes, and
\eqref{matrix-curvature-bounds} show all integrands are in $L^1(\mathbb T^d)$.
\end{proof}

The last spatial estimate absorbs the area remainder into curvature
dissipation for $d\geq3$.  In dimension two, the remainder is zero.

\begin{samepage}
\begin{lemma}\label{lem-area-spatial-coercivity}
Under Assumption~\ref{constant-noise-covariance}, let
$u\in H^2(\mathbb T^d)$.  Then
\begin{equation}\label{area-spatial-coercivity}
 H(\nabla u)\kappa(u)^2-\lambda\mathcal J_d(u)
 \geq a_{d,\lambda}H(\nabla u)\kappa(u)^2
 \quad\text{a.e.},\quad
 a_{d,\lambda}:=
 \begin{cases}
 1,&d=2,\\
 1-\lambda(d-1)/(de),&d\geq3,
 \end{cases}
 \quad a_{d,\lambda}>0.
\end{equation}
Moreover,
\begin{equation}\label{area-covariance-remainder-bound}
 0\leq\operatorname{tr}C_K-
 \frac{\nabla u\cdot C_K\nabla u}{H(\nabla u)^2}
 \leq\operatorname{tr}C_K.
\end{equation}
\end{lemma}
\end{samepage}

\begin{proof}
By \eqref{matrix-curvature-bounds} and the boundedness
of $\log(r)/r^2$ on $[1,\infty)$, we have
\begin{equation*}
 H(\nabla u)\kappa(u)^2\leq d|D^2u|^2\in L^1(\mathbb T^d),
 \qquad \mathcal J_d(u)\in L^1(\mathbb T^d).
\end{equation*}
The quantities in the two inequalities are therefore well defined at
$H^2(\mathbb T^d)$ regularity.
The covariance expression
is bounded by the smooth function $\operatorname{tr}C_K$ and is therefore
integrable on $\mathbb T^d$.

For $d=2$, $\mathcal J_d(u)=0$.  For $d\geq3$,
\eqref{matrix-curvature-bounds}, namely
$\kappa(u)^2\leq d|\mathcal B(u)|_g^2$, gives
\begin{equation*}
 \mathcal J_d(u)
 \leq\frac{d-1}{d}\log H(\nabla u)\,\kappa(u)^2
 \leq\frac{d-1}{de}H(\nabla u)\kappa(u)^2.
\end{equation*}
This proves \eqref{area-spatial-coercivity}.  Under the assumed
small-noise condition,
$a_{d,\lambda}>1-(d+1)/(2de)>0$ for $d\geq3$.
Finally, $C_K$ is a sum of nonnegative
rank-one matrices, so
$0\leq\nabla u\cdot C_K\nabla u
\leq |\nabla u|^2\operatorname{tr}C_K
\leq H(\nabla u)^2\operatorname{tr}C_K$, then this proves
\eqref{area-covariance-remainder-bound}.
\end{proof}

\begin{remark}
Since $(d+1)/(2(d-1))>1/2$ for every $d\geq3$, the condition
$\lambda<1/2$ is sufficient for Theorem~\ref{thm-target-existence}
in every fixed dimension $d\geq3$.  For $d=3$, the dimension-dependent
condition is $\lambda<1$.  The absence of a small-noise restriction
in dimension two relies on the exact determinant cancellation.
\end{remark}

%%%%%%%%%%%%%%%%%%%%% 
\section{Construction of the viscous equation}
\label{sec-viscous-construction}

For $\nu>0$, we consider the viscous approximation equation: 
\begin{equation}\label{viscous-equation}
 du^\nu=\left[\nu\Delta u^\nu
 +H(\nabla u^\nu)\kappa(u^\nu)
 +\lambda\frac{\nabla u^\nu\cdot D^2u^\nu\nabla u^\nu}
 {H(\nabla u^\nu)^3}\right]dt+\sqrt{\varepsilon}\sum_{k=1}^K
 f_kH(\nabla u^\nu)^{1/2}\,d\beta_t^k.
\end{equation}

The definition of martingale solutions to \eqref{viscous-equation} is similar to Definition~\ref{def-energy-martingale-solution}, requiring an additional 
regularity condition 
$u^\nu\in L^2(\Omega\times(0,T);H^2(\mathbb T^d))$
and with an extra drift term
\begin{equation*}
 \nu\int_0^t\langle\Delta u^\nu(s),\varphi\rangle\,ds
 =-\nu\int_0^t\int_{\mathbb T^d}
 \nabla u^\nu(s)\cdot\nabla\varphi\,dx\,ds
\end{equation*}
on the right-hand side of \eqref{weak-martingale-formulation}.

Recall the noise tensors $F_2=0$ and $C_K$ from
\eqref{noise-gradient-covariance}.  We first state the existence result
and then construct the solution in the following subsections.
Unless otherwise specified, the convergence arrows $\longrightarrow$
and $\rightharpoonup$ for the Galerkin approximations in this section
refer to $N\to\infty$ at fixed $\nu>0$.

\begin{theorem}\label{thm-viscous-equation}
Suppose Assumption~\ref{constant-noise-covariance} holds.  Fix an
$\mathcal F_0$-measurable $H^1(\mathbb T^d)$-valued initial datum
$u_0$ with $\mathbb E\|u_0\|_{H^1(\mathbb T^d)}^2<\infty$.
For every $\nu>0$, there is a martingale solution
$u^\nu$ of \eqref{viscous-equation} with initial law
$\mathcal L(u_0)$.  Recall that $m=3+d/2$, the solution belongs to
$
 L^2\bigl(\Omega;L^2(0,T;H^2(\mathbb T^d))\bigr)
 \cap L^2\bigl(\Omega;C([0,T];H^{-m}(\mathbb T^d))\bigr),
$
and its trajectories belong to $C([0,T];H^1(\mathbb T^d))$ almost surely.
Moreover,
\begin{equation}\label{fixed-nu-solution-estimate}
\begin{aligned}
 &\sup_{0\leq t\leq T}\mathbb E\|\nabla u^\nu(t)\|_{L^2(\mathbb T^d)}^2
 +\nu\mathbb E\int_0^T\|\Delta u^\nu\|_{L^2(\mathbb T^d)}^2\,dt+\mathbb E\int_0^T\!\int_{\mathbb T^d}
 \frac{|D^2u^\nu|^2}{H(\nabla u^\nu)^2}\,dx\,dt\\
 &\leq C_T\bigl(1+\mathbb E\|u_0\|_{H^1(\mathbb T^d)}^2\bigr),
 \end{aligned}
\end{equation}
and
\begin{equation}\label{fixed-nu-solution-nonuniform-estimate}
 \begin{aligned}
 \mathbb E\sup_{0\leq t\leq T}\|u^\nu(t)\|_{L^2(\mathbb T^d)}^2
 +\mathbb E\int_0^T\|u^\nu(t)\|_{H^2(\mathbb T^d)}^2\,dt
 &\leq C_{\nu,T}\bigl(1+\mathbb E\|u_0\|_{H^1(\mathbb T^d)}^2\bigr).
 \end{aligned}
\end{equation}
Here $C_T$ is independent of $\nu$, whereas $C_{\nu,T}$ may depend on
$\nu$.
\end{theorem}

We fix the initial datum $u_0$ independently of $\nu$ and approximate
it by the Fourier projections $\Pi_Nu_0$ in the Galerkin construction.

%%%%%%%%%%%%%%%%%%%%
\subsection{Galerkin construction and uniform estimates}

We first derive estimates independent of the Galerkin dimension, with
constants which may depend on the fixed viscosity $\nu$.
Throughout this section, we fix the data of
Theorem~\ref{thm-viscous-equation} on the original stochastic basis.
For each $N\in\mathbb N^+$, let $\Pi_N$ be the real Fourier projection onto modes of size at most $N$,
and let $u_0^N=\Pi_Nu_0$. Then 
\begin{equation}\label{fixed-nu-initial-convergence}
 \mathbb E\|u_0^N-u_0\|_{H^1(\mathbb T^d)}^2\longrightarrow0.
\end{equation}
Consequently,
\begin{equation}\label{fixed-nu-initial-uniform-bound}
 C_0:=\sup_{N\geq1}\mathbb E\|u_0^N\|_{H^1(\mathbb T^d)}^2<\infty.
\end{equation}
The Galerkin equation of \eqref{viscous-equation} reads 
\begin{equation}\label{fully-regularized-equation}
  \begin{aligned}
 du^{N,\nu}=&\Pi_N\Big[
 \nu\Delta u^{N,\nu}
 +H(\nabla u^{N,\nu})\kappa(u^{N,\nu})
 +\lambda\frac{\nabla u^{N,\nu}\cdot D^2u^{N,\nu}
 \nabla u^{N,\nu}}{H(\nabla u^{N,\nu})^3}\Big]dt\\
 &+\sqrt{\varepsilon}\sum_{k=1}^K\Pi_N\left[
 f_kH(\nabla u^{N,\nu})^{1/2}\right]d\beta_t^k, 
\end{aligned}
\end{equation}
with initial condition $u^{N,\nu}(0)=u_0^N$. 
The projected equation is a finite-dimensional SDE, and we construct a
local strong solution and then apply the spatial coercivity from
Proposition~\ref{lem-spatial-dirichlet-coercivity} in the stopped
finite-dimensional It\^o formula.

\begin{lemma}\label{lem-projected-dirichlet-energy}
Suppose Assumption~\ref{constant-noise-covariance} holds.  Then for each
$N$ and $\nu>0$, equation \eqref{fully-regularized-equation} with
initial condition $u_0^N$ has a pathwise unique maximal local
strong solution $u^{N,\nu}$ on the given stochastic basis.  It is
adapted, has continuous $\Pi_N L^2(\mathbb T^d)$-valued paths up to its lifetime
$\zeta^{N,\nu}$, and satisfies the blow-up alternative
\begin{equation}\label{fixed-nu-blow-up-alternative}
 \{\zeta^{N,\nu}<\infty\}
 \subseteq
 \left\{\limsup_{t\uparrow\zeta^{N,\nu}}
 \|u^{N,\nu}(t)\|_{H^1(\mathbb T^d)}=\infty\right\}
 \quad\text{almost surely}.
\end{equation}
For every $R>0$, define the exit time
\begin{equation*}
 \tau_R:=\inf\bigl\{t\in[0,\zeta^{N,\nu}):
 \|u^{N,\nu}(t)\|_{H^1(\mathbb T^d)}\geq R\bigr\}\wedge T,
\end{equation*}
where $\inf\varnothing=\infty$.
Recall the constant $c_{d,\lambda}>0$ from
\eqref{dirichlet-coercivity-constant}.
Then, for every $t\in[0,T]$, almost surely,
\begin{equation}\label{projected-dirichlet-energy}
\begin{aligned}
 &\frac12\|\nabla u^{N,\nu}(t\wedge\tau_R)\|_{L^2(\mathbb T^d)}^2
 +\nu\int_0^{t\wedge\tau_R}\|\Delta u^{N,\nu}(s)\|_{L^2(\mathbb T^d)}^2\,ds+c_{d,\lambda}\int_0^{t\wedge\tau_R}\int_{\mathbb T^d}
 \frac{|D^2u^{N,\nu}|^2}{H(\nabla u^{N,\nu})^2}\,dx\,ds\\
 \leq& \frac12\|\nabla u^{N,\nu}(0)\|_{L^2(\mathbb T^d)}^2
 +\frac{\varepsilon}{2}\int_0^{t\wedge\tau_R}\int_{\mathbb T^d}
 H(\nabla u^{N,\nu})\operatorname{tr}C_K\,dx\,ds
 +\mathcal M_t^R,
\end{aligned}
\end{equation}
where
\begin{equation*}
 \mathcal M_t^R=-\sqrt\varepsilon\sum_{k=1}^K
 \int_0^{t\wedge\tau_R}
 \left(\int_{\mathbb T^d}\Delta u^{N,\nu}(s)\,
 f_kH(\nabla u^{N,\nu}(s))^{1/2}\,dx\right)d\beta_s^k
\end{equation*}
is a continuous square-integrable martingale with $\mathbb E\mathcal M_t^R=0$, for all $t\in[0,T]$.
\end{lemma}

\begin{proof}
Since $\Pi_N L^2(\mathbb T^d)$ is finite dimensional and $H(q)\geq1$, the projected drift
and diffusion coefficients in \eqref{fully-regularized-equation} are
smooth maps on $\Pi_N L^2(\mathbb T^d)$, hence locally Lipschitz.  The finite-dimensional
SDE theorem applied with the $\mathcal F_0$-measurable datum
$u_0^N$ yields an adapted, pathwise unique maximal local strong
solution with continuous paths and the blow-up alternative
\eqref{fixed-nu-blow-up-alternative}.  Continuity makes $\tau_R$ a
stopping time.  Since all norms on $\Pi_N L^2(\mathbb T^d)$ are equivalent,
$\tau_R<\zeta^{N,\nu}$ on $\{\zeta^{N,\nu}\leq T\}$, while
$\tau_R\leq T<\zeta^{N,\nu}$ on $\{\zeta^{N,\nu}>T\}$.
Thus every term below is defined on $[0,t\wedge\tau_R]$.
The identities $F_1=\Lambda_K$ and $F_2=0$ follow from
Assumption~\ref{constant-noise-covariance}, as recorded at the start
of this section.
Because $\Pi_N$ is the orthogonal Fourier projection, it is self-adjoint
and commutes with spatial derivatives.  For $z\in\Pi_N L^2(\mathbb T^d)$, we have
\begin{align*}
 \langle-\Delta z,\Pi_Nh\rangle_{L^2(\mathbb T^d)}
 &=\langle\Pi_N(-\Delta z),h\rangle_{L^2(\mathbb T^d)}
 =\langle-\Delta z,h\rangle_{L^2(\mathbb T^d)},
 &&h\in L^2(\mathbb T^d),\\
 \|\nabla\Pi_Nh\|_{L^2(\mathbb T^d)}^2
 &=\sum_{j=1}^d\|\Pi_N\partial_jh\|_{L^2(\mathbb T^d)}^2
 \leq\sum_{j=1}^d\|\partial_jh\|_{L^2(\mathbb T^d)}^2
 =\|\nabla h\|_{L^2(\mathbb T^d)}^2,
 &&h\in H^1(\mathbb T^d).
\end{align*}

Applying the finite-dimensional It\^o formula on
$[0,t\wedge\tau_R]$, we have almost surely, 
\begin{align*}
 \frac12\|\nabla u^{N,\nu}(t\wedge\tau_R)\|_{L^2(\mathbb T^d)}^2
 &=\frac12\|\nabla u^{N,\nu}(0)\|_{L^2(\mathbb T^d)}^2
 -\nu\int_0^{t\wedge\tau_R}\|\Delta u^{N,\nu}\|_{L^2(\mathbb T^d)}^2\,ds\\
 &\quad-\int_0^{t\wedge\tau_R}\!\int_{\mathbb T^d}
 \Delta u^{N,\nu}H(\nabla u^{N,\nu})
 \kappa(u^{N,\nu})\,dxds\\
 &\quad-\frac{\varepsilon\Lambda_K}{8}
 \int_0^{t\wedge\tau_R}\!\int_{\mathbb T^d}
 \frac{\Delta u^{N,\nu}
 (\nabla u^{N,\nu}\cdot D^2u^{N,\nu}\nabla u^{N,\nu})}
 {H(\nabla u^{N,\nu})^3}\,dxds\\
 &\quad+\frac\varepsilon2\sum_{k=1}^K
 \int_0^{t\wedge\tau_R}
 \|\nabla\Pi_NG_k(u^{N,\nu})\|_{L^2(\mathbb T^d)}^2\,ds
 +\mathcal M_t^R\\
 &=:\frac12\|\nabla u^{N,\nu}(0)\|_{L^2(\mathbb T^d)}^2
 -\nu\int_0^{t\wedge\tau_R}\|\Delta u^{N,\nu}\|_{L^2(\mathbb T^d)}^2 ds
 -J_1-J_2+J_3+\mathcal M_t^R.
\end{align*}

For term $J_3$, 
by the spatial chain rule, we have
\begin{equation*}
 \nabla G_k(u^{N,\nu})=H(\nabla u^{N,\nu})^{1/2}\nabla f_k
 +\frac12f_kH(\nabla u^{N,\nu})^{-3/2}D^2u^{N,\nu}\,\nabla u^{N,\nu}.
\end{equation*}
Squaring this identity and summing over $k$ yields, pointwise on
$\mathbb T^d$,
\begin{align*}
 \sum_{k=1}^K|\nabla G_k(u^{N,\nu})|^2
 =&H(\nabla u^{N,\nu})\sum_{k=1}^K|\nabla f_k|^2+H(\nabla u^{N,\nu})^{-1}
 \left(\sum_{k=1}^Kf_k\nabla f_k\right)\cdot
 D^2u^{N,\nu}\nabla u^{N,\nu}\\
 &+\frac{\sum_{k=1}^Kf_k^2}
 {4}
 \frac{|D^2u^{N,\nu}\nabla u^{N,\nu}|^2}{H(\nabla u^{N,\nu})^3}\\
 =&H(\nabla u^{N,\nu})\operatorname{tr}C_K
 +\frac{\Lambda_K}{4}\frac{|D^2u^{N,\nu}\nabla u^{N,\nu}|^2}{H(\nabla u^{N,\nu})^3}.
\end{align*}
The last equality uses
$\sum_k|\nabla f_k|^2=\operatorname{tr}C_K$,
$\sum_kf_k\nabla f_k=F_2=0$, and
$\sum_kf_k^2=F_1=\Lambda_K$.  The gradient contraction of $\Pi_N$ therefore
gives
\begin{align*}
 J_3
 &\leq\frac{\varepsilon}{2}\int_0^{t\wedge\tau_R}\sum_{k=1}^K\|\nabla G_k(u^{N,\nu})\|_{L^2(\mathbb T^d)}^2 ds\\
 &=\frac{\varepsilon}{2}\int_0^{t\wedge\tau_R}\int_{\mathbb T^d}
 H(\nabla u^{N,\nu})\operatorname{tr}C_Kdxds
 +\lambda\int_0^{t\wedge\tau_R}\int_{\mathbb T^d}\frac{|D^2u^{N,\nu}\,\nabla u^{N,\nu}|^2}
 {H(\nabla u^{N,\nu})^3}dx ds\\
 &=:J_{3,1}+J_{3,2}.
\end{align*}

We next estimate for the term $-J_1-J_2+J_{3,2}$. 
For each $s<\tau_R$, $u^{N,\nu}(s)\in\Pi_N L^2(\mathbb T^d)\subset C^\infty(\mathbb T^d)$.
Apply the appropriate case of
Proposition~\ref{lem-spatial-dirichlet-coercivity} with
$\lambda=\varepsilon\Lambda_K/8$.  For $d\geq3$,
Assumption~\ref{constant-noise-covariance} implies
$\lambda<(d+1)/(2(d-1))$.
Dropping the nonnegative term involving
$|D^2u(\nabla u)^\perp|^2$ in the two-dimensional identity, both cases give
\begin{align*}
 &-\int_{\mathbb T^d}\Delta u^{N,\nu}
 H(\nabla u^{N,\nu})\kappa(u^{N,\nu})\,dx+\lambda\int_{\mathbb T^d}
 \frac{|D^2u^{N,\nu}\nabla u^{N,\nu}|^2
 -\Delta u^{N,\nu}(\nabla u^{N,\nu}\cdot
 D^2u^{N,\nu}\nabla u^{N,\nu})}
 {H(\nabla u^{N,\nu})^3}\,dx\\
 &\leq-c_{d,\lambda}\int_{\mathbb T^d}
 \frac{|D^2u^{N,\nu}|^2}{H(\nabla u^{N,\nu})^2}\,dx.
\end{align*}
Integrating the preceding spatial inequality over $[0,t\wedge\tau_R]$
therefore yields
\begin{equation*}
 -J_1-J_2+J_{3,2}
\leq
 -c_{d,\lambda}\int_0^{t\wedge\tau_R}\!\int_{\mathbb T^d}
 \frac{|D^2u^{N,\nu}|^2}{H(\nabla u^{N,\nu})^2}\,dx\,ds.
\end{equation*}
Substitution into the stopped It\^o identity yields the claimed inequality \eqref{projected-dirichlet-energy}.

In the following, we estimate the martingale term. For each $k$, the map
\begin{equation*}
 s\longmapsto-\sqrt\varepsilon\int_{\mathbb T^d}
 \Delta u^{N,\nu}(s)f_kH(\nabla u^{N,\nu}(s))^{1/2}\,dx
\end{equation*}
is adapted and continuous on $[0,\tau_R]$.  Its stopped version is therefore
progressively measurable and has a predictable version.  The quadratic
variation of the resulting local martingale is
\begin{align*}
 \langle\mathcal M^R\rangle_t
 =\varepsilon\sum_{k=1}^K\int_0^{t\wedge\tau_R}
 \left|\int_{\mathbb T^d}\Delta u^{N,\nu}(s)\,
 f_kH(\nabla u^{N,\nu}(s))^{1/2}\,dx\right|^2ds.
\end{align*}
By the Cauchy--Schwarz inequality and Assumption~\ref{constant-noise-covariance},
\begin{align*}
 \sum_{k=1}^K\left|\int_{\mathbb T^d}\Delta u^{N,\nu}\,
 f_kH(\nabla u^{N,\nu})^{1/2}\,dx\right|^2
 &\leq\|\Delta u^{N,\nu}\|_{L^2(\mathbb T^d)}^2
 \sum_{k=1}^K\|f_kH(\nabla u^{N,\nu})^{1/2}\|_{L^2(\mathbb T^d)}^2\\
 &=\|\Delta u^{N,\nu}\|_{L^2(\mathbb T^d)}^2
 \int_{\mathbb T^d}\left(\sum_{k=1}^Kf_k^2\right)
 H(\nabla u^{N,\nu})\,dx\\
 &=\Lambda_K\|\Delta u^{N,\nu}\|_{L^2(\mathbb T^d)}^2
 \mathcal A(u^{N,\nu}).
\end{align*}
Since $\Pi_N L^2(\mathbb T^d)$ is finite dimensional, there is $c_N<\infty$ such that
$\|\Delta z\|_{L^2(\mathbb T^d)}\leq c_N\|z\|_{H^1(\mathbb T^d)}$ for every $z\in\Pi_N L^2(\mathbb T^d)$.  On
$\{s<\tau_R\}$,
\begin{equation*}
 \|\Delta u^{N,\nu}(s)\|_{L^2(\mathbb T^d)}\leq c_NR.
\end{equation*}
Moreover,
\begin{align*}
 \mathcal A(u^{N,\nu}(s))
 &\leq\int_{\mathbb T^d}(1+|\nabla u^{N,\nu}(s)|)\,dx\lesssim 1+
 \|\nabla u^{N,\nu}(s)\|_{L^2(\mathbb T^d)}\lesssim1+R.
\end{align*}
Consequently,
\begin{align*}
 \mathbb E\langle\mathcal M^R\rangle_T
 &\lesssim\varepsilon\Lambda_KTc_N^2R^2(1+R)<\infty.
\end{align*}
Using It\^o's isometry, we can see that
\begin{equation*}
 \mathbb E|\mathcal M_t^R|^2
 =\mathbb E\langle\mathcal M^R\rangle_t
 \leq\mathbb E\langle\mathcal M^R\rangle_T<\infty.
\end{equation*}
Hence $\mathcal M^R$ is a continuous square-integrable martingale and
$\mathbb E\mathcal M_t^R=0$ for every $t\in[0,T]$.
\end{proof}

For every $z\in H^2(\mathbb T^d)$, we define the drift of \eqref{viscous-equation} by 
\begin{equation}\label{fixed-nu-drift-definition}
  \begin{aligned}
 b^\nu(z):=&\nu\Delta z
 +H(\nabla z)\kappa(z)
 +\lambda\frac{\nabla z\cdot D^2z\nabla z}
 {H(\nabla z)^3}\\
 =&\nu\Delta z
 +\left(I-\frac{\nabla z\otimes\nabla z}{H(\nabla z)^2}\right):D^2z
 +\lambda
 \frac{\nabla z\otimes\nabla z}{H(\nabla z)^3}:D^2z.
\end{aligned}
\end{equation}
In the following, we derive stopped
estimates for the solution and its gradient, and use them to exclude
finite-time explosion later in
Proposition~\ref{lem-global-galerkin-fixed-nu}.

\begin{lemma}\label{lem-stopped-L2-estimate}
Let $u^{N,\nu}$ be the unique maximal local strong solution from
Lemma~\ref{lem-projected-dirichlet-energy}, with lifetime
$\zeta^{N,\nu}$ and exit time $\tau_R$ as defined there.  Then
\begin{equation}\label{fixed-nu-stopped-L2-path-bound}
 \sup_{\substack{N\geq1, R>0}}\mathbb E\sup_{0\leq t\leq T}
 \|u^{N,\nu}(t\wedge\tau_R)\|_{L^2(\mathbb T^d)}^2
 \leq C_{\nu,T}(1+C_0),
\end{equation}
and
\begin{equation}\label{fixed-nu-stopped-gradient-bound}
 \sup_{\substack{N\geq1, R>0}}\sup_{0\leq t\leq T}
 \mathbb E\|\nabla u^{N,\nu}(t\wedge\tau_R)\|_{L^2(\mathbb T^d)}^2
 \leq C_T(1+C_0).
\end{equation}
Moreover,
\begin{equation}\label{fixed-nu-stopped-dissipation-bound}
 \sup_{\substack{N\geq1, R>0}}\Bigg\{
 \nu\mathbb E\int_0^{T\wedge\tau_R}
 \|\Delta u^{N,\nu}(s)\|_{L^2(\mathbb T^d)}^2\,ds+\mathbb E\int_0^{T\wedge\tau_R}\int_{\mathbb T^d}
 \frac{|D^2u^{N,\nu}|^2}{H(\nabla u^{N,\nu})^2}\,dx\,ds\Bigg\}
 \leq C_T(1+C_0),
\end{equation}
where $C_0$ is defined in \eqref{fixed-nu-initial-uniform-bound} and the
constants are independent of $N$ and $R$.
\end{lemma}

\begin{proof}
All constants below are independent of $N$, $R$, and $C_0$, for simplicity, their dependence
on the fixed structural parameters is suppressed.  The local solution,
its blow-up alternative, and $\tau_R$ are supplied by
Lemma~\ref{lem-projected-dirichlet-energy}.
Taking expectations in \eqref{projected-dirichlet-energy} is legitimate because its stopped
martingale is square integrable.  Dropping its nonnegative dissipation terms
and using
\begin{equation}\label{fixed-nu-area-gradient-bound}
 \mathcal A(z)\lesssim1+\|\nabla z\|_{L^2(\mathbb T^d)},
 \qquad z\in H^1(\mathbb T^d),
\end{equation}
we obtain
\begin{align*}
 \mathbb E\|\nabla u^{N,\nu}(t\wedge\tau_R)\|_{L^2(\mathbb T^d)}^2
 &\leq \mathbb E\|\nabla u_0^N\|_{L^2(\mathbb T^d)}^2
 +C\int_0^t\mathbb E\!\left[
 \mathbf 1_{\{s<\tau_R\}}
 \mathcal A(u^{N,\nu}(s))\right]ds\\
 &\leq C_0+C\int_0^t
 \left(1+\mathbb E\!\left[
 \mathbf 1_{\{s<\tau_R\}}
 \|\nabla u^{N,\nu}(s)\|_{L^2(\mathbb T^d)}\right]\right)ds\\
 &\leq C_0+C\int_0^t\left(
 1+\left(\mathbb E\!\left[
 \mathbf 1_{\{s<\tau_R\}}
 \|\nabla u^{N,\nu}(s)\|_{L^2(\mathbb T^d)}^2\right]\right)^{1/2}\right)ds\\
 &\leq C_0+C\int_0^t\left(
 1+\left(\mathbb E
 \|\nabla u^{N,\nu}(s\wedge\tau_R)\|_{L^2(\mathbb T^d)}^2\right)^{1/2}\right)ds\\
 &\leq C(1+C_0)+C\int_0^t
 \mathbb E\|\nabla u^{N,\nu}(s\wedge\tau_R)\|_{L^2(\mathbb T^d)}^2\,ds.
\end{align*}
Using Gronwall's lemma, we obtain \eqref{fixed-nu-stopped-gradient-bound}.
Returning to all nonnegative terms in
Lemma~\ref{lem-projected-dirichlet-energy}, taking $t=T$, and using
\eqref{fixed-nu-area-gradient-bound} and
\eqref{fixed-nu-stopped-gradient-bound}, we first note that
\begin{align*}
 \mathbb E\int_0^{T\wedge\tau_R}\mathcal A(u^{N,\nu}(s))\,ds
 &=\int_0^T\mathbb E\left[
 \mathbf 1_{\{s<\tau_R\}}\mathcal A(u^{N,\nu}(s))\right]ds\\
 &\lesssim1+
 \int_0^T\left(\mathbb E
 \|\nabla u^{N,\nu}(s\wedge\tau_R)\|_{L^2(\mathbb T^d)}^2\right)^{1/2}ds\\
 &\leq C_T(1+C_0).
\end{align*}
The expected stopped energy inequality \eqref{projected-dirichlet-energy}, with $c_{d,\lambda}^{-1}$
absorbed into $C_T$, therefore yields
\eqref{fixed-nu-stopped-dissipation-bound}.

The two symmetric matrix maps in \eqref{fixed-nu-drift-definition} satisfy
\begin{equation*}
 \left|I-\frac{q\otimes q}{H(q)^2}\right|\leq\sqrt d,
 \qquad
 \left|\frac{q\otimes q}{H(q)^3}\right|
 =\frac{|q|^2}{(1+|q|^2)^{3/2}}\leq1.
\end{equation*}
So we have
\begin{equation}\label{fixed-nu-drift-L2-bound}
 \|b^\nu(z)\|_{L^2(\mathbb T^d)}\leq C_{\nu,\varepsilon,\Lambda_K}
 \|D^2z\|_{L^2(\mathbb T^d)},
 \qquad z\in H^2(\mathbb T^d).
\end{equation}
Since $\Pi_N$ is an $L^2(\mathbb T^d)$ contraction,
\eqref{periodic-hessian-laplacian-equality} and
\eqref{fixed-nu-stopped-dissipation-bound} imply
\begin{equation}\label{fixed-nu-stopped-drift-bound}
 \sup_{\substack{N\geq1, R>0}}\mathbb E\int_0^{T\wedge\tau_R}
 \|\Pi_Nb^\nu(u^{N,\nu}(s))\|_{L^2(\mathbb T^d)}^2\,ds
 \leq C_{\nu,T}(1+C_0).
\end{equation}

For the diffusion, it follows from Assumption~\ref{constant-noise-covariance} that
\begin{equation}\label{galerkin-noise-HS-bound}
 \sum_{k=1}^K\|G_k(z)\|_{L^2(\mathbb T^d)}^2
 =\int_{\mathbb T^d}\left(\sum_{k=1}^Kf_k^2\right)
 H(\nabla z)\,dx
 =\Lambda_K\mathcal A(z),
 \qquad z\in H^1(\mathbb T^d).
\end{equation}
The $L^2(\mathbb T^d)$ contraction of $\Pi_N$, the stopped area estimate above, and
\eqref{galerkin-noise-HS-bound} therefore yield
\begin{equation}\label{fixed-nu-stopped-noise-bound}
 \sup_{\substack{N\geq1, R>0}}\mathbb E\int_0^{T\wedge\tau_R}
 \sum_{k=1}^K\|\Pi_NG_k(u^{N,\nu}(s))\|_{L^2(\mathbb T^d)}^2\,ds
 \leq C_T(1+C_0).
\end{equation}

Applying the finite-dimensional It\^o formula to
the quadratic function on $[0,t\wedge\tau_R]$, we have 
\begin{align}\label{fixed-nu-stopped-L2-Ito}
 \frac12\|u^{N,\nu}(t\wedge\tau_R)\|_{L^2(\mathbb T^d)}^2
 ={}&\frac12\|u_0^N\|_{L^2(\mathbb T^d)}^2
 +\int_0^{t\wedge\tau_R}
 \langle u^{N,\nu}(s),\Pi_Nb^\nu(u^{N,\nu}(s))\rangle_{L^2(\mathbb T^d)}\,ds\notag\\
 &+\frac\varepsilon2\int_0^{t\wedge\tau_R}
 \sum_{k=1}^K\|\Pi_NG_k(u^{N,\nu}(s))\|_{L^2(\mathbb T^d)}^2\,ds
 +\mathcal N_t^R,
\end{align}
where
\begin{equation*}
 \mathcal N_t^R:=\sqrt\varepsilon\sum_{k=1}^K
 \int_0^{t\wedge\tau_R}
 \langle u^{N,\nu}(s),\Pi_NG_k(u^{N,\nu}(s))\rangle_{L^2(\mathbb T^d)}
 \,d\beta_s^k.
\end{equation*}
On $\{s<\tau_R\}$, we have $\|u^{N,\nu}(s)\|_{L^2(\mathbb T^d)}\leq R$, and it follows from \eqref{fixed-nu-stopped-noise-bound} that  
\begin{align*}
 \mathbb E\langle\mathcal N^R\rangle_T
 &=\varepsilon\mathbb E\int_0^{T\wedge\tau_R}
 \sum_{k=1}^K
 |\langle u^{N,\nu}(s),\Pi_NG_k(u^{N,\nu}(s))\rangle_{L^2(\mathbb T^d)}|^2\,ds\\
 &\leq\varepsilon R^2\mathbb E\int_0^{T\wedge\tau_R}
 \sum_{k=1}^K\|\Pi_NG_k(u^{N,\nu}(s))\|_{L^2(\mathbb T^d)}^2\,ds<\infty. 
\end{align*}
Thus $\mathcal N^R$ is a
square-integrable martingale.  Moreover, the Burkholder--Davis--Gundy inequality and
Young inequality give, for every $t\in[0,T]$,
\begin{align*}
 \mathbb E\sup_{0\leq r\leq t}|\mathcal N_r^R|
 &\leq C\sqrt\varepsilon\,\mathbb E\left[
 \sup_{0\leq r\leq t}\|u^{N,\nu}(r\wedge\tau_R)\|_{L^2(\mathbb T^d)}
 \left(\int_0^{t\wedge\tau_R}\sum_{k=1}^K
 \|\Pi_NG_k(u^{N,\nu}(s))\|_{L^2(\mathbb T^d)}^2\,ds\right)^{1/2}\right]\\
 &\leq\frac14\mathbb E\sup_{0\leq r\leq t}
 \|u^{N,\nu}(r\wedge\tau_R)\|_{L^2(\mathbb T^d)}^2
 +C\varepsilon\mathbb E\int_0^{t\wedge\tau_R}\sum_{k=1}^K
 \|\Pi_NG_k(u^{N,\nu}(s))\|_{L^2(\mathbb T^d)}^2\,ds.
\end{align*}
The drift term in \eqref{fixed-nu-stopped-L2-Ito} satisfies
\begin{align*}
 |\langle u^{N,\nu},\Pi_Nb^\nu(u^{N,\nu})\rangle_{L^2(\mathbb T^d)}|
 \leq\frac12\|u^{N,\nu}\|_{L^2(\mathbb T^d)}^2
 +\frac12\|\Pi_Nb^\nu(u^{N,\nu})\|_{L^2(\mathbb T^d)}^2.
\end{align*}
Taking the supremum over $r\in[0,t]$ in
\eqref{fixed-nu-stopped-L2-Ito}, then taking expectations and absorbing the
first term in the preceding estimate from the
Burkholder--Davis--Gundy inequality, we obtain
\begin{align*}
 &\mathbb E\sup_{0\leq r\leq t}
 \|u^{N,\nu}(r\wedge\tau_R)\|_{L^2(\mathbb T^d)}^2\\
 \leq& C\mathbb E\|u_0^N\|_{L^2(\mathbb T^d)}^2
 +C\int_0^t\mathbb E\sup_{0\leq\ell\leq s}
 \|u^{N,\nu}(\ell\wedge\tau_R)\|_{L^2(\mathbb T^d)}^2\,ds\\
 &+C\mathbb E\int_0^{t\wedge\tau_R}
 \|\Pi_Nb^\nu(u^{N,\nu}(s))\|_{L^2(\mathbb T^d)}^2\,ds
 +C\varepsilon\mathbb E\int_0^{t\wedge\tau_R}\sum_{k=1}^K
 \|\Pi_NG_k(u^{N,\nu}(s))\|_{L^2(\mathbb T^d)}^2\,ds.
\end{align*}
Equations \eqref{fixed-nu-initial-uniform-bound},
\eqref{fixed-nu-stopped-drift-bound}, and
\eqref{fixed-nu-stopped-noise-bound}, followed by Gronwall's lemma, prove
\eqref{fixed-nu-stopped-L2-path-bound}.
\end{proof}

The stopped estimates give global Galerkin solutions and bounds uniform
in $N$, which we use for compactness at fixed viscosity.

\begin{proposition}\label{lem-global-galerkin-fixed-nu}
Under Assumption~\ref{constant-noise-covariance}, the maximal local
solution of \eqref{fully-regularized-equation} is global.  There are
constants $C_T,C_{\nu,T}<\infty$, independent of $N$ and $C_0$, with
$C_T$ also independent of $\nu$, such that
\begin{equation}\label{fixed-nu-galerkin-estimates}
\begin{aligned}
 &\sup_{N\geq1}\Bigg\{\sup_{0\leq t\leq T}\mathbb E
 \|\nabla u^{N,\nu}(t)\|_{L^2(\mathbb T^d)}^2
 +\nu\mathbb E\int_0^T\|\Delta u^{N,\nu}(s)\|_{L^2(\mathbb T^d)}^2\,ds\\
 &+\mathbb E\int_0^T\!\int_{\mathbb T^d}
 \frac{|D^2u^{N,\nu}|^2}{H(\nabla u^{N,\nu})^2}\,dx\,ds\Bigg\}
 \leq C_T(1+C_0).
\end{aligned}
\end{equation}
The remaining path and full $H^2(\mathbb T^d)$ bounds are
\begin{equation}\label{fixed-nu-galerkin-L2-bound}
 \sup_{N\geq1}\mathbb E\sup_{0\leq t\leq T}
 \|u^{N,\nu}(t)\|_{L^2(\mathbb T^d)}^2
 \leq C_{\nu,T}(1+C_0)
\end{equation}
and
\begin{equation}\label{fixed-nu-full-H2}
 \sup_{N\geq1}\mathbb E\int_0^T
 \|u^{N,\nu}(s)\|_{H^2(\mathbb T^d)}^2\,ds
 \leq C_{\nu,T}(1+C_0).
\end{equation}
Here $C_0$ is defined in \eqref{fixed-nu-initial-uniform-bound}
and is independent of $\nu$. 
\end{proposition}

\begin{proof}

We use the three stopped estimates of Lemma~\ref{lem-stopped-L2-estimate}.
On $\{\tau_R<T\}$, continuity and the definition of $\tau_R$ imply
$\|u^{N,\nu}(\tau_R)\|_{H^1(\mathbb T^d)}\geq R$.  Therefore
\eqref{fixed-nu-stopped-L2-path-bound}
and
\eqref{fixed-nu-stopped-gradient-bound} 
 give
\begin{align*}
 R^2\mathbb P(\tau_R<T)
 &\leq\mathbb E\!\left[
 \mathbf 1_{\{\tau_R<T\}}
 \|u^{N,\nu}(\tau_R)\|_{H^1(\mathbb T^d)}^2\right]\\
 &\leq\mathbb E\|\nabla u^{N,\nu}(T\wedge\tau_R)\|_{L^2(\mathbb T^d)}^2
 +\mathbb E\sup_{0\leq t\leq T}
 \|u^{N,\nu}(t\wedge\tau_R)\|_{L^2(\mathbb T^d)}^2\leq C_{\nu,T}(1+C_0).
\end{align*}
Consequently,
\begin{align*}
 \mathbb P(\tau_R<T)&\leq
 \frac{C_{\nu,T}(1+C_0)}{R^2}\longrightarrow0, \quad\text{as} \ R\to\infty,\\
 \{\zeta^{N,\nu}\leq T\}&\subseteq\{\tau_R<T\},
 \qquad R>0,\quad \mathbb P(\zeta^{N,\nu}\leq T)
 \leq\inf_{R>0}\mathbb P(\tau_R<T)=0.
\end{align*}
The second line follows from the blow-up alternative
\eqref{fixed-nu-blow-up-alternative}.  Hence the solution is global on
$[0,T]$.

Estimate \eqref{fixed-nu-stopped-dissipation-bound} directly controls
the weighted-Hessian term in \eqref{fixed-nu-galerkin-estimates},
uniformly in $N$ and $R$.
For every global continuous finite-dimensional trajectory,
\begin{equation*}
 \sup_{0\leq s\leq T}\|u^{N,\nu}(s)\|_{H^1(\mathbb T^d)}<\infty,
 \qquad
 R>\sup_{0\leq s\leq T}\|u^{N,\nu}(s)\|_{H^1(\mathbb T^d)}
 \Longrightarrow\ \tau_R=T.
\end{equation*}
Therefore $\tau_R\uparrow T$ almost surely and, for every $t\in[0,T]$,
\begin{align*}
 \mathbb E\sup_{0\leq r\leq T}\|u^{N,\nu}(r)\|_{L^2(\mathbb T^d)}^2
 &=\lim_{R\to\infty}\mathbb E\sup_{0\leq r\leq T}
 \|u^{N,\nu}(r\wedge\tau_R)\|_{L^2(\mathbb T^d)}^2,\\
 \mathbb E\|\nabla u^{N,\nu}(t)\|_{L^2(\mathbb T^d)}^2
 &\leq\liminf_{R\to\infty}
 \mathbb E\|\nabla u^{N,\nu}(t\wedge\tau_R)\|_{L^2(\mathbb T^d)}^2,\\
 \mathbb E\int_0^T\|\Delta u^{N,\nu}(s)\|_{L^2(\mathbb T^d)}^2\,ds
 &=\lim_{R\to\infty}\mathbb E
 \int_0^{T\wedge\tau_R}\|\Delta u^{N,\nu}(s)\|_{L^2(\mathbb T^d)}^2\,ds,\\
 \mathbb E\int_0^T\!\int_{\mathbb T^d}
 \frac{|D^2u^{N,\nu}|^2}{H(\nabla u^{N,\nu})^2}\,dx\,ds
 &=\lim_{R\to\infty}\mathbb E
 \int_0^{T\wedge\tau_R}\!\int_{\mathbb T^d}
 \frac{|D^2u^{N,\nu}|^2}{H(\nabla u^{N,\nu})^2}\,dx\,ds.
\end{align*}
Monotone convergence applies in the first line because
\begin{equation*}
 \sup_{0\leq t\leq T}\|u^{N,\nu}(t\wedge\tau_R)\|_{L^2(\mathbb T^d)}^2
 =\sup_{0\leq s\leq\tau_R}\|u^{N,\nu}(s)\|_{L^2(\mathbb T^d)}^2
\end{equation*}
increases with $R$.  Fatou's lemma applies to the gradient line, and monotone
convergence applies to the nonnegative dissipation integrals.  Combining
these limits with \eqref{fixed-nu-stopped-L2-path-bound},
\eqref{fixed-nu-stopped-gradient-bound}, and
\eqref{fixed-nu-stopped-dissipation-bound} proves
\eqref{fixed-nu-galerkin-estimates} and
\eqref{fixed-nu-galerkin-L2-bound}.

Finally, periodic integration by parts and Young's inequality give
\begin{equation*}
 \|\nabla z\|_{L^2(\mathbb T^d)}^2
 =-\langle z,\Delta z\rangle_{L^2(\mathbb T^d)}
 \leq\frac12\|z\|_{L^2(\mathbb T^d)}^2+
 \frac12\|\Delta z\|_{L^2(\mathbb T^d)}^2.
\end{equation*}
Together with \eqref{periodic-hessian-laplacian-equality}, this yields
\begin{equation*}
 \|z\|_{H^2(\mathbb T^d)}^2\leq C\left(
 \|z\|_{L^2(\mathbb T^d)}^2+\|D^2z\|_{L^2(\mathbb T^d)}^2\right),
 \qquad z\in H^2(\mathbb T^d).
\end{equation*}
Consequently, \eqref{fixed-nu-galerkin-estimates} and
\eqref{fixed-nu-galerkin-L2-bound} give
\begin{align*}
 \sup_{N\geq1}\mathbb E\int_0^T\|u^{N,\nu}(s)\|_{H^2(\mathbb T^d)}^2\,ds
 &\leq C T\sup_{N\geq1}\mathbb E\sup_{0\leq s\leq T}
 \|u^{N,\nu}(s)\|_{L^2(\mathbb T^d)}^2\\
 &+C\sup_{N\geq1}\mathbb E\int_0^T
 \|\Delta u^{N,\nu}(s)\|_{L^2(\mathbb T^d)}^2\,ds
 \leq C\left(TC_{\nu,T}+\frac{C_T}{\nu}\right)(1+C_0),
\end{align*}
which is \eqref{fixed-nu-full-H2}.
\end{proof}

We also derive a uniform temporal H\"older estimate for
$u^{N,\nu}$.

\begin{proposition}\label{lem-fixed-nu-holder}
Let $u^{N,\nu}$ be the global solutions of
Proposition~\ref{lem-global-galerkin-fixed-nu}, at a fixed $\nu>0$.
For $m=3+d/2$ and every $0<\gamma<1/4$,
\begin{equation}\label{fixed-nu-holder-estimate}
 \sup_{N\geq1}\mathbb E\|u^{N,\nu}\|_{C^\gamma([0,T];H^{-m}(\mathbb T^d))}^2
 \leq C_{\nu,T,m,\gamma}(1+C_0).
\end{equation}
\end{proposition}

\begin{proof}
We write the Galerkin equation \eqref{fully-regularized-equation} as
\begin{equation*}
 \begin{aligned}
 u^{N,\nu}(t)=u_0^N+\int_0^t\Pi_Nb^\nu(u^{N,\nu}(r))\,dr+\sqrt\varepsilon\sum_{k=1}^K\int_0^t
 \Pi_NG_k(u^{N,\nu}(r))\,d\beta_r^k.
 \end{aligned}
\end{equation*}
For $0\leq s<t\leq T$, the Cauchy--Schwarz inequality yields
\begin{equation*}
 \left\|\int_s^t\Pi_Nb^\nu(u^{N,\nu}(r))\,dr\right\|_{H^{-m}(\mathbb T^d)}^2
 \leq (t-s)\int_s^t
 \|\Pi_Nb^\nu(u^{N,\nu}(r))\|_{H^{-m}(\mathbb T^d)}^2\,dr.
\end{equation*}
Hence, for every $0<\gamma<1/4$, \eqref{fixed-nu-drift-L2-bound} and
\eqref{fixed-nu-full-H2} imply
\begin{align*}
 &\mathbb E\left[\int_0^\cdot\Pi_Nb^\nu(u^{N,\nu}(r))\,dr\right]_{C^\gamma([0,T];H^{-m}(\mathbb T^d))}^2\\
 &\leq T^{1-2\gamma}\mathbb E\int_0^T
 \|\Pi_Nb^\nu(u^{N,\nu}(r))\|_{H^{-m}(\mathbb T^d)}^2\,dr\leq C_{\nu,T,m,\gamma}(1+C_0).
\end{align*}
Since the drift integral vanishes at $t=0$, this also controls its full
$C^\gamma([0,T];H^{-m}(\mathbb T^d))$ norm.
The area bound \eqref{fixed-nu-area-gradient-bound} and the pointwise
gradient moment in \eqref{fixed-nu-galerkin-estimates} imply 
\begin{equation}\label{fixed-nu-area-second-moment}
 \sup_{N\geq1}\sup_{0\leq r\leq T}\mathbb E\mathcal A(u^{N,\nu}(r))^2
 +\sup_{N\geq1}\mathbb E\int_0^T\mathcal A(u^{N,\nu}(r))^2\,dr
 \leq C_T(1+C_0).
\end{equation}
The Burkholder--Davis--Gundy inequality, the contraction
of $\Pi_N$, and \eqref{galerkin-noise-HS-bound} yield, for $0\leq s<t\leq T$,
\begin{align*}
 \mathbb E\left\|\sqrt\varepsilon\sum_{k=1}^K\int_s^t
 \Pi_NG_k(u^{N,\nu}(r))\,d\beta_r^k\right\|_{H^{-m}(\mathbb T^d)}^4\leq C\varepsilon^2\Lambda_K^2\mathbb E
 \left(\int_s^t\mathcal A(u^{N,\nu}(r))\,dr\right)^2
 \leq C_T(1+C_0)|t-s|^2.
\end{align*}
Choose $\alpha$ with $\gamma+1/4<\alpha<1/2$.  Tonelli's theorem
and the last estimate imply
\begin{align*}
 \sup_{N\geq1}\mathbb E\left[\sqrt\varepsilon\sum_{k=1}^K\int_0^\cdot
 \Pi_NG_k(u^{N,\nu}(r))\,d\beta_r^k\right]_{W^{\alpha,4}(0,T;H^{-m}(\mathbb T^d))}^4\leq C_T(1+C_0)\int_0^T\!\int_0^t
 (t-s)^{1-4\alpha}\,ds\,dt<\infty.
\end{align*}
The same increment bound with $s=0$ controls the $L^4(0,T;H^{-m}(\mathbb T^d))$ part of this
Sobolev norm because the stochastic integral vanishes at $t=0$.  The embedding
\[
W^{\alpha,4}(0,T;H^{-m}(\mathbb T^d))\hookrightarrow
C^\gamma([0,T];H^{-m}(\mathbb T^d))
\]
therefore controls the full H\"older norm of the stochastic integral in
$L^4(\Omega)$.  Finally,
$\mathbb E\|u_0^N\|_{H^{-m}(\mathbb T^d)}^2\leq C_0$ and the preceding drift
bound prove \eqref{fixed-nu-holder-estimate}.
\end{proof}

\begin{lemma}\label{lem-fixed-nu-tightness}
Let $m=3+d/2$ and fix $0<\gamma<1/4$.  A set bounded in
$$L^2(0,T;H^2(\mathbb T^d))\cap C^\gamma([0,T];H^{-m}(\mathbb T^d))$$ is relatively compact in
$L^2(0,T;H^1(\mathbb T^d))$.  An additional uniform bound in
$L^\infty(0,T;L^2(\mathbb T^d))$ implies relative compactness in
\begin{equation*}
 C([0,T];H^{-m}(\mathbb T^d)).
\end{equation*}

For the global Galerkin solutions at fixed $\nu>0$, define
\begin{equation*}
 \mathcal X_{u,\nu}:=L^2(0,T;H^1(\mathbb T^d))
 \cap L^2(0,T;H^2(\mathbb T^d))_{\mathrm w}
 \cap C([0,T];H^{-m}(\mathbb T^d)),
\end{equation*}
endowed with the topology induced by the diagonal embedding into the
product of the three displayed spaces, and set
\begin{equation}\label{fixed-nu-joint-space}
 \mathcal X_{\nu}:=\mathcal X_{u,\nu}
 \times C([0,T];\mathbb R^K)\times H^1(\mathbb T^d).
\end{equation}
Then the joint laws of
\begin{equation}\label{fixed-nu-joint-vector}
 Z^{N,\nu}:=\bigl(u^{N,\nu},\beta,u_0^N\bigr)
\end{equation}
are tight on $\mathcal X_{\nu}$.
\end{lemma}

\begin{proof}
This follows from the Aubin--Lions--Simon theorem \cite{Simon1987},
Arzel\`a--Ascoli, and Markov's inequality applied to
\eqref{fixed-nu-galerkin-L2-bound}, \eqref{fixed-nu-full-H2}, and
\eqref{fixed-nu-holder-estimate}, together with the initial convergence
\eqref{fixed-nu-initial-convergence} and the fixed Brownian law.
\end{proof}

\subsection{Identification of the viscous limit equation}
\label{subsec-fixed-nu-equation-identification}

\begin{theorem}\label{thm-fixed-nu-identification}
Suppose Assumption~\ref{constant-noise-covariance} holds.
Let $u_0$ be an $H^1(\mathbb T^d)$-valued random variable with
$\mathbb E\|u_0\|_{H^1(\mathbb T^d)}^2<\infty$.
For every $\nu>0$, equation~\eqref{viscous-equation} has a martingale
solution $u^\nu$ with initial law $\mathcal L(u_0)$.
The solution and its driving $K$-dimensional Brownian motion $\beta$
satisfy, for every $\varphi\in C^\infty(\mathbb T^d)$, almost surely for all
$t\in[0,T]$,
\begin{align}\label{fixed-nu-limit-equation}
 \langle u^\nu(t),\varphi\rangle
 ={}&\langle u^\nu(0),\varphi\rangle
 +\int_0^t\Big\langle\nu\Delta u^\nu
 +H(\nabla u^\nu)\kappa(u^\nu)
 +\lambda\frac{\nabla u^\nu\cdot D^2u^\nu\nabla u^\nu}
 {H(\nabla u^\nu)^3},\varphi\Big\rangle\,dr\notag\\
 &+\sqrt\varepsilon\sum_{k=1}^K\int_0^t
 \langle f_kH(\nabla u^\nu)^{1/2},\varphi\rangle\,d\beta_r^k.
\end{align}
\end{theorem}

\begin{proof}
\medskip
\noindent\textit{Step 1. The Jakubowski--Skorokhod argument.}
Fix $\nu>0$ and recall that $m=3+d/2$.
We consider the global Galerkin solutions from
Proposition~\ref{lem-global-galerkin-fixed-nu} as $N\to\infty$.
By Lemma~\ref{lem-fixed-nu-tightness}, the joint laws of
$Z^{N,\nu}$ are tight on $\mathcal X_\nu$.
The weak $L^2(0,T;H^2(\mathbb T^d))$ component has a countable separating
family of continuous linear functionals because this Hilbert space is
separable, and the other factors are Polish.
Thus $\mathcal X_\nu$ has a countable family of continuous functions
separating points and meets the Jakubowski--Skorokhod hypothesis.
The Jakubowski--Skorokhod theorem \cite{Jakubowski1997} therefore
yields a subsequence and random variables
\[
 (\widetilde u^{N,\nu},\widetilde\beta^N,\widetilde u_0^N),
 \qquad(\widetilde u^\nu,\widetilde\beta,\widetilde u_0)
\]
on a new probability space
$(\widetilde\Omega,\widetilde{\mathcal F},\widetilde{\mathbb P})$.
For every $N$, the represented triple has the same joint law as the
original $(u^{N,\nu},\beta,u_0^N)$, and the triples converge almost surely
in $\mathcal X_\nu$.  We retain tildes for all quantities on the new
probability space throughout this proof.  Thus, almost surely,
\begin{align}
 \widetilde u^{N,\nu}&\longrightarrow \widetilde u^\nu
 &&\text{strongly in }L^2(0,T;H^1(\mathbb T^d)),
 \label{fixed-viscosity-strong-H1}\\
 \widetilde u^{N,\nu}&\rightharpoonup \widetilde u^\nu
 &&\text{weakly in }L^2(0,T;H^2(\mathbb T^d)),
 \label{fixed-viscosity-weak-H2}\\
 \widetilde u^{N,\nu}&\longrightarrow \widetilde u^\nu
 &&\text{strongly in }C([0,T];H^{-m}(\mathbb T^d)),
 \label{fixed-viscosity-strong-path}\\
 \widetilde\beta^N&\longrightarrow\widetilde\beta
 &&\text{strongly in }C([0,T];\mathbb R^K),
 \label{fixed-viscosity-brownian-convergence}\\
 \widetilde u_0^N&\longrightarrow \widetilde u_0
 &&\text{strongly in }H^1(\mathbb T^d).
 \label{fixed-viscosity-initial-convergence}
\end{align}
Moreover,
\begin{equation}\label{fixed-viscosity-limit-initial-value}
 \widetilde u^\nu(0)=\widetilde u_0\qquad\text{in }H^{-m}(\mathbb T^d),\quad\text{almost surely}.
\end{equation}

Evaluation at time
zero is continuous on
$C([0,T];H^{-m}(\mathbb T^d))$, and $H^1(\mathbb T^d)\hookrightarrow H^{-m}(\mathbb T^d)$ continuously.
Equality of joint laws transfers the relation
$u^{N,\nu}(0)=u_0^N$ to the new space.  Hence equations
\eqref{fixed-viscosity-strong-path} and
\eqref{fixed-viscosity-initial-convergence} yield, almost surely,
\begin{align*}
 \|\widetilde u^\nu(0)-\widetilde u_0\|_{H^{-m}(\mathbb T^d)}
 &\leq\|\widetilde u^\nu(0)-\widetilde u^{N,\nu}(0)\|_{H^{-m}(\mathbb T^d)}
 +\|\widetilde u_0^N-\widetilde u_0\|_{H^{-m}(\mathbb T^d)}\longrightarrow0, \quad\text{as } N\to\infty.
\end{align*}
This completes the proof of \eqref{fixed-viscosity-limit-initial-value}. 

Let $(\widetilde{\mathcal F}_t)$ be the usual augmentation of the raw canonical
filtration generated by $(\widetilde u^\nu,\widetilde\beta,\widetilde u_0)$.
For $\varphi\in C^\infty(\mathbb T^d)$, we define the represented
Galerkin residual by
\begin{equation}\label{fixed-nu-test-martingale}
\begin{aligned}
 \widetilde M_t^{N,\nu,\varphi}
 :={}&\langle \widetilde u^{N,\nu}(t)-\widetilde u_0^N,\varphi\rangle
 -\int_0^t\langle b^\nu(\widetilde u^{N,\nu}(r)),\Pi_N\varphi\rangle\,dr.
\end{aligned}
\end{equation}

Define also the limiting residual by
\begin{equation}\label{fixed-nu-limit-martingale}
\begin{aligned}
 \widetilde M_t^{\nu,\varphi}:={}&
 \langle \widetilde u^\nu(t)-\widetilde u^\nu(0),\varphi\rangle
 -\int_0^t\langle b^\nu(\widetilde u^\nu(r)),\varphi\rangle\,dr.
\end{aligned}
\end{equation}
The drift and noise coefficients are viewed as Borel functions on $H^{-m}(\mathbb T^d)$,
extended by zero outside their Sobolev domains.  

\medskip
\noindent\textit{Step 2. Passage to the limit in the drift.}
We show that, almost surely,
\begin{align}
 \left(I-\frac{\nabla \widetilde u^{N,\nu}\otimes\nabla \widetilde u^{N,\nu}}{H(\nabla \widetilde u^{N,\nu})^2}\right):D^2\widetilde u^{N,\nu}
 &\rightharpoonup \left(I-\frac{\nabla \widetilde u^\nu\otimes\nabla \widetilde u^\nu}{H(\nabla \widetilde u^\nu)^2}\right):D^2\widetilde u^\nu,
 \label{fixed-nu-mean-curvature-convergence}\\
 \frac{\nabla \widetilde u^{N,\nu}\cdot D^2\widetilde u^{N,\nu}\nabla \widetilde u^{N,\nu}}
 {H(\nabla \widetilde u^{N,\nu})^3}
 &\rightharpoonup
 \frac{\nabla \widetilde u^\nu\cdot D^2\widetilde u^\nu\nabla \widetilde u^\nu}
 {H(\nabla \widetilde u^\nu)^3}
 \label{fixed-nu-correction-convergence}
\end{align}
weakly in $L^2((0,T)\times\mathbb T^d)$, and verify that the limiting
drift has a predictable, square-integrable $L^2(\mathbb T^d)$-valued modification. 

Due to the expression of \eqref{fixed-nu-drift-definition}, we define the following two matrix maps
$$
 A_1(q):=I-\frac{q\otimes q}{H(q)^2}\in C_b(\mathbb R^d;\mathbb R^{d\times d}),
 \qquad
 A_2(q):=\frac{q\otimes q}{H(q)^3}\in C_b(\mathbb R^d;\mathbb R^{d\times d}),
 \qquad q\in\mathbb R^d, 
$$
for convenience. By \eqref{fixed-viscosity-strong-H1}, every subsequence has a further
subsequence along which the gradients converge almost everywhere.
Continuity and boundedness of the matrix maps, followed by dominated
convergence on the finite time--space domain, therefore give
\begin{equation}\label{fixed-nu-matrix-coefficient-convergence}
\begin{aligned}
 I-\frac{\nabla \widetilde u^{N,\nu}\otimes\nabla \widetilde u^{N,\nu}}
 {H(\nabla \widetilde u^{N,\nu})^2}
 &\longrightarrow
 I-\frac{\nabla \widetilde u^\nu\otimes\nabla \widetilde u^\nu}
 {H(\nabla \widetilde u^\nu)^2},\\
 \frac{\nabla \widetilde u^{N,\nu}\otimes\nabla \widetilde u^{N,\nu}}
 {H(\nabla \widetilde u^{N,\nu})^3}
 &\longrightarrow
 \frac{\nabla \widetilde u^\nu\otimes\nabla \widetilde u^\nu}
 {H(\nabla \widetilde u^\nu)^3}
\end{aligned}
\end{equation}
strongly in $L^2((0,T)\times\mathbb T^d)$, almost surely.
For every scalar test function $\phi\in L^2((0,T)\times\mathbb T^d)$, combining the same
subsequence argument and dominated convergence theorem with 
$$
 |[A_i(\nabla \widetilde u^{N,\nu})-A_i(\nabla \widetilde u^\nu)]\phi|^2
 \leq4\|A_i\|_{L^\infty(\mathbb R^d)}^2|\phi|^2,
 \qquad i=1,2,
$$
we have 
\begin{equation*}
 \|[A_i(\nabla \widetilde u^{N,\nu})-A_i(\nabla \widetilde u^\nu)]\phi\|_{L^2((0,T)\times\mathbb T^d)}
 \longrightarrow0,\qquad i=1,2.
\end{equation*}
Since $D^2\widetilde u^{N,\nu}\rightharpoonup D^2\widetilde u^\nu$ in $L^2((0,T)\times\mathbb T^d)$ by
\eqref{fixed-viscosity-weak-H2}, for each $i=1,2$ we obtain
\begin{align*}
 \int_0^T\!\int_{\mathbb T^d}
 [A_i(\nabla \widetilde u^{N,\nu}):D^2\widetilde u^{N,\nu}]\phi\,dx\,dt&=\int_0^T\!\int_{\mathbb T^d}
 D^2\widetilde u^{N,\nu}:[A_i(\nabla \widetilde u^{N,\nu})\phi]\,dx\,dt\\
 &\longrightarrow\int_0^T\!\int_{\mathbb T^d}
 D^2\widetilde u^\nu:[A_i(\nabla \widetilde u^\nu)\phi]\,dx\,dt.
\end{align*}
Taking $i=1,2$, together with
\eqref{fixed-nu-drift-definition}, we have 
\eqref{fixed-nu-mean-curvature-convergence} and
\eqref{fixed-nu-correction-convergence}.  The linear Laplacian term
converges weakly in the same space by
\eqref{fixed-viscosity-weak-H2}.  Therefore
\begin{equation}\label{fixed-nu-drift-convergence}
 b^\nu(\widetilde u^{N,\nu})\rightharpoonup b^\nu(\widetilde u^\nu)
 \quad\text{weakly in }L^2((0,T)\times\mathbb T^d)
 \quad\text{almost surely}.
\end{equation}

By \eqref{matrix-curvature-formula}, we have
\begin{equation*}
 H(\nabla z)\kappa(z)
 =\left(I-\frac{\nabla z\otimes\nabla z}{H(\nabla z)^2}\right):D^2z,
 \qquad z\in H^2(\mathbb T^d).
\end{equation*}
Both matrix coefficients are bounded, so the nonlinear maps take $H^2(\mathbb T^d)$
into $L^2(\mathbb T^d)$.  Their continuity from $H^2(\mathbb T^d)$ to $L^2(\mathbb T^d)$ follows directly from
Lemma~\ref{lem-quadratic-hessian-continuity}, applied to the two bounded
continuous matrix maps above.  Since $\Delta:H^2(\mathbb T^d)\to L^2(\mathbb T^d)$ is continuous,
so is $b^\nu:H^2(\mathbb T^d)\to L^2(\mathbb T^d)$.  The Lusin--Souslin theorem
\cite[Section~15]{Kechris1995} yields Borel extensions of these maps to
$H^{-m}(\mathbb T^d)$. Composing them with the adapted, continuous process $\widetilde u^\nu$
provides predictable versions of the coefficients.
The bounded matrix coefficients, weak lower semicontinuity in
\eqref{fixed-nu-mean-curvature-convergence}--\eqref{fixed-nu-drift-convergence}, Fatou's lemma, equality of joint laws, and
\eqref{fixed-nu-full-H2} yield
\begin{align*}
 &\widetilde{\mathbb E}\int_0^T\|H(\nabla \widetilde u^\nu)\kappa(\widetilde u^\nu)\|_{L^2(\mathbb T^d)}^2\,dr
 +\widetilde{\mathbb E}\int_0^T\|b^\nu(\widetilde u^\nu)\|_{L^2(\mathbb T^d)}^2\,dr+\widetilde{\mathbb E}\int_0^T
 \left\|\frac{\nabla \widetilde u^\nu\cdot D^2\widetilde u^\nu\nabla \widetilde u^\nu}
 {H(\nabla \widetilde u^\nu)^3}\right\|_{L^2(\mathbb T^d)}^2\,dr\\
 \leq& C_\nu\liminf_{N\to\infty}\widetilde{\mathbb E}\int_0^T
 \|D^2\widetilde u^{N,\nu}(r)\|_{L^2(\mathbb T^d)}^2\,dr
 \leq C_{\nu,T}(1+C_0).
\end{align*}

\medskip
\noindent\textit{Step 3. Identification and passage to the limit of the martingales.}
We first claim that for every $k=1,\ldots,K$, almost surely,
\begin{equation}\label{fixed-nu-projected-noise-convergence}
 \Pi_NG_k(\widetilde u^{N,\nu})\longrightarrow G_k(\widetilde u^\nu)
 \quad\text{strongly in }L^2(0,T;L^2(\mathbb T^d)).
\end{equation}
Here $G_k$ is defined in \eqref{galerkin-noise-definition}.

The scalar map $q\mapsto H(q)^{1/2}$ is globally Lipschitz because
\begin{equation*}
 \left|\nabla_qH(q)^{1/2}\right|
 =\frac{|q|}{2H(q)^{3/2}}\leq\frac12.
\end{equation*}
Thus, for every $k$, it follows from \eqref{fixed-viscosity-strong-H1} that
\begin{align*}
 \|G_k(\widetilde u^{N,\nu})-G_k(\widetilde u^\nu)\|_{L^2(0,T;L^2(\mathbb T^d))}
 &\leq\frac12\|f_k\|_{L^\infty(\mathbb T^d)}
 \|\nabla \widetilde u^{N,\nu}-\nabla \widetilde u^\nu\|_{L^2(0,T;L^2(\mathbb T^d))}
 \longrightarrow0.
\end{align*}
Since $\Pi_N$ is an $L^2(\mathbb T^d)$ contraction and converges strongly to the identity,
\begin{align*}
 &\|\Pi_NG_k(\widetilde u^{N,\nu})-G_k(\widetilde u^\nu)\|_{L^2(0,T;L^2(\mathbb T^d))}\\
 \leq&\|G_k(\widetilde u^{N,\nu})-G_k(\widetilde u^\nu)\|_{L^2(0,T;L^2(\mathbb T^d))}
 +\|\Pi_NG_k(\widetilde u^\nu)-G_k(\widetilde u^\nu)\|_{L^2(0,T;L^2(\mathbb T^d))}
 \longrightarrow0.
\end{align*}

We next transfer the Galerkin martingale identities to the new
space and establish uniform integrability.
On the original probability space, write $\beta^N:=\beta$ for the noise
driving the $N$th Galerkin equation, and define $M^{N,\nu,\varphi}$
by the same formula as in
\eqref{fixed-nu-test-martingale}, with
$(\widetilde u^{N,\nu},\widetilde u_0^N)$ replaced by $(u^{N,\nu},u_0^N)$.  Then
\begin{equation}\label{fixed-nu-galerkin-martingale-representation}
 M_t^{N,\nu,\varphi}
 =\sqrt\varepsilon\sum_{k=1}^K
 \int_0^t\langle G_k(u^{N,\nu}(r)),\Pi_N\varphi\rangle\,d\beta_{r}^{N,k}.
\end{equation}
The It\^o product rule yields the identities below on the original
space.  For fixed $N$, the Fourier projection extends continuously
from $H^{-m}(\mathbb T^d)$ to $\Pi_N L^2(\mathbb T^d)$ and hence from
$C([0,T];H^{-m}(\mathbb T^d))$ to $C([0,T];\Pi_N L^2(\mathbb T^d))$.
Since $b^\nu$ and $G_k$ are continuous on $\Pi_N L^2(\mathbb T^d)$,
the residuals and compensators are Borel functions of the projected
path and initial datum.  They depend only on past values, so equality
of the joint laws transfers the martingale identities to the represented
copies in their raw canonical filtration.
In particular, all of the following are martingales:
\begin{align}\label{fixed-nu-product-martingales}
 &\widetilde M^{N,\nu,\varphi},\qquad \widetilde\beta^{N,k},\qquad
 \widetilde\beta^{N,k}\widetilde\beta^{N,\ell}-\delta_{k\ell}(\cdot),\notag\\
 &\widetilde M^{N,\nu,\varphi}\widetilde\beta^{N,k}
 -\sqrt\varepsilon\int_0^{\cdot}\langle G_k(\widetilde u^{N,\nu}(r)),\Pi_N\varphi\rangle\,dr,\notag\\
 &\widetilde M^{N,\nu,\varphi}\widetilde M^{N,\nu,\psi}
 -\varepsilon\sum_{k=1}^K\int_0^{\cdot}
 \langle G_k(\widetilde u^{N,\nu}(r)),\Pi_N\varphi\rangle\langle G_k(\widetilde u^{N,\nu}(r)),\Pi_N\psi\rangle\,dr.
\end{align}

We verify uniform integrability on the represented space, using equality
of joint laws to transfer the moment bounds from the original space.  By the contraction of
$\Pi_N$, \eqref{galerkin-noise-HS-bound}, and the Cauchy--Schwarz inequality,
\begin{align*}
 \sum_{k=1}^K|\langle G_k(\widetilde u^{N,\nu}(r)),\Pi_N\varphi\rangle|^2
 \leq\|\varphi\|_{L^2(\mathbb T^d)}^2
 \sum_{k=1}^K\|G_k(\widetilde u^{N,\nu}(r))\|_{L^2(\mathbb T^d)}^2=\Lambda_K\|\varphi\|_{L^2(\mathbb T^d)}^2
 \mathcal A(\widetilde u^{N,\nu}(r)).
\end{align*}
Furthermore, \eqref{fixed-nu-area-second-moment} gives
\begin{align*}
 \sup_{N\geq1}\widetilde{\mathbb E}\left(\int_0^T
 \mathcal A(\widetilde u^{N,\nu}(r))\,dr\right)^2
 &\leq T\sup_{N\geq1}\widetilde{\mathbb E}\int_0^T
 \mathcal A(\widetilde u^{N,\nu}(r))^2\,dr<\infty.
\end{align*}
The Burkholder--Davis--Gundy inequality applied to
\eqref{fixed-nu-galerkin-martingale-representation} and equality of joint laws yield
\begin{align}\label{fixed-nu-martingale-uniform-integrability}
 \sup_{N\geq1}\widetilde{\mathbb E}\sup_{0\leq t\leq T}
 |\widetilde M_t^{N,\nu,\varphi}|^4
 &\leq C\varepsilon^2\sup_{N\geq1}\widetilde{\mathbb E}\left(
 \sum_{k=1}^K\int_0^T|\langle G_k(\widetilde u^{N,\nu}(r)),\Pi_N\varphi\rangle|^2\,dr
 \right)^2\notag\\
 &\leq C_{\varphi,\varepsilon,\Lambda_K}
 \sup_{N\geq1}\widetilde{\mathbb E}\left(\int_0^T
 \mathcal A(\widetilde u^{N,\nu}(r))\,dr\right)^2<\infty.
\end{align}
The squared suprema of the two compensators are bounded respectively by
\begin{equation*}
 T\int_0^T|\langle G_k(\widetilde u^{N,\nu}(r)),\Pi_N\varphi\rangle|^2\,dr
 \quad\text{and}\quad
 C_{\varphi,\psi}\left(\int_0^T\mathcal A(\widetilde u^{N,\nu}(r))\,dr\right)^2.
\end{equation*}
The latter has bounded expectation by \eqref{fixed-nu-area-second-moment}.
Using these bounds, the Cauchy--Schwarz inequality, the Brownian fourth-moment bound, and
\eqref{fixed-nu-martingale-uniform-integrability}, we obtain
\begin{align*}
 &\sup_{N\geq1}\widetilde{\mathbb E}\sup_{0\leq t\leq T}
 \bigl(|\widetilde M_t^{N,\nu,\varphi}\widetilde\beta^{N,k}(t)|^2
 +|\widetilde M_t^{N,\nu,\varphi}\widetilde M_t^{N,\nu,\psi}|^2\bigr)<\infty,\\
 &\sup_{N\geq1}\widetilde{\mathbb E}\sup_{0\leq t\leq T}
 \left|\int_0^t \langle G_k(\widetilde u^{N,\nu}(r)),\Pi_N\varphi\rangle\,dr\right|^2<\infty,\\
 &\sup_{N\geq1}\widetilde{\mathbb E}\sup_{0\leq t\leq T}
 \left|\sum_{j=1}^K\int_0^t
 \langle G_j(\widetilde u^{N,\nu}(r)),\Pi_N\varphi\rangle\langle G_j(\widetilde u^{N,\nu}(r)),\Pi_N\psi\rangle\,dr\right|^2<\infty.
\end{align*}
The same Brownian moment bound applies to
$\widetilde\beta^{N,k}\widetilde\beta^{N,\ell}-\delta_{k\ell}(\cdot)$.
Thus every process in \eqref{fixed-nu-product-martingales} has an
integrable time supremum uniformly in $N$, and its values at fixed
times are uniformly integrable.

We now pass the residuals and their compensators to the limit.
For every $t\in[0,T]$ and every smooth $\varphi$,
\eqref{fixed-nu-drift-convergence} and
$\mathbf 1_{[0,t]}\Pi_N\varphi\to\mathbf 1_{[0,t]}\varphi$
strongly in $L^2((0,T)\times\mathbb T^d)$ derive 
\begin{align}\label{fixed-nu-tested-drift-convergence}
 &\int_0^t\langle b^\nu(\widetilde u^{N,\nu}(r)),\Pi_N\varphi\rangle\,dr
 \longrightarrow
 \int_0^t\langle b^\nu(\widetilde u^\nu(r)),\varphi\rangle\,dr.
\end{align}
Equations
\eqref{fixed-viscosity-strong-path} and
\eqref{fixed-viscosity-initial-convergence} then give
$\widetilde M_t^{N,\nu,\varphi}\to \widetilde M_t^{\nu,\varphi}$ almost surely.
By \eqref{fixed-nu-projected-noise-convergence},
\begin{equation*}
 \langle G_k(\widetilde u^{N,\nu}),\Pi_N\varphi\rangle\longrightarrow
 \langle G_k(\widetilde u^\nu),\varphi\rangle
 \quad\text{strongly in }L^2(0,T)\quad\text{almost surely}.
\end{equation*}
Consequently, the Cauchy--Schwarz inequality yields, almost surely and uniformly in
$t\in[0,T]$,
\begin{align*}
 \int_0^t \langle G_k(\widetilde u^{N,\nu}(r)),\Pi_N\varphi\rangle\,dr
 &\longrightarrow\int_0^t \langle G_k(\widetilde u^\nu(r)),\varphi\rangle\,dr,\\
 \sum_{k=1}^K\int_0^t
 \langle G_k(\widetilde u^{N,\nu}(r)),\Pi_N\varphi\rangle\langle G_k(\widetilde u^{N,\nu}(r)),\Pi_N\psi\rangle\,dr
 &\longrightarrow\sum_{k=1}^K\int_0^t
 \langle G_k(\widetilde u^\nu(r)),\varphi\rangle\langle G_k(\widetilde u^\nu(r)),\psi\rangle\,dr.
\end{align*}

For the noise, the global Lipschitz continuity of
$q\mapsto H(q)^{1/2}$ implies
that
\begin{equation*}
 H^1(\mathbb T^d)\ni z\longmapsto
 \langle f_kH(\nabla z)^{1/2},\varphi\rangle\in\mathbb R
\end{equation*}
is continuous.  The Borel-extension argument of
Step 2, applied to
$H^1(\mathbb T^d)\hookrightarrow H^{-m}(\mathbb T^d)$, gives a predictable version of
$\langle G_k(\widetilde u^\nu),\varphi\rangle$ because $\widetilde u^\nu$ has continuous adapted
$H^{-m}(\mathbb T^d)$-paths and belongs to $H^1(\mathbb T^d)$ almost everywhere in time.
Fatou's lemma and the coefficient bound in
the preceding Galerkin estimates give
\begin{align*}
 \widetilde{\mathbb E}\sum_{k=1}^K\int_0^T
 |\langle G_k(\widetilde u^\nu(r)),\varphi\rangle|^2\,dr
 &\leq\liminf_{N\to\infty}
 \widetilde{\mathbb E}\sum_{k=1}^K\int_0^T
 |\langle G_k(\widetilde u^{N,\nu}(r)),\Pi_N\varphi\rangle|^2\,dr\\
 &\leq C_{\varphi,\Lambda_K}
 \liminf_{N\to\infty}\widetilde{\mathbb E}\int_0^T
 \mathcal A(\widetilde u^{N,\nu}(r))\,dr<\infty.
\end{align*}
Thus the stochastic integrals in Theorem~\ref{thm-fixed-nu-identification}
are well defined.

In the following, we pass the martingale identities. Let $0\leq r_1,\ldots,r_J\leq s<t\leq T$, let
$\varphi_1,\ldots,\varphi_J\in C^\infty(\mathbb T^d)$, and let
$\Gamma\in C_b(\mathbb R^{J(K+2)})$.  The first martingale identity in
\eqref{fixed-nu-product-martingales} shows that 
\begin{equation*}
 \widetilde{\mathbb E}\left[
 (\widetilde M_t^{N,\nu,\varphi}-\widetilde M_s^{N,\nu,\varphi})
 \Gamma\left(
 \bigl(\langle \widetilde u^{N,\nu}(r_j),\varphi_j\rangle,
 \widetilde\beta^N(r_j),\langle \widetilde u_0^N,\varphi_j\rangle\bigr)_{j=1}^J
 \right)\right]=0.
\end{equation*}
The almost sure convergences above and
\eqref{fixed-nu-martingale-uniform-integrability} allow us to pass to
the limit by Vitali's theorem and obtain
\begin{equation*}
 \widetilde{\mathbb E}\left[
 (\widetilde M_t^{\nu,\varphi}-\widetilde M_s^{\nu,\varphi})
 \Gamma\left(
 \bigl(\langle\widetilde u^\nu(r_j),\varphi_j\rangle,
 \widetilde\beta(r_j),\langle\widetilde u_0,\varphi_j\rangle\bigr)_{j=1}^J
 \right)\right]=0.
\end{equation*}
The same argument, using the product estimates, applies to the Brownian
identities and all products in \eqref{fixed-nu-product-martingales}.

The limiting processes have continuous paths
and, by Fatou's lemma, \eqref{fixed-nu-martingale-uniform-integrability},
and the product estimates, integrable time suprema.
Their $L^1(\widetilde\Omega)$ continuity allows us to let $s_j\downarrow s$
in the increment identities, so the martingale property holds in the
right-continuous filtration as well,  then these identities hold in the usual augmentation. 

In particular, L\'evy's characterization shows that $\widetilde\beta$ is a standard
$K$-dimensional Brownian motion in that filtration, and
\begin{align}
 \langle \widetilde M^{\nu,\varphi},\widetilde\beta^k\rangle_t
 &=\sqrt\varepsilon\int_0^t
 \langle G_k(\widetilde u^\nu(r)),\varphi\rangle\,dr,\label{fixed-nu-limit-cross-variation}\\
 \langle \widetilde M^{\nu,\varphi},\widetilde M^{\nu,\psi}\rangle_t
 &=\varepsilon\sum_{k=1}^K\int_0^t
 \langle G_k(\widetilde u^\nu(r)),\varphi\rangle\langle G_k(\widetilde u^\nu(r)),\psi\rangle\,dr.
 \label{fixed-nu-limit-quadratic-variation}
\end{align}

Equations \eqref{fixed-nu-limit-cross-variation} and
\eqref{fixed-nu-limit-quadratic-variation} give
\begin{align*}
 &\left\langle \widetilde M^{\nu,\varphi}
 -\sqrt\varepsilon\sum_{k=1}^K
 \int_0^{\cdot}\langle G_k(\widetilde u^\nu(r)),\varphi\rangle\,d\widetilde\beta_r^k\right\rangle_t\\
 =&\langle \widetilde M^{\nu,\varphi}\rangle_t
 +\varepsilon\sum_{k=1}^K\int_0^t|\langle G_k(\widetilde u^\nu(r)),\varphi\rangle|^2\,dr
 -2\varepsilon\sum_{k=1}^K\int_0^t|\langle G_k(\widetilde u^\nu(r)),\varphi\rangle|^2\,dr=0.
\end{align*}
This shows \eqref{fixed-nu-limit-equation}. 

Mean-square convergence in \eqref{fixed-nu-initial-convergence} implies
$\mathcal L(u_0^N)\Rightarrow\mathcal L(u_0)$ on $H^1(\mathbb T^d)$.
Since the represented initial coordinates have the same laws as the
original ones and converge almost surely in $H^1(\mathbb T^d)$, their limit has law
$\mathcal L(u_0)$.  By
\eqref{fixed-viscosity-limit-initial-value}, this is also the law of
$\widetilde u^\nu(0)$.  This completes the proof. 
\end{proof}

For the remainder of this section, we omit tildes from quantities on
this new probability space.

\subsection{Regularity and energy estimates for the viscous solution}
\label{subsec-fixed-nu-regularity}

We first bound the coefficients to justify the Hilbert-space stochastic
integral and gradient equation.  We then prove path continuity and the
Dirichlet inequality for the limit.

\begin{lemma}\label{lem-fixed-nu-coefficient-regularity}
The process obtained in Theorem~\ref{thm-fixed-nu-identification} satisfies
$u^\nu\in L^2(\Omega;L^2(0,T;H^2(\mathbb T^d)))$ and
\begin{equation}\label{fixed-nu-limit-drift-regularity}
 b^\nu(u^\nu)\in
 L^2(\Omega\times(0,T);L^2(\mathbb T^d)),
\end{equation}
and
\begin{equation}\label{fixed-nu-limit-noise-H1-regularity}
 \sum_{k=1}^K\mathbb E\int_0^T
 \|G_k(u^\nu(s))\|_{H^1(\mathbb T^d)}^2\,ds<\infty.
\end{equation}
\end{lemma}

\begin{proof}
We again use the representation of
Step 1 of the proof of Theorem~\ref{thm-fixed-nu-identification}.
For almost every represented outcome, weak lower semicontinuity in
\eqref{fixed-viscosity-weak-H2} gives
\begin{equation*}
 \int_0^T\|u^\nu(s)\|_{H^2(\mathbb T^d)}^2\,ds
 \leq\liminf_{N\to\infty}\int_0^T
 \|u^{N,\nu}(s)\|_{H^2(\mathbb T^d)}^2\,ds.
\end{equation*}
Fatou's lemma and \eqref{fixed-nu-full-H2} therefore give
\begin{equation}\label{fixed-nu-limit-H2-moment}
 \mathbb E\int_0^T\|u^\nu(s)\|_{H^2(\mathbb T^d)}^2\,ds
 \leq\liminf_{N\to\infty}\mathbb E\int_0^T
 \|u^{N,\nu}(s)\|_{H^2(\mathbb T^d)}^2\,ds<\infty.
\end{equation}
The drift assertion \eqref{fixed-nu-limit-drift-regularity} was proved in
Step 2 of the proof of Theorem~\ref{thm-fixed-nu-identification}.
The Sobolev chain rule extends the noise-gradient identity in
Lemma~\ref{lem-projected-dirichlet-energy} to $z\in H^2(\mathbb T^d)$
and gives
\begin{align*}
 \sum_{k=1}^K|\nabla G_k(z)|^2
 &\leq\|\operatorname{tr}C_K\|_{L^\infty(\mathbb T^d)}H(\nabla z)
 +\frac{\Lambda_K}{4}|D^2z|^2.
\end{align*}
Together with \eqref{galerkin-noise-HS-bound}, it follows that 
\begin{align*}
 \sum_{k=1}^K\|G_k(z)\|_{H^1(\mathbb T^d)}^2
 &\leq C_{f_1,\ldots,f_K}
 \left(\mathcal A(z)+\|D^2z\|_{L^2(\mathbb T^d)}^2\right),
 \qquad z\in H^2(\mathbb T^d).
\end{align*}
Substituting $z=u^\nu(s)$, integrating over
$\Omega\times(0,T)$, and using
\eqref{fixed-nu-limit-H2-moment} and the bound
\begin{equation*}
 \mathcal A(z)\lesssim1+
 \|\nabla z\|_{L^2(\mathbb T^d)}
\end{equation*}
proves
\eqref{fixed-nu-limit-noise-H1-regularity}.
The same chain rule proves continuity of
$G_k:H^2(\mathbb T^d)\to H^1(\mathbb T^d)$.  If $z_n\to z$ in $H^2(\mathbb T^d)$, the global Lipschitz property of
$q\mapsto H(q)^{1/2}$ implies $G_k(z_n)\to G_k(z)$ in $L^2(\mathbb T^d)$.
The vector map $q\mapsto q/H(q)^{3/2}$ is bounded and continuous.
Using this property and dominated convergence, we obtain
\begin{align*}
 &\left\|D^2z_n\frac{\nabla z_n}{H(\nabla z_n)^{3/2}}
 -D^2z\frac{\nabla z}{H(\nabla z)^{3/2}}\right\|_{L^2(\mathbb T^d)}\\
\leq& C\|D^2z_n-D^2z\|_{L^2(\mathbb T^d)}
 +\left\|D^2z\left(
 \frac{\nabla z_n}{H(\nabla z_n)^{3/2}}
 -\frac{\nabla z}{H(\nabla z)^{3/2}}
 \right)\right\|_{L^2(\mathbb T^d)}\longrightarrow0.
\end{align*}
Using the gradient formula \eqref{noise-spatial-gradient-formula}, we obtain
\begin{equation*}
 G_k(z_n)\longrightarrow G_k(z)\quad\text{in }H^1(\mathbb T^d).
\end{equation*}
As in Step 2 of the proof of Theorem~\ref{thm-fixed-nu-identification}, the
Lusin--Souslin theorem for $H^2(\mathbb T^d)\hookrightarrow H^{-m}(\mathbb T^d)$ yields a Borel
extension of $G_k:H^2(\mathbb T^d)\to H^1(\mathbb T^d)$.  Composing it with the adapted,
$H^{-m}(\mathbb T^d)$-continuous process $u^\nu$ yields a predictable $H^1(\mathbb T^d)$-valued
version agreeing with $G_k(u^\nu)$ for $dt\otimes d\mathbb P$-almost
every $(t,\omega)$.
\end{proof}

Using these coefficient bounds, we prove the path regularity needed
to evaluate the energy at every time.

\begin{proposition}\label{lem-fixed-nu-path-regularity}
Let $u^\nu$ be the represented solution of
Theorem~\ref{thm-fixed-nu-identification}.
The process $u^\nu$ has an almost surely $H^1(\mathbb T^d)$-continuous modification and
\begin{equation}\label{fixed-viscosity-regularity}
 u^\nu\in L^2\bigl(\Omega;L^2(0,T;H^2(\mathbb T^d))\bigr)
 \cap L^2\bigl(\Omega;C([0,T];H^{-m}(\mathbb T^d))\bigr),
 \qquad m=3+d/2.
\end{equation}
Almost surely, \eqref{fixed-nu-limit-equation} holds in its
$L^2(\mathbb T^d)$-valued integral form for every $t\in[0,T]$.
\end{proposition}

\begin{proof}
By Lemma~\ref{lem-fixed-nu-coefficient-regularity}, we have
\begin{equation*}
 b^\nu(u^\nu)\in L^2(\Omega\times(0,T);L^2(\mathbb T^d)),
 \qquad
 G_k(u^\nu)\in L^2(\Omega\times(0,T);H^1(\mathbb T^d)),
\end{equation*}
with predictable noise coefficients.  Their Bochner and
Hilbert-space stochastic integrals are therefore, respectively,
$L^2(\mathbb T^d)$- and $H^1(\mathbb T^d)$-continuous.  Thus the integral form of
\eqref{fixed-nu-limit-equation} has an adapted, $L^2(\mathbb T^d)$-continuous right-hand side.
Choose a countable dense subset of $H^m(\mathbb T^d)$ and first impose
\eqref{fixed-nu-limit-equation} for these test functions at rational times.
Both sides are continuous in $H^{-m}(\mathbb T^d)$, so
\begin{align*}
 u^\nu(t)={}&u^\nu(0)+\int_0^t b^\nu(u^\nu(s))\,ds
 +\sqrt\varepsilon\sum_{k=1}^K
 \int_0^tG_k(u^\nu(s))\,d\beta_s^k
 \quad\text{in }H^{-m}(\mathbb T^d)
\end{align*}
for every $t\in[0,T]$ on one probability-one set.  Thus the right-hand side
defines an $L^2(\mathbb T^d)$-continuous modification of $u^\nu$, and the
displayed identity also holds in $L^2(\mathbb T^d)$.
We use this modification below.
Doob's inequality, \eqref{fixed-nu-limit-drift-regularity} and
\eqref{fixed-nu-limit-noise-H1-regularity} give
\begin{align}
 \mathbb E\sup_{0\leq t\leq T}\|u^\nu(t)\|_{L^2(\mathbb T^d)}^2
 &\leq 3\mathbb E\|u_0\|_{L^2(\mathbb T^d)}^2
 +3T\mathbb E\int_0^T\|b^\nu(u^\nu(s))\|_{L^2(\mathbb T^d)}^2\,ds\notag\\
 &\quad+12\varepsilon\sum_{k=1}^K\mathbb E\int_0^T
 \|G_k(u^\nu(s))\|_{L^2(\mathbb T^d)}^2\,ds<\infty.\label{fixed-nu-limit-path-moment}
\end{align}
Applying the bounded map $\nabla:L^2(\mathbb T^d)\to H^{-1}(\mathbb T^d;\mathbb R^d)$ to the
$L^2(\mathbb T^d)$-valued integral form of \eqref{fixed-nu-limit-equation}, we have that almost surely for every $t\in[0,T]$,
\begin{equation}\label{fixed-nu-gradient-variational-equation}
 \nabla u^\nu(t)
 =\nabla u^\nu(0)
 +\int_0^t\nabla b^\nu(u^\nu(s))\,ds+\sqrt\varepsilon\sum_{k=1}^K
 \int_0^t\nabla G_k(u^\nu(s))\,d\beta_s^k
\end{equation}
holds in $H^{-1}(\mathbb T^d;\mathbb R^d)$. Furthermore, by \eqref{fixed-nu-limit-drift-regularity},
\eqref{fixed-nu-limit-noise-H1-regularity} and \eqref{fixed-nu-limit-H2-moment}, we have 
\begin{align*}
 &\nabla u^\nu\in L^2(\Omega\times(0,T);H^1(\mathbb T^d)),\quad
 \nabla u^\nu\in C([0,T];H^{-1}(\mathbb T^d))\quad\text{a.s.},\\
 &\nabla b^\nu(u^\nu)\in L^2(\Omega\times(0,T);H^{-1}(\mathbb T^d)),\quad
 \nabla G_k(u^\nu)\in L^2(\Omega\times(0,T);L^2(\mathbb T^d)),\quad
 \nabla u_0\in L^2(\Omega;L^2(\mathbb T^d)).
\end{align*}
Since $\nabla u^\nu$ is adapted and continuous in $H^{-1}(\mathbb T^d)$, it is
progressively measurable as an $H^{-1}(\mathbb T^d)$-valued process.
The first bound above also shows that it belongs to $H^1(\mathbb T^d)$ for
$dt\otimes d\mathbb P$-almost every $(t,\omega)$.
To view it as an $H^1(\mathbb T^d)$-valued progressively measurable process, let
$R:H^{-1}(\mathbb T^d;\mathbb R^d)\to H^1(\mathbb T^d;\mathbb R^d)$
be given by
\[
 R(z):=
 \begin{cases}
 z,&z\in H^1(\mathbb T^d;\mathbb R^d),\\
 0,&\text{otherwise}.
 \end{cases}
\]
The Lusin--Souslin theorem applied to the continuous injection
$H^1(\mathbb T^d)\hookrightarrow H^{-1}(\mathbb T^d)$ ensures that $R$ is Borel measurable.
Thus $R(\nabla u^\nu)$ is progressively measurable with values in
$H^1(\mathbb T^d)$ and agrees with $\nabla u^\nu$ almost everywhere; we use this
representative below.  The drift and
diffusion are progressively measurable by
Step 2 of the proof of Theorem~\ref{thm-fixed-nu-identification} and
Lemma~\ref{lem-fixed-nu-coefficient-regularity}.  For the vector-valued Gelfand
triple $H^1(\mathbb T^d)\hookrightarrow L^2(\mathbb T^d)\hookrightarrow H^{-1}(\mathbb T^d)$,
\cite[Theorem~4.2.5]{PrevotRockner2007} applied to
\eqref{fixed-nu-gradient-variational-equation} therefore yields
\begin{equation}\label{fixed-nu-gradient-continuity}
 \nabla u^\nu\in C([0,T];L^2(\mathbb T^d;\mathbb R^d))
 \quad\text{almost surely}.
\end{equation}
The gradient of the $L^2(\mathbb T^d)$-continuous modification of $u^\nu$ equals
the $L^2(\mathbb T^d)$-continuous gradient representative from
\eqref{fixed-nu-gradient-continuity} in $H^{-1}(\mathbb T^d)$ at every time: the two
$H^{-1}(\mathbb T^d)$-continuous processes agree almost everywhere in time.
This completes the proof. 
\end{proof}

\begin{proposition}\label{lem-fixed-nu-energy-limit}
Let $u^\nu$ be the represented viscous solution with the modification
from Proposition~\ref{lem-fixed-nu-path-regularity}.
For every $t\in[0,T]$, the $H^1(\mathbb T^d)$-continuous modification of $u^\nu$
satisfies
\begin{equation}\label{dirichlet-energy-estimate}
\begin{aligned}
 &\frac12\mathbb E\|\nabla u^\nu(t)\|_{L^2(\mathbb T^d)}^2
 +c_{d,\lambda}\mathbb E\int_0^t\!\int_{\mathbb T^d}
 \frac{|D^2u^\nu|^2}{H(\nabla u^\nu)^2}\,dx\,ds
 +\nu\mathbb E\int_0^t\|\Delta u^\nu\|_{L^2(\mathbb T^d)}^2\,ds\\
\leq&\frac12\mathbb E\|\nabla u_0\|_{L^2(\mathbb T^d)}^2
 +\frac{\varepsilon}{2}\|\operatorname{tr}C_K\|_{L^\infty(\mathbb T^d)}
 \mathbb E\int_0^t\mathcal A(u^\nu(s))\,ds.
\end{aligned}
\end{equation}
\end{proposition}

\begin{proof}
By Proposition~\ref{lem-global-galerkin-fixed-nu}, $\tau_R\uparrow T$
almost surely.  Taking expectations in
\eqref{projected-dirichlet-energy}, the expectation of the stopped martingale term is zero.  Fatou's lemma for the terminal gradient and monotone
convergence for the time integrals, as in the proof of that proposition,
give for every $t\in[0,T]$,
\begin{equation}\label{fixed-nu-unstopped-energy-inequality}
\begin{aligned}
 &\frac12\mathbb E\|\nabla u^{N,\nu}(t)\|_{L^2(\mathbb T^d)}^2
 +\nu\mathbb E\int_0^t\|\Delta u^{N,\nu}(s)\|_{L^2(\mathbb T^d)}^2\,ds+c_{d,\lambda}\mathbb E\int_0^t\!\int_{\mathbb T^d}
 \frac{|D^2u^{N,\nu}|^2}{H(\nabla u^{N,\nu})^2}\,dx\,ds\\
 \leq&\frac12\mathbb E\|\nabla u_0^N\|_{L^2(\mathbb T^d)}^2
 +\frac\varepsilon2\|\operatorname{tr}C_K\|_{L^\infty(\mathbb T^d)}
 \mathbb E\int_0^t\mathcal A(u^{N,\nu}(s))\,ds.
\end{aligned}
\end{equation}

We verify each limiting operation in
\eqref{fixed-nu-unstopped-energy-inequality}.  First,
it follows from \eqref{fixed-viscosity-weak-H2} that
\begin{equation}\label{fixed-nu-laplacian-weak-convergence}
 \Delta u^{N,\nu}\rightharpoonup\Delta u^\nu
 \quad\text{weakly in }L^2((0,T)\times\mathbb T^d)
 \quad\text{almost surely}.
\end{equation}
Next, let
$\Phi\in L^2((0,T)\times\mathbb T^d;\mathbb R^{d\times d})$.  Thanks to
\eqref{fixed-viscosity-strong-H1} and the fact that $0<H(\nabla u^{N,\nu})^{-1}\leq1$, dominated convergence theorem yields
\begin{equation*}
 \frac{\Phi}{H(\nabla u^{N,\nu})}
 \longrightarrow\frac{\Phi}{H(\nabla u^\nu)}
 \quad\text{strongly in }L^2((0,T)\times\mathbb T^d).
\end{equation*}
Therefore \eqref{fixed-viscosity-weak-H2} yields
\begin{align*}
 \int_0^T\!\int_{\mathbb T^d}
 \frac{D^2u^{N,\nu}}{H(\nabla u^{N,\nu})}:\Phi\,dx\,ds
 &=\int_0^T\!\int_{\mathbb T^d}
 D^2u^{N,\nu}:\frac{\Phi}{H(\nabla u^{N,\nu})}\,dx\,ds\longrightarrow\int_0^T\!\int_{\mathbb T^d}
 \frac{D^2u^\nu}{H(\nabla u^\nu)}:\Phi\,dx\,ds.
\end{align*}
Hence
\begin{equation}\label{fixed-nu-weighted-hessian-weak-convergence}
 \frac{D^2u^{N,\nu}}{H(\nabla u^{N,\nu})}
 \rightharpoonup\frac{D^2u^\nu}{H(\nabla u^\nu)}
 \quad\text{weakly in }L^2((0,T)\times\mathbb T^d)
 \quad\text{almost surely}.
\end{equation}
Weak lower semicontinuity therefore applies to the weighted-Hessian and
Laplacian terms, also after restriction to $(0,t)\times\mathbb T^d$.

The area functional is $1$-Lipschitz with respect to the $L^1(\mathbb T^d)$ norm of the
gradient because $q\mapsto H(q)$ is $1$-Lipschitz.  Thus
\begin{align}\label{fixed-nu-area-pathwise-convergence}
 \int_0^T|\mathcal A(u^{N,\nu}(s))-\mathcal A(u^\nu(s))|\,ds
 &\leq\int_0^T\!\int_{\mathbb T^d}
 |\nabla u^{N,\nu}-\nabla u^\nu|\,dx\,ds
 \longrightarrow0
\end{align}
almost surely.  Furthermore,
\eqref{fixed-nu-area-second-moment} gives a uniform $L^2(\Omega)$ bound
for their time integrals.  Vitali's theorem
applied to \eqref{fixed-nu-area-pathwise-convergence} yields, for every
$t\in[0,T]$,
\begin{equation}\label{fixed-nu-expected-area-convergence}
 \mathbb E\int_0^t\mathcal A(u^{N,\nu}(s))\,ds
 \longrightarrow\mathbb E\int_0^t\mathcal A(u^\nu(s))\,ds.
\end{equation}
Here $u_0^N$ and $u_0$ denote the represented initial
variables.  By Theorem~\ref{thm-fixed-nu-identification}, each has the law of its
corresponding initial variable on the original space.  Hence
\eqref{fixed-nu-initial-convergence} and the Cauchy--Schwarz inequality give, on the
representation space,
\begin{equation}\label{fixed-nu-represented-initial-energy-convergence}
 \mathbb E\|\nabla u_0^N\|_{L^2(\mathbb T^d)}^2
 \longrightarrow\mathbb E\|\nabla u_0\|_{L^2(\mathbb T^d)}^2.
\end{equation}

Finally, \eqref{fixed-viscosity-strong-path} and lower semicontinuity
of the Dirichlet energy give
$$\|\nabla u^\nu(t)\|_{L^2(\mathbb T^d)}^2
\leq\liminf_{N\to\infty}\|\nabla u^{N,\nu}(t)\|_{L^2(\mathbb T^d)}^2$$
almost surely for every $t\in[0,T]$.
Fatou's lemma and the weak lower semicontinuity of the dissipation terms
allow passage to the limit on the left-hand side of
\eqref{fixed-nu-unstopped-energy-inequality}.
Together with \eqref{fixed-nu-expected-area-convergence} and
\eqref{fixed-nu-represented-initial-energy-convergence} on the right-hand
side, this proves \eqref{dirichlet-energy-estimate} for every $t\in[0,T]$.
\end{proof}

We can now complete the proof of the existence theorem for the viscous
equation.

\begin{proof}[Proof of Theorem~\ref{thm-viscous-equation}]
Take $u_0^N=\Pi_Nu_0$, so
$C_0\leq\mathbb E\|u_0\|_{H^1(\mathbb T^d)}^2$.
Theorem~\ref{thm-fixed-nu-identification} provides a martingale solution
with initial law $\mathcal L(u_0)$, and
Proposition~\ref{lem-fixed-nu-path-regularity} gives the regularity
\eqref{fixed-viscosity-regularity} and $H^1(\mathbb T^d)$ path continuity.
Using $\mathcal A(z)\lesssim1+\|\nabla z\|_{L^2(\mathbb T^d)}$
in \eqref{dirichlet-energy-estimate}, we obtain
\begin{equation*}
 \begin{aligned}
 &\mathbb E\|\nabla u^\nu(t)\|_{L^2(\mathbb T^d)}^2
 +2c_{d,\lambda}\mathbb E\int_0^t\!\int_{\mathbb T^d}
 \frac{|D^2u^\nu|^2}{H(\nabla u^\nu)^2}\,dx\,ds+2\nu\mathbb E\int_0^t\|\Delta u^\nu\|_{L^2(\mathbb T^d)}^2\,ds\\
\leq&\mathbb E\|\nabla u_0\|_{L^2(\mathbb T^d)}^2
 +C\int_0^t\bigl(1+\mathbb E\|\nabla u^\nu(s)\|_{L^2(\mathbb T^d)}^2\bigr)\,ds,
 \end{aligned}
\end{equation*}
where $C$ is independent of $\nu$.  Applying Gronwall's inequality, we obtain \eqref{fixed-nu-solution-estimate}. Moreover, \eqref{fixed-nu-solution-nonuniform-estimate} follows from
\eqref{fixed-nu-limit-path-moment} and
\eqref{fixed-nu-limit-H2-moment}. This completes the proof. \end{proof}

\section{Vanishing-viscosity limit}
\label{sec-vanishing-viscosity}

We prove Theorem~\ref{thm-target-existence} by letting $\nu\downarrow0$
in the solutions from Theorem~\ref{thm-viscous-equation}. 

\subsection{Estimates uniform in viscosity}
\label{subsec-uniform-viscosity-estimates}

We first establish a rigorous area estimate as introduced in
Subsection~\ref{subsec-area-spatial}. 

\begin{proposition}\label{prop-viscous-area-identity}
Fix $0<\nu\leq1$, and let $u^\nu$ be a solution supplied by
Theorem~\ref{thm-viscous-equation} with initial datum $u_0$.
Almost surely, for all $t\in[0,T]$,
\begin{equation}\label{viscous-area-identity}
\begin{aligned}
 &\mathcal A(u^\nu(t))
 +\int_0^t\!\int_{\mathbb T^d}
 H(\nabla u^\nu)|\kappa(u^\nu)|^2\,dx\,ds+\nu\int_0^t\!\int_{\mathbb T^d}
 \left(\frac{|D^2u^\nu|^2}{H(\nabla u^\nu)}
 -\frac{|D^2u^\nu\nabla u^\nu|^2}{H(\nabla u^\nu)^3}\right)dx\,ds\\
 =&\mathcal A(u_0)+M_{\mathcal A,t}^\nu
 +\lambda\int_0^t\!\int_{\mathbb T^d}\mathcal J_d(u^\nu)\,dx\,ds+\frac{\varepsilon}{2}\int_0^t\!\int_{\mathbb T^d}
 \left(\operatorname{tr}C_K-
 \frac{\nabla u^\nu\cdot C_K\nabla u^\nu}{H(\nabla u^\nu)^2}\right)dx\,ds,
\end{aligned}
\end{equation}
where $\mathcal J_d$ is defined in \eqref{dimension-area-remainder} and
\begin{equation}\label{viscous-area-martingale}
 M_{\mathcal A,t}^\nu:=-\sqrt\varepsilon\sum_{k=1}^K\int_0^t
 \left(\int_{\mathbb T^d}f_kH(\nabla u^\nu)^{1/2}
 \kappa(u^\nu)\,dx\right)d\beta_s^k
\end{equation}
is a continuous square-integrable martingale.
Recall that $a_{d,\lambda}$ is defined in
\eqref{area-spatial-coercivity}.  Moreover,
\begin{equation}\label{coercive-area-estimate}
 \begin{aligned}
 \mathbb E\sup_{0\leq t\leq T}\mathcal A(u^\nu(t))
 +\frac{a_{d,\lambda}}2\mathbb E\int_0^T\!\int_{\mathbb T^d}
 H(\nabla u^\nu)|\kappa(u^\nu)|^2\,dx\,dt\leq C_T\left(1+\mathbb E\|u_0\|_{H^1(\mathbb T^d)}^2\right).
 \end{aligned}
\end{equation}
\end{proposition}

\begin{proof}
By Lemma~\ref{lem-fixed-nu-coefficient-regularity} and
Proposition~\ref{lem-fixed-nu-path-regularity},
the solution has $H^2(\mathbb T^d)$ regularity in mean square and almost surely
continuous $H^1(\mathbb T^d)$ paths. More precisely, the drift is square-integrable in $L^2(\mathbb T^d)$, and
the predictable noise coefficients are square-integrable in $H^1(\mathbb T^d)$.
The gradient equation \eqref{fixed-nu-gradient-variational-equation} reads,
in $H^{-1}(\mathbb T^d;\mathbb R^d)$,
\begin{equation}\label{area-gradient-equation}
 \nabla u^\nu(t)
 =\nabla u_0+\int_0^t\nabla b^\nu(u^\nu(s))\,ds+\sqrt\varepsilon\sum_{k=1}^K
 \int_0^t\nabla G_k(u^\nu(s))\,d\beta_s^k.
\end{equation}
By the Sobolev chain rule, we have
$\nabla u^\nu/H(\nabla u^\nu)\in
L^2(\Omega\times(0,T);H^1(\mathbb T^d))$.
By \eqref{matrix-curvature-bounds},
\begin{equation*}
 \mathbb E\int_0^T\!\int_{\mathbb T^d}
 H(\nabla u^\nu)|\kappa(u^\nu)|^2\,dx\,ds
 \leq d\mathbb E\int_0^T\|D^2u^\nu\|_{L^2(\mathbb T^d)}^2\,ds<\infty.
\end{equation*}
This fixed-$\nu$ bound justifies the area martingale before any uniform
curvature estimate is derived.

We adapt the convolution proof of
\cite[Proposition~A.1]{DebusscheHofmanovaVovelle2016}, replacing its
scalar It\^o formula by the finite-dimensional formula for
$H:\mathbb R^d\to\mathbb R$; compare also
\cite[proof of Proposition~5.1]{HofmanovaRogerVonRenesse2017}.
Let $\rho_\delta$ be a smooth periodic approximate identity and
$u_\delta^\nu=\rho_\delta*u^\nu$.
Taking convolution in \eqref{area-gradient-equation}, then almost surely for all $t\in[0,T]$, 
\begin{equation}\label{mollified-gradient-area-Ito-formula}
\begin{aligned}
 \mathcal A(u_\delta^\nu(t))
 ={}&\mathcal A(\rho_\delta*u_0)
 -\int_0^t\!\int_{\mathbb T^d}
 \kappa(u_\delta^\nu(s))
 (\rho_\delta*b^\nu(u^\nu(s)))\,dx\,ds\\
 &+\sqrt\varepsilon\sum_{k=1}^K\int_0^t
 \left(\int_{\mathbb T^d}
 DH(\nabla u_\delta^\nu(s))\cdot
 (\rho_\delta*\nabla G_k(u^\nu(s)))\,dx\right)d\beta_s^k\\
 &+\frac{\varepsilon}{2}\sum_{k=1}^K
 \int_0^t\!\int_{\mathbb T^d}
 (\rho_\delta*\nabla G_k(u^\nu(s)))\cdot
 D^2H(\nabla u_\delta^\nu(s))
 (\rho_\delta*\nabla G_k(u^\nu(s)))\,dx\,ds.
\end{aligned}
\end{equation}
The moment bounds already established at fixed $\nu$ give
\begin{align*}
 u_\delta^\nu&\xrightarrow{\delta\downarrow0}
 u^\nu
 &&\text{in }L^2(\Omega\times(0,T);H^2(\mathbb T^d)),\\
 \rho_\delta*b^\nu(u^\nu)&\xrightarrow{\delta\downarrow0}
 b^\nu(u^\nu)
 &&\text{in }L^2(\Omega\times(0,T)\times\mathbb T^d),\\
 \rho_\delta*\nabla G_k(u^\nu)&\xrightarrow{\delta\downarrow0}
 \nabla G_k(u^\nu)
 &&\text{in }L^2(\Omega\times(0,T)\times\mathbb T^d).
\end{align*}
The argument of Lemma~\ref{lem-quadratic-hessian-continuity},
applied on the product space, gives
$\kappa(u_\delta^\nu)\to\kappa(u^\nu)$ strongly in
$L^2(\Omega\times(0,T)\times\mathbb T^d)$.
Consequently,
\begin{align*}
 \kappa(u_\delta^\nu)(\rho_\delta*b^\nu(u^\nu))
 &\xrightarrow{\delta\downarrow0}
 \kappa(u^\nu)b^\nu(u^\nu)
 \quad\text{in }L^1(\Omega\times(0,T)\times\mathbb T^d).
\end{align*}
The first convergence above implies
$\nabla u_\delta^\nu\to\nabla u^\nu$ strongly in
$L^2(\Omega\times(0,T)\times\mathbb T^d)$.
Together with the last strong $L^2(\Omega\times(0,T)\times\mathbb T^d)$ convergence and the bounds
$|DH|\leq1$ and $|D^2H|\leq\sqrt d$, this gives
\begin{align*}
 DH(\nabla u_\delta^\nu)\cdot
  (\rho_\delta*\nabla G_k(u^\nu))&\xrightarrow{\delta\downarrow0}
  DH(\nabla u^\nu)\cdot\nabla G_k(u^\nu),\\
 (\rho_\delta*\nabla G_k(u^\nu))\cdot
  D^2H(\nabla u_\delta^\nu)(\rho_\delta*\nabla G_k(u^\nu))&\xrightarrow{\delta\downarrow0}
  \nabla G_k(u^\nu)\cdot D^2H(\nabla u^\nu)\nabla G_k(u^\nu), 
\end{align*}
in $L^2(\Omega\times(0,T)\times\mathbb T^d)$ and $L^1(\Omega\times(0,T)\times\mathbb T^d)$, respectively. The It\^o isometry and Doob's inequality therefore give convergence of the
stochastic integration terms in $L^2(\Omega;C([0,T]))$.
The quadratic variation terms converge in $L^1(\Omega;C([0,T]))$.
Pathwise $H^1(\mathbb T^d)$ continuity and boundedness of convolution on $H^1(\mathbb T^d)$
also give $u_\delta^\nu\to u^\nu$ uniformly in time in $H^1(\mathbb T^d)$
almost surely, hence the same convergence for their areas. Consequently, almost surely for every $t\in[0,T]$,
\begin{equation}\label{gradient-area-Ito-formula}
  \begin{aligned}
 \mathcal A(u^\nu(t))
 ={}&\mathcal A(u_0)
 -\int_0^t\!\int_{\mathbb T^d}b^\nu(u^\nu)\kappa(u^\nu)\,dx\,ds
 +M_{\mathcal A,t}^\nu\\
 &+\frac{\varepsilon}{2}\sum_{k=1}^K
 \int_0^t\!\int_{\mathbb T^d}
 \nabla G_k(u^\nu)\cdot D^2H(\nabla u^\nu)\nabla G_k(u^\nu)\,dx\,ds,
  \end{aligned}
\end{equation}
with $M_{\mathcal A}^\nu$ as in \eqref{viscous-area-martingale}.
Using the Cauchy--Schwarz inequality in space, we obtain the bound needed below:
\begin{align}\label{area-martingale-bracket}
 d\langle M_{\mathcal A}^\nu\rangle_t
 =\varepsilon\sum_{k=1}^K\left|
 \int_{\mathbb T^d}f_kH(\nabla u^\nu)^{1/2}\kappa(u^\nu)\,dx
 \right|^2dt\leq\varepsilon\left(\sum_{k=1}^K\|f_k\|_{L^2(\mathbb T^d)}^2\right)
 \int_{\mathbb T^d}H(\nabla u^\nu)|\kappa(u^\nu)|^2\,dx\,dt.
\end{align}
By the fixed-$\nu$ integrability established at the start,
$\mathbb E\langle M_{\mathcal A}^\nu\rangle_T<\infty$, therefore
 $M_{\mathcal A}^\nu$ is a continuous square-integrable martingale.

For almost every $(\omega,s)$, $u^\nu(s)\in H^2(\mathbb T^d)$.  By \eqref{fixed-nu-drift-definition},
\begin{align*}
 -\int_{\mathbb T^d}b^\nu(u^\nu)\kappa(u^\nu)\,dx
 ={}&-\nu\int_{\mathbb T^d}\kappa(u^\nu)\Delta u^\nu\,dx
 -\int_{\mathbb T^d}H(\nabla u^\nu)|\kappa(u^\nu)|^2\,dx\\
 &-\lambda\int_{\mathbb T^d}\kappa(u^\nu)
 \frac{\nabla u^\nu\cdot D^2u^\nu\nabla u^\nu}
 {H(\nabla u^\nu)^3}\,dx.
\end{align*}
Since $\lambda=\varepsilon\Lambda_K/8$,
\eqref{simplified-area-remainder} gives
\begin{align*}
 &\frac{\varepsilon}{2}\sum_{k=1}^K\int_{\mathbb T^d}
 \nabla G_k(u^\nu)\cdot D^2H(\nabla u^\nu)
 \nabla G_k(u^\nu)\,dx-\lambda\int_{\mathbb T^d}\kappa(u^\nu)
 \frac{\nabla u^\nu\cdot D^2u^\nu\nabla u^\nu}
 {H(\nabla u^\nu)^3}\,dx\\
 &=\frac{\varepsilon}{2}\int_{\mathbb T^d}
 \left(\operatorname{tr}C_K-
 \frac{\nabla u^\nu\cdot C_K\nabla u^\nu}
 {H(\nabla u^\nu)^2}\right)dx
 +\lambda\int_{\mathbb T^d}\mathcal J_d(u^\nu)\,dx.
\end{align*}
Moreover, \eqref{viscosity-area-compatibility} yields
\begin{align*}
 \nu\int_{\mathbb T^d}\kappa(u^\nu)\Delta u^\nu\,dx
 &=\nu\int_{\mathbb T^d}\left(
 \frac{|D^2u^\nu|^2}{H(\nabla u^\nu)}
 -\frac{|D^2u^\nu\nabla u^\nu|^2}
 {H(\nabla u^\nu)^3}\right)dx\geq0.
\end{align*}
Integrating these identities in $s$ and substituting into
\eqref{gradient-area-Ito-formula}, we obtain
\begin{align*}
 &\mathcal A(u^\nu(t))
 +\int_0^t\!\int_{\mathbb T^d}
 H(\nabla u^\nu)|\kappa(u^\nu)|^2\,dx\,ds+\nu\int_0^t\!\int_{\mathbb T^d}
 \left(\frac{|D^2u^\nu|^2}{H(\nabla u^\nu)}
 -\frac{|D^2u^\nu\nabla u^\nu|^2}
 {H(\nabla u^\nu)^3}\right)dx\,ds\\
 &=\mathcal A(u_0)+M_{\mathcal A,t}^\nu
 +\lambda\int_0^t\!\int_{\mathbb T^d}\mathcal J_d(u^\nu)\,dx\,ds+\frac{\varepsilon}{2}\int_0^t\!\int_{\mathbb T^d}
 \left(\operatorname{tr}C_K-
 \frac{\nabla u^\nu\cdot C_K\nabla u^\nu}
 {H(\nabla u^\nu)^2}\right)dx\,ds,
\end{align*}
which is \eqref{viscous-area-identity}.  We remark that the above terms are integrable, due to the fact that 
\eqref{matrix-curvature-bounds} and
$\sup_{r\geq1}(\log r)/r^2<\infty$ give
\begin{align*}
 &H(\nabla u^\nu)|\kappa(u^\nu)|^2
 +|\mathcal J_d(u^\nu)|+\left|\frac{|D^2u^\nu|^2}{H(\nabla u^\nu)}
 -\frac{|D^2u^\nu\nabla u^\nu|^2}
 {H(\nabla u^\nu)^3}\right|
 \leq C|D^2u^\nu|^2,
\end{align*}
and \eqref{area-covariance-remainder-bound} controls the covariance
term.  

Lemma~\ref{lem-area-spatial-coercivity}, applied almost everywhere,
allows us to move $\lambda\mathcal J_d$ to the left and discard the
nonnegative viscosity term.  Thus
\begin{align*}
 &\mathcal A(u^\nu(t))
 +a_{d,\lambda}\int_0^t\!\int_{\mathbb T^d}
 H(\nabla u^\nu)|\kappa(u^\nu)|^2\,dx\,ds\leq\mathcal A(u_0)+M_{\mathcal A,t}^\nu
 +\frac{\varepsilon t}{2}\int_{\mathbb T^d}\operatorname{tr}C_K\,dx.
\end{align*}
Taking the time supremum, we estimate the martingale by
\eqref{area-martingale-bracket} and the Burkholder--Davis--Gundy inequality 
and Young inequality:
\begin{align*}
 2\mathbb E\sup_{0\leq t\leq T}|M_{\mathcal A,t}^\nu|
 &\leq\frac{a_{d,\lambda}}2\mathbb E\int_0^T\!\int_{\mathbb T^d}
 H(\nabla u^\nu)|\kappa(u^\nu)|^2\,dx\,ds+C_T.
\end{align*}
Combining this with the preceding inequality at $t=T$ and absorbing yields
\eqref{coercive-area-estimate}, since
$\mathcal A(u_0)\lesssim1+\|\nabla u_0\|_{L^2(\mathbb T^d)}$.
\end{proof}

We need the following second-moment bounds for products of the area
and curvature dissipation in the uniform path estimates and martingale
limit.

\begin{lemma}\label{lem-second-moment-area}
Let $u^\nu$ be as in Proposition~\ref{prop-viscous-area-identity}.
Then
\begin{equation}\label{second-moment-area-estimate}
 \begin{aligned}
 &\mathbb E\left[\sup_{0\leq t\leq T}\mathcal A(u^\nu(t))^2\right]
 +\mathbb E\left[\left(
 \int_0^T\!\int_{\mathbb T^d}H(\nabla u^\nu)
 |\kappa(u^\nu)|^2\,dx\,dt\right)^2\right]\leq C_T\left(1+\mathbb E\|u_0\|_{H^1(\mathbb T^d)}^2\right), 
 \end{aligned}
\end{equation}
with a constant independent of $\nu\in(0,1]$.
\end{lemma}

\begin{proof}
Using the Cauchy--Schwarz inequality to see that 
\begin{equation*}
 \mathcal A(u_0)^2
 \lesssim1+\|\nabla u_0\|_{L^2(\mathbb T^d)}^2,
\end{equation*}
so the initial area has a finite second moment.  The first-moment
estimate \eqref{coercive-area-estimate}, the bracket bound
\eqref{area-martingale-bracket}, and Doob's inequality give
\begin{align*}
 \mathbb E\sup_{0\leq t\leq T}|M_{\mathcal A,t}^\nu|^2
 &\leq4\mathbb E\langle M_{\mathcal A}^\nu\rangle_T
 \leq C\mathbb E\int_0^T\!\int_{\mathbb T^d}
 H(\nabla u^\nu)|\kappa(u^\nu)|^2\,dx\,dt\leq C_T\bigl(1+\mathbb E\|u_0\|_{H^1(\mathbb T^d)}^2\bigr).
\end{align*}
Taking the time supremum in the pathwise area inequality from the proof
of Proposition~\ref{prop-viscous-area-identity}, squaring, and taking
expectations, this martingale bound and the initial-area estimate give
\eqref{second-moment-area-estimate}.
\end{proof}

In the following, we collect the Dirichlet, area, and noise estimates. 

\begin{proposition}
\label{prop-uniform-geometric-estimates-plan}
Suppose Assumption~\ref{constant-noise-covariance} holds.  For each $\nu\in(0,1]$, let
$u^\nu$ be a solution from
Theorem~\ref{thm-viscous-equation}, and suppose that
\begin{equation*}
 \mathbb E\|u_0\|_{H^1(\mathbb T^d)}^2<\infty.
\end{equation*}
Then, for every $0<\nu\leq1$,
\begin{equation}\label{target-gradient-curvature-estimate}
  \begin{aligned}
 &\sup_{0\leq t\leq T}\mathbb E
 \|\nabla u^\nu(t)\|_{L^2(\mathbb T^d)}^2
 +\mathbb E\int_0^T\!\int_{\mathbb T^d}
 \frac{|D^2u^\nu|^2}{H(\nabla u^\nu)^2}\,dx\,dt\\
 &+\mathbb E\int_0^T\!\int_{\mathbb T^d}
 H(\nabla u^\nu)|\kappa(u^\nu)|^2\,dx\,dt
 +\nu\mathbb E\int_0^T\|\Delta u^\nu\|_{L^2(\mathbb T^d)}^2\,dt\\
 &\leq C_T\bigl(1+\mathbb E\|u_0\|_{H^1(\mathbb T^d)}^2\bigr).
\end{aligned}
\end{equation}
The spatial mean estimate also yields
\begin{equation}\label{uniform-u-estimate}
 \sup_{0\leq t\leq T}\mathbb E\|u^\nu(t)\|_{L^2(\mathbb T^d)}^2
 \leq C_T\bigl(1+\mathbb E\|u_0\|_{H^1(\mathbb T^d)}^2\bigr).
\end{equation}
For every $\sigma>d/2$, the full drift satisfies
\begin{equation}\label{uniform-drift-negative-bound}
 \mathbb E\int_0^T\|b^\nu(u^\nu(t))\|_{H^{-\sigma}(\mathbb T^d)}^2\,dt
 \leq C_T\bigl(1+\mathbb E\|u_0\|_{H^1(\mathbb T^d)}^2\bigr).
\end{equation}
In particular,
\begin{equation}\label{uniform-viscosity-vanishing-bound}
 \mathbb E\int_0^T\|\nu\Delta u^\nu(t)\|_{L^2(\mathbb T^d)}^2\,dt
 \leq C_T\nu\bigl(1+\mathbb E\|u_0\|_{H^1(\mathbb T^d)}^2\bigr).
\end{equation}
The constants are independent of $\nu$.
\end{proposition}

\begin{proof}
The first, second, and fourth terms of
\eqref{target-gradient-curvature-estimate} can be estimated by \eqref{fixed-nu-solution-estimate}, and the third term of \eqref{target-gradient-curvature-estimate} can be estimated by \eqref{coercive-area-estimate}.
Both can be bounded by the right-hand side
$C_T(1+\mathbb E\|u_0\|_{H^1(\mathbb T^d)}^2)$.
The Laplacian estimate also yields
\begin{align*}
 \mathbb E\int_0^T\|\nu\Delta u^\nu\|_{L^2(\mathbb T^d)}^2\,dt
 &=\nu\left(\nu\mathbb E\int_0^T
 \|\Delta u^\nu\|_{L^2(\mathbb T^d)}^2\,dt\right)\\
 &\leq C_T\nu\bigl(1+\mathbb E\|u_0\|_{H^1(\mathbb T^d)}^2\bigr)
 \leq C_T\bigl(1+\mathbb E\|u_0\|_{H^1(\mathbb T^d)}^2\bigr),
\end{align*}
which proves \eqref{uniform-viscosity-vanishing-bound}.  This uses the
$C_T$ estimate in \eqref{fixed-nu-solution-estimate}, not its
$C_{\nu,T}$ bound \eqref{fixed-nu-solution-nonuniform-estimate}
for the unweighted $H^2(\mathbb T^d)$-norm.

To bound the drift, fix $\sigma>d/2$.
The continuous
embedding $L^1(\mathbb T^d)\hookrightarrow H^{-\sigma}(\mathbb T^d)$ and
the Cauchy--Schwarz inequality in space give
\begin{align*}
 \|H(\nabla u^\nu)\kappa(u^\nu)\|_{H^{-\sigma}(\mathbb T^d)}^2
 &\leq C_\sigma\left(\int_{\mathbb T^d}
 |H(\nabla u^\nu)\kappa(u^\nu)|\,dx\right)^2\leq C_\sigma\mathcal A(u^\nu)
 \int_{\mathbb T^d}H(\nabla u^\nu)|\kappa(u^\nu)|^2\,dx.
\end{align*}
Consequently, Lemma~\ref{lem-second-moment-area} implies
\begin{align*}
 &\mathbb E\int_0^T
 \|H(\nabla u^\nu)\kappa(u^\nu)\|_{H^{-\sigma}(\mathbb T^d)}^2\,dt\\
 \leq& C_\sigma
 \left(\mathbb E\sup_{0\leq t\leq T}\mathcal A(u^\nu(t))^2\right)^{1/2}
 \left[\mathbb E\left(\int_0^T\!\int_{\mathbb T^d}
 H(\nabla u^\nu)|\kappa(u^\nu)|^2\,dx\,dt\right)^2\right]^{1/2}
 \\
 \leq& C_T\bigl(1+\mathbb E\|u_0\|_{H^1(\mathbb T^d)}^2\bigr).
\end{align*}
The remaining drift terms satisfy
\begin{align*}
 &\mathbb E\int_0^T\left\|
 \frac{\nabla u^\nu\cdot D^2u^\nu\nabla u^\nu}
 {H(\nabla u^\nu)^3}\right\|_{H^{-\sigma}(\mathbb T^d)}^2dt
 \leq C_\sigma\mathbb E\int_0^T\!\int_{\mathbb T^d}
 \frac{|D^2u^\nu|^2}{H(\nabla u^\nu)^2}\,dx\,dt
 \leq C_T\bigl(1+\mathbb E\|u_0\|_{H^1(\mathbb T^d)}^2\bigr),\\
 &\mathbb E\int_0^T\|\nu\Delta u^\nu\|_{H^{-\sigma}(\mathbb T^d)}^2dt
 \leq\nu^2\mathbb E\int_0^T\|\Delta u^\nu\|_{L^2(\mathbb T^d)}^2dt
 \leq C_T\nu\bigl(1+\mathbb E\|u_0\|_{H^1(\mathbb T^d)}^2\bigr)
 \leq C_T\bigl(1+\mathbb E\|u_0\|_{H^1(\mathbb T^d)}^2\bigr).
\end{align*}
Here the first inequality uses
\begin{equation*}
 \left|\frac{\nabla u^\nu\cdot D^2u^\nu\nabla u^\nu}
 {H(\nabla u^\nu)^3}\right|
 \leq\frac{|D^2u^\nu|}{H(\nabla u^\nu)},
\end{equation*}
and the last inequalities follow from
\eqref{target-gradient-curvature-estimate} and
\eqref{uniform-viscosity-vanishing-bound}, proved above.
Summing the three terms in \eqref{fixed-nu-drift-definition} now proves
\eqref{uniform-drift-negative-bound}.

To prove \eqref{uniform-u-estimate}, we
adapt the mean-value argument from
\cite[Proof of Corollary~5.2]{HofmanovaRogerVonRenesse2017}.  Testing the
$L^2(\mathbb T^d)$-valued identity of Proposition~\ref{lem-fixed-nu-path-regularity}
against the constant function $1$ gives
\begin{equation}\label{uniform-spatial-mean-equation}
\begin{aligned}
 \int_{\mathbb T^d}u^\nu(t)\,dx
 &=\int_{\mathbb T^d}u_0\,dx+\int_0^t\langle b^\nu(u^\nu(s)),1\rangle_{H^{-\sigma}(\mathbb T^d),H^\sigma(\mathbb T^d)}\,ds\\
 &\quad+\sqrt\varepsilon\sum_{k=1}^K\int_0^t
 \langle G_k(u^\nu(s)),1\rangle\,d\beta_s^k.
\end{aligned}
\end{equation}
The duality estimate and \eqref{uniform-drift-negative-bound} give
\begin{align*}
 \mathbb E\sup_{0\leq t\leq T}
 \left|\int_0^t
 \langle b^\nu(u^\nu(s)),1\rangle_{H^{-\sigma}(\mathbb T^d),H^\sigma(\mathbb T^d)}\,ds
 \right|^2&\leq T\|1\|_{H^\sigma(\mathbb T^d)}^2\mathbb E\int_0^T
 \|b^\nu(u^\nu(s))\|_{H^{-\sigma}(\mathbb T^d)}^2\,ds\\
 &\leq C_T\bigl(1+\mathbb E\|u_0\|_{H^1(\mathbb T^d)}^2\bigr).
\end{align*}
Doob's inequality, $\sum_k\|G_k(u^\nu)\|_{L^2(\mathbb T^d)}^2
=\Lambda_K\mathcal A(u^\nu)$, and \eqref{coercive-area-estimate} yield
\begin{align*}
 &\mathbb E\sup_{0\leq t\leq T}\left|
 \sqrt\varepsilon\sum_{k=1}^K\int_0^t
 \langle G_k(u^\nu(s)),1\rangle\,d\beta_s^k\right|^2\\
 &\quad\leq4\varepsilon\|1\|_{L^2(\mathbb T^d)}^2\Lambda_K
 \mathbb E\int_0^T\mathcal A(u^\nu(s))\,ds
 \leq C_T\bigl(1+\mathbb E\|u_0\|_{H^1(\mathbb T^d)}^2\bigr).
\end{align*}
Moreover,
\begin{equation*}
 \mathbb E\left|\int_{\mathbb T^d}u_0\,dx\right|^2
 \lesssim\mathbb E\|u_0\|_{L^2(\mathbb T^d)}^2.
\end{equation*}
Combining these estimates in \eqref{uniform-spatial-mean-equation}, we obtain
\begin{equation}\label{uniform-spatial-mean-bound}
 \mathbb E\sup_{0\leq t\leq T}
 \left|\int_{\mathbb T^d}u^\nu(t)\,dx\right|^2
 \leq C_T\bigl(1+\mathbb E\|u_0\|_{H^1(\mathbb T^d)}^2\bigr).
\end{equation}
Finally, the Poincar\'e inequality gives, for every
$z\in H^1(\mathbb T^d)$,
\begin{equation*}
 \|z\|_{L^2(\mathbb T^d)}^2
 \lesssim\|\nabla z\|_{L^2(\mathbb T^d)}^2
 +\left|\int_{\mathbb T^d}z\,dx\right|^2.
\end{equation*}
Applying this inequality to $u^\nu(t)$ and using
\eqref{target-gradient-curvature-estimate} and
\eqref{uniform-spatial-mean-bound} proves \eqref{uniform-u-estimate}.
\end{proof}

We next provide an estimate of H\"older-norm in time. 

\begin{lemma}\label{lem-uniform-viscosity-holder}
Under the hypotheses of Proposition~\ref{prop-uniform-geometric-estimates-plan},
with $m=3+d/2$, for every $d/2<\sigma<m$ and every $0<\gamma<1/4$,
\begin{equation}\label{uniform-viscosity-holder-estimate}
 \mathbb E
 \|u^\nu\|_{C^\gamma([0,T];H^{-\sigma}(\mathbb T^d))}^2
 +\mathbb E
 \|u^\nu\|_{C^\gamma([0,T];H^{-m}(\mathbb T^d))}^2
 \leq C_T\bigl(1+\mathbb E\|u_0\|_{H^1(\mathbb T^d)}^2\bigr),
 \qquad 0<\nu\leq1.
\end{equation}
\end{lemma}

\begin{proof}
Use the $H^{-\sigma}(\mathbb T^d)$-valued identity
\begin{equation*}
 u^\nu(t)=u_0+
 \int_0^t b^\nu(u^\nu(r))\,dr
 +\sqrt\varepsilon\sum_{k=1}^K\int_0^tG_k(u^\nu(r))\,d\beta_r^k
\end{equation*}
from Proposition~\ref{lem-fixed-nu-path-regularity}.  The drift bound
\eqref{uniform-drift-negative-bound} and the Cauchy--Schwarz inequality give, for
$0\leq s<t\leq T$,
\begin{equation*}
 \left\|\int_s^t b^\nu(u^\nu(r))\,dr\right\|_{H^{-\sigma}(\mathbb T^d)}
 \leq |t-s|^{1/2}
 \left(\int_0^T\|b^\nu(u^\nu(r))\|_{H^{-\sigma}(\mathbb T^d)}^2\,dr\right)^{1/2}.
\end{equation*}
Since the integral vanishes at $t=0$, its full $C^\gamma([0,T];H^{-\sigma}(\mathbb T^d))$ norm satisfies
\begin{align*}
 &\mathbb E\left\|
\int_0^\cdot b^\nu(u^\nu(r))\,dr
 \right\|_{C^\gamma([0,T];H^{-\sigma}(\mathbb T^d))}^2\\
 &\quad\leq 2\bigl(T+T^{1-2\gamma}\bigr)\mathbb E\int_0^T
 \|b^\nu(u^\nu(r))\|_{H^{-\sigma}(\mathbb T^d)}^2\,dr
 \leq C_T\bigl(1+\mathbb E\|u_0\|_{H^1(\mathbb T^d)}^2\bigr).
\end{align*}
For $s<t$, the Hilbert-space Burkholder--Davis--Gundy inequality,
the identity
$\sum_k\|G_k(u^\nu)\|_{L^2(\mathbb T^d)}^2=\Lambda_K\mathcal A(u^\nu)$, and
\eqref{second-moment-area-estimate} imply
\begin{align*}
 \mathbb E\left\|\sqrt\varepsilon\sum_{k=1}^K\int_s^t
 G_k(u^\nu(r))\,d\beta_r^k\right\|_{L^2(\mathbb T^d)}^4
 &\leq C\mathbb E\left(\int_s^t
 \mathcal A(u^\nu(r))\,dr\right)^2\\
 &\leq C|t-s|^2\mathbb E\sup_{0\leq r\leq T}
 \mathcal A(u^\nu(r))^2
 \leq C_T\bigl(1+\mathbb E\|u_0\|_{H^1(\mathbb T^d)}^2\bigr)|t-s|^2.
\end{align*}
Choose $\alpha$ with $\gamma+1/4<\alpha<1/2$.  As in
Proposition~\ref{lem-fixed-nu-holder}, Tonelli's theorem yields
\begin{align*}
 \mathbb E\left[
 \sqrt\varepsilon\sum_{k=1}^K\int_0^\cdot
 G_k(u^\nu(r))\,d\beta_r^k
 \right]_{W^{\alpha,4}(0,T;L^2(\mathbb T^d))}^4\leq C_T\bigl(1+\mathbb E\|u_0\|_{H^1(\mathbb T^d)}^2\bigr)
 \int_0^T\!\int_0^t(t-s)^{1-4\alpha}\,ds\,dt<\infty.
\end{align*}
The $L^4(0,T;L^2(\mathbb T^d))$ part is bounded by
$C_T(1+\mathbb E\|u_0\|_{H^1(\mathbb T^d)}^2)$ using the same increment
estimate with $s=0$, since the stochastic integral vanishes at $t=0$.  The embedding
$W^{\alpha,4}(0,T;L^2(\mathbb T^d))\hookrightarrow C^\gamma([0,T];L^2(\mathbb T^d))$
therefore bounds the full stochastic H\"older norm in $L^4(\Omega)$.
Finally, $H^1(\mathbb T^d)\hookrightarrow H^{-\sigma}(\mathbb T^d)$ and
$\mathbb E\|u_0\|_{H^1(\mathbb T^d)}^2<\infty$ control the constant
initial path.  Adding the three parts proves the first bound in
\eqref{uniform-viscosity-holder-estimate}; the second follows from
$H^{-\sigma}(\mathbb T^d)\hookrightarrow H^{-m}(\mathbb T^d)$.
\end{proof}

\subsection{Compactness and identification of geometric terms}
\label{subsec-vanishing-compactness}

With $m=3+d/2$, define the path spaces
\begin{equation*}
\begin{aligned}
 \mathcal X_u&=L^2(0,T;L^2(\mathbb T^d))
 \cap C([0,T];H^{-m}(\mathbb T^d)),\\
 \mathcal X_{\nabla}&=
 L^2((0,T)\times\mathbb T^d;\mathbb R^d)_{\mathrm w},
 &\mathcal X_n&=
 L^2(0,T;H^1(\mathbb T^d;\mathbb R^d))_{\mathrm w},\\
 \mathcal X_D&=
 L^2((0,T)\times\mathbb T^d;\mathbb R^{d\times d})_{\mathrm w},
 &\mathcal X_\kappa&=
 L^2((0,T)\times\mathbb T^d)_{\mathrm w},\\
 \mathcal X_\beta&=C([0,T];\mathbb R^K),
 &\mathcal X_0&=H^1(\mathbb T^d)_{\mathrm w},
\end{aligned}
\end{equation*}
and set
\begin{equation}\label{joint-tightness-space}
 \mathcal X=\mathcal X_u\times\mathcal X_{\nabla}
 \times\mathcal X_n\times\mathcal X_D\times\mathcal X_\kappa
 \times\mathcal X_\beta\times\mathcal X_0.
\end{equation}
Here the subscript
$\mathrm w$ denotes the weak topology. For each solution $u^\nu$, let $\beta^\nu$ denote its driving Brownian
motion and define the following coordinates
\begin{equation}\label{vanishing-viscosity-spatial-coordinates}
\begin{aligned}
 Y_{\nabla}^\nu&:=\nabla u^\nu,
 &Y_n^\nu&:=\frac{\nabla u^\nu}{H(\nabla u^\nu)},\\
 Y_D^\nu&:=\frac{D^2u^\nu}{H(\nabla u^\nu)},
 &Y_\kappa^\nu&:=H(\nabla u^\nu)^{1/2}\kappa(u^\nu).
\end{aligned}
\end{equation}
Set
\begin{equation}\label{joint-tightness-vector}
 Z^\nu:=\left(u^\nu,Y_{\nabla}^\nu,Y_n^\nu,
 Y_D^\nu,Y_\kappa^\nu,\beta^\nu,u_0\right),
 \qquad \mu^\nu:=\mathcal L(Z^\nu).
\end{equation}
The spatial derivatives are understood almost everywhere in the weak
sense, using the $H^2(\mathbb T^d)$ regularity of the solution for fixed $\nu$.

\begin{proposition}\label{prop-time-tightness-plan}
Under Assumption~\ref{constant-noise-covariance}, let $u^\nu$ be viscous
solutions from Theorem~\ref{thm-viscous-equation} such that
\begin{equation*}
 \mathbb E\|u_0\|_{H^1(\mathbb T^d)}^2<\infty.
\end{equation*}
Then $\{\mu^\nu:0<\nu\leq1\}$ is tight on $\mathcal X$.
\end{proposition}

\begin{proof}
Fix $d/2<\sigma<m$ and $0<\gamma<1/4$.
Equations~\eqref{uniform-u-estimate} and
\eqref{target-gradient-curvature-estimate}, together with
Lemma~\ref{lem-uniform-viscosity-holder}, give
\begin{align*}
 &\sup_{0<\nu\leq1}\mathbb E\left(
 \|u^\nu\|_{L^2(0,T;H^1(\mathbb T^d))}^2
 +\|u^\nu\|_{C^\gamma([0,T];H^{-\sigma}(\mathbb T^d))}^2\right)\leq C_T\left(1+
 \mathbb E\|u_0\|_{H^1(\mathbb T^d)}^2\right)<\infty.
\end{align*}
By the Aubin--Lions--Simon theorem \cite{Simon1987}, the
Arzel\`a--Ascoli theorem, and Markov's inequality, this bound implies
tightness of the laws of $u^\nu$ on $\mathcal X_u$.

For the remaining spatial coordinates, the fixed-viscosity $H^2(\mathbb T^d)$
regularity and Sobolev chain rule give
\begin{equation*}
 \nabla Y_n^\nu=(I-Y_n^\nu\otimes Y_n^\nu)Y_D^\nu.
\end{equation*}
Since $|Y_n^\nu|\leq1$, the Frobenius norm satisfies
\begin{equation*}
 |I-Y_n^\nu\otimes Y_n^\nu|^2
 =d-2|Y_n^\nu|^2+|Y_n^\nu|^4\leq d.
\end{equation*}
Thus $|\nabla Y_n^\nu|
\leq\sqrt d\,|Y_D^\nu|$, and
\eqref{target-gradient-curvature-estimate} yields
\begin{align*}
 &\sup_{0<\nu\leq1}\mathbb E\left(
 \|Y_{\nabla}^\nu\|_{L^2(0,T;L^2(\mathbb T^d))}^2
 +\|Y_n^\nu\|_{L^2(0,T;H^1(\mathbb T^d))}^2
 +\|Y_D^\nu\|_{L^2(0,T;L^2(\mathbb T^d))}^2
 +\|Y_\kappa^\nu\|_{L^2(0,T;L^2(\mathbb T^d))}^2\right)\\
 &\leq C_T\left(1+
 \mathbb E\|u_0\|_{H^1(\mathbb T^d)}^2\right)<\infty.
\end{align*}
The initial
coordinate has the assumed $H^1(\mathbb T^d)$ moment bound, and therefore tight in $H^1(\mathbb T^d)_w$. Each $\beta^\nu$ has the
same tight Wiener law. This completes the proof. 

\end{proof}

In the following, we carry out the Jakubowski-Skorokhod argument again to construct an almost surely convergent representation. 

\begin{proposition}\label{prop-vanishing-viscosity-representation}
Assume the hypotheses of Proposition~\ref{prop-time-tightness-plan}.  There exist a subsequence $(\nu_j)_{j\geq1}$ with $\nu_j\downarrow0$, a probability space $(\widetilde\Omega,\widetilde{\mathcal F},
\widetilde{\mathbb P})$, $\mathcal X$-valued random variables
$\widetilde Z^j$, and an $\mathcal X$-valued random variable
$\widetilde Z$ such that, for a probability measure $\mu$ on $\mathcal X$,
\begin{equation}\label{vanishing-viscosity-representation-convergence}
\begin{aligned}
 \mu^{\nu_j}&\Longrightarrow\mu,\quad
 \mathcal L_{\widetilde{\mathbb P}}(\widetilde Z^j)
 =\mu^{\nu_j},\quad
 \mathcal L_{\widetilde{\mathbb P}}(\widetilde Z)=\mu,\\
 \widetilde Z^j&\longrightarrow\widetilde Z
 \quad\text{in }\mathcal X,
 \qquad\widetilde{\mathbb P}\text{-almost surely}.
\end{aligned}
\end{equation}
\end{proposition}

\begin{proof}
Tightness follows from Proposition~\ref{prop-time-tightness-plan}. Consequently, the
Jakubowski--Skorokhod representation theorem
\cite{Jakubowski1997} proves the result. 
\end{proof}

Write the components of $\widetilde Z^j$ as
$(\widetilde u^{\nu_j},\widetilde Y_{\nabla}^{\nu_j},
\widetilde Y_n^{\nu_j},\widetilde Y_D^{\nu_j},
\widetilde Y_\kappa^{\nu_j},\widetilde\beta^{\nu_j},
\widetilde u_0^j)$, and denote the components of $\widetilde Z$ by
$(\widetilde u,\widetilde Y_{\nabla},\widetilde Y_n,
\widetilde Y_D,\widetilde Y_\kappa,\widetilde\beta,\widetilde u_0)$. We first verify that the spatial coordinates retain their geometric
relations to the viscous solution after representation.

\begin{lemma}\label{lem-vanishing-viscosity-coordinate-identities}
For every $j\geq1$, the represented coordinates furnished by
Proposition~\ref{prop-vanishing-viscosity-representation} satisfy, almost
surely,
\begin{equation}\label{represented-spatial-coordinate-identities}
\begin{aligned}
 \widetilde Y_{\nabla}^{\nu_j}
 &=\nabla\widetilde u^{\nu_j},
 &\widetilde Y_n^{\nu_j}
 &=\frac{\nabla\widetilde u^{\nu_j}}
 {H(\nabla\widetilde u^{\nu_j})},\\
 \widetilde Y_D^{\nu_j}
 &=\frac{D^2\widetilde u^{\nu_j}}
 {H(\nabla\widetilde u^{\nu_j})},
 &\widetilde Y_\kappa^{\nu_j}
 &=H(\nabla\widetilde u^{\nu_j})^{1/2}
 \kappa(\widetilde u^{\nu_j}).
\end{aligned}
\end{equation}
In particular, $\widetilde u^{\nu_j}\in L^2(0,T;H^2(\mathbb T^d))$ and has a
$C([0,T];H^1(\mathbb T^d))$ version almost surely.  The area-path functionals of
this version have the same laws as those of the original solution.
\end{lemma}

\begin{proof}
Set
$$
 \mathcal E:=\{z\in\mathcal X_u:z\in L^2(0,T;H^2(\mathbb T^d))\}.
$$
By the Lusin--Souslin theorem applied to the continuous injection
\begin{equation*}
 L^2(0,T;H^2(\mathbb T^d))\hookrightarrow L^2(0,T;L^2(\mathbb T^d)),
\end{equation*}
the set $\mathcal E$ is Borel in $\mathcal X_u$, and the reconstruction map
\begin{equation}\label{represented-H2-coordinate-map}
 \mathcal E\ni z\longmapsto z\in L^2(0,T;H^2(\mathbb T^d))
\end{equation}
is Borel.  The Sobolev chain rule and the bounded-coefficient
continuity argument of Lemma~\ref{lem-quadratic-hessian-continuity}
show that
\[
 \mathcal C(z):=\left(
 \nabla z,\frac{\nabla z}{H(\nabla z)},
 \frac{D^2z}{H(\nabla z)},H(\nabla z)^{1/2}\kappa(z)\right)
\]
is continuous from $L^2(0,T;H^2(\mathbb T^d))$ into the product of the corresponding
norm Hilbert spaces.  Composing with
\eqref{represented-H2-coordinate-map} and extending by zero outside
$\mathcal E$ makes $\mathcal C$ Borel on $\mathcal X_u$.
It is also Borel with values in the weak Hilbert factors, since
their weak and strong Borel $\sigma$-fields agree.

We now use equality of the joint laws.  For any Borel map $F$ between
the spaces considered here,
\[
 \mathcal L(\bar X,\bar Y)=\mathcal L(X,F(X))
 \quad\Longrightarrow\quad
 \bar Y=F(\bar X)\ \text{almost surely},
\]
because the graph of $F$ is Borel.
On the original space, $u^{\nu_j}\in\mathcal E$ almost surely and the
four spatial coordinates equal $\mathcal C(u^{\nu_j})$.
Equality of the original and represented joint laws implies
$\widetilde u^{\nu_j}\in\mathcal E$ and the same coordinate identities
on the new space, proving \eqref{represented-spatial-coordinate-identities}.

Finally, apply Lusin--Souslin to the embedding
\[
 C([0,T];H^1(\mathbb T^d))\hookrightarrow C([0,T];H^{-m}(\mathbb T^d)).
\]
Its image and inverse are Borel, so membership in $C([0,T];H^1(\mathbb T^d))$
and reconstruction in that space are Borel.
Equality of the solution-coordinate laws therefore preserves the
$H^1(\mathbb T^d)$ path continuity.  Since $\mathcal A:H^1(\mathbb T^d)\to\mathbb R$ is
Lipschitz, the map $z\mapsto\mathcal A(z(\cdot))$ is continuous from
$C([0,T];H^1(\mathbb T^d))$ to $C([0,T])$.
Thus the area paths also retain their laws, and hence their
measurable functionals have the same laws.
\end{proof}

For the rest of this section, we work on the representation space
$(\widetilde\Omega,\widetilde{\mathcal F},\widetilde{\mathbb P})$. 
All subsequent convergences of the represented sequence are taken as
$j\to\infty$. We write $\widetilde{\mathbb E}$ for expectation on this space. 

We first prove strong convergence of the normalized gradient and the curvature. 

\begin{proposition}\label{lem-strong-normal-limit}
Assume the hypotheses of Proposition~\ref{prop-time-tightness-plan}.
Let $\widetilde Z^j\to \widetilde Z$ be the almost surely convergent representation furnished
by Proposition~\ref{prop-vanishing-viscosity-representation}.
The first two spatial limit coordinates satisfy
\begin{equation}\label{gradient-coordinate-identification}
 \widetilde Y_{\nabla}=\nabla \widetilde u,
 \qquad
 \widetilde Y_n=\frac{\nabla \widetilde u}{H(\nabla \widetilde u)}
 \quad\text{almost everywhere}.
\end{equation}
Moreover,
\begin{align*}
 \frac{\nabla \widetilde u}{H(\nabla \widetilde u)}
 &\in L^2\bigl(\widetilde\Omega;L^2(0,T;H^1(\mathbb T^d))\bigr),\\
 \kappa(\widetilde u)=\operatorname{div}
 \left(\frac{\nabla \widetilde u}{H(\nabla \widetilde u)}\right)
 &\in L^2(\widetilde\Omega\times(0,T)\times\mathbb T^d).
\end{align*}
Along the entire represented sequence,
\begin{equation}\label{strong-normal-convergence}
 \widetilde Y_n^{\nu_j}\longrightarrow \widetilde Y_n
 \quad\text{strongly in }L^2((0,T)\times\mathbb T^d)
 \quad\text{almost surely}.
\end{equation}
The same normalized gradients converge strongly in
$L^2(\widetilde\Omega\times(0,T)\times\mathbb T^d)$.
\end{proposition}

\begin{proof}
By Lemma~\ref{lem-vanishing-viscosity-coordinate-identities}, on one common
event of probability one each $\widetilde u^{\nu_j}$ has
$L^2(0,T;H^2(\mathbb T^d))$ regularity and the represented coordinates satisfy
\eqref{represented-spatial-coordinate-identities}.  Since
$\widetilde u^{\nu_j}\to \widetilde u$ strongly in $L^2(0,T;L^2(\mathbb T^d))$ and
$\widetilde Y_{\nabla}^{\nu_j}\rightharpoonup \widetilde Y_{\nabla}$ in $L^2((0,T)\times\mathbb T^d)$, integration by parts
identifies the weak limit:
\begin{equation}\label{gradient-weak-limit}
 \widetilde u\in L^2(0,T;H^1(\mathbb T^d)),\qquad
 \widetilde Y_{\nabla}=\nabla \widetilde u,\qquad
 \widetilde Y_{\nabla}^{\nu_j}\rightharpoonup\nabla \widetilde u
 \quad\text{in }L^2((0,T)\times\mathbb T^d)
\end{equation}
almost surely.  Weak convergence of the gradient coordinate implies
$\sup_{j\geq1}\|\widetilde Y_{\nabla}^{\nu_j}\|_{L^2((0,T)\times\mathbb T^d)}<\infty$.
Similarly, weak $L^2((0,T)\times\mathbb T^d)$ convergence of $\widetilde Y_\kappa^{\nu_j}$ and
$H(\nabla \widetilde u^{\nu_j})\geq1$ give
$$
 \sup_{j\geq1}\|\kappa(\widetilde u^{\nu_j})\|_{L^2(0,T;L^2(\mathbb T^d))}
 \leq\sup_{j\geq1}\|\widetilde Y_\kappa^{\nu_j}\|_{L^2(0,T;L^2(\mathbb T^d))}<\infty.
$$
Thus, on this event,
\begin{equation}\label{pathwise-gradient-curvature-bounds}
 \sup_{j\geq1}\|\nabla \widetilde u^{\nu_j}\|_{L^2(0,T;L^2(\mathbb T^d))}
 +\sup_{j\geq1}\|\kappa(\widetilde u^{\nu_j})\|_{L^2(0,T;L^2(\mathbb T^d))}<\infty.
\end{equation}
Moreover, \eqref{target-gradient-curvature-estimate}, equality of laws,
weak lower semicontinuity and Fatou's lemma imply
\begin{equation}\label{expected-gradient-bounds-represented}
 \sup_{j\geq1}\widetilde{\mathbb E}\|\nabla \widetilde u^{\nu_j}\|_{L^2(0,T;L^2(\mathbb T^d))}^2
 +\widetilde{\mathbb E}\|\nabla \widetilde u\|_{L^2(0,T;L^2(\mathbb T^d))}^2
 \leq C_T\left(1+
 \mathbb E\|u_0\|_{H^1(\mathbb T^d)}^2\right).
\end{equation}

For $H(\xi)=(1+|\xi|^2)^{1/2}$,
\begin{equation*}
 \nabla H(\xi)=\frac{\xi}{H(\xi)},\qquad
 D^2H(\xi)=\frac1{H(\xi)}I-\frac{\xi\otimes\xi}{H(\xi)^3}.
\end{equation*}
The eigenvalues of $D^2H(\xi)$ are $H(\xi)^{-1}$ and
$H(\xi)^{-3}$ (with the former of multiplicity $d-1$), therefore
$0\leq D^2H(\xi)\leq I$ and $D^2H(\xi)^2\leq D^2H(\xi)$.
Here inequalities between symmetric matrices are understood in the
quadratic-form sense: $A\leq B$ means
$v\cdot Av\leq v\cdot Bv$ for every $v\in\mathbb R^d$,
or equivalently that $B-A$ is positive semidefinite.
In particular, the preceding inequalities mean
$$
 0\leq v\cdot D^2H(\xi)v\leq |v|^2,
 \qquad
 |D^2H(\xi)v|^2\leq v\cdot D^2H(\xi)v,
 \qquad v\in\mathbb R^d.
$$
By the fundamental theorem of calculus and Jensen's inequality, we have 
\begin{equation}\label{normal-cocoercivity}
\begin{aligned}
 \left|\frac{\xi}{H(\xi)}-\frac{\eta}{H(\eta)}\right|^2
 &\leq\int_0^1
 |D^2H(\eta+r(\xi-\eta))(\xi-\eta)|^2\,dr\\
 &\leq\int_0^1
 (\xi-\eta)\cdot D^2H(\eta+r(\xi-\eta))
 (\xi-\eta)\,dr\\
 &=\left(\frac{\xi}{H(\xi)}-\frac{\eta}{H(\eta)}\right)
 \cdot(\xi-\eta).
\end{aligned}
\end{equation}
Applying \eqref{normal-cocoercivity} to $\nabla \widetilde u^{\nu_j}$ and $\nabla \widetilde u$,
then by the integration by parts formula, it follows that 
\begin{equation}\label{normal-limit-estimate}
\begin{aligned}
 &\left\|\frac{\nabla \widetilde u^{\nu_j}}{H(\nabla \widetilde u^{\nu_j})}
 -\frac{\nabla \widetilde u}{H(\nabla \widetilde u)}
 \right\|_{L^2((0,T)\times\mathbb T^d)}^2\\
\leq&-\int_0^T\!\int_{\mathbb T^d}
 \kappa(\widetilde u^{\nu_j})(\widetilde u^{\nu_j}-\widetilde u)\,dx\,dt
 -\int_0^T\!\int_{\mathbb T^d}
 \frac{\nabla \widetilde u}{H(\nabla \widetilde u)}
 \cdot(\nabla \widetilde u^{\nu_j}-\nabla \widetilde u)\,dx\,dt.
\end{aligned}
\end{equation}
The first term tends to zero by
\eqref{pathwise-gradient-curvature-bounds} and the strong convergence of
$\widetilde u^{\nu_j}$ in $L^2(0,T;L^2(\mathbb T^d))$, and the second tends to zero by
\eqref{gradient-weak-limit}.  This proves
\eqref{strong-normal-convergence}.  Since the squared norm on the left of
\eqref{normal-limit-estimate} is bounded by a deterministic constant, dominated
convergence also yields
\begin{equation}\label{strong-normal-product-convergence}
 \frac{\nabla \widetilde u^{\nu_j}}{H(\nabla \widetilde u^{\nu_j})}
 \longrightarrow\frac{\nabla \widetilde u}{H(\nabla \widetilde u)}
 \quad\text{strongly in }
 L^2(\widetilde\Omega\times(0,T)\times\mathbb T^d).
\end{equation}
Thus the represented weak limit of the normalized-gradient coordinate is
\begin{equation*}
 \widetilde Y_n=\frac{\nabla \widetilde u}{H(\nabla \widetilde u)}
 \quad\text{in }L^2(0,T;H^1(\mathbb T^d)),
\end{equation*}
which completes \eqref{gradient-coordinate-identification}.
Weak lower semicontinuity in $L^2(0,T;H^1(\mathbb T^d))$, Fatou's lemma, equality of
laws in Proposition~\ref{prop-vanishing-viscosity-representation}, and
the uniform bound in the proof of
Proposition~\ref{prop-time-tightness-plan} show that 
\begin{align*}
 \widetilde{\mathbb E}\left\|\frac{\nabla \widetilde u}{H(\nabla \widetilde u)}
 \right\|_{L^2(0,T;H^1(\mathbb T^d))}^2
 &\leq\liminf_{j\to\infty}
 \widetilde{\mathbb E}\|\widetilde Y_n^{\nu_j}\|_{L^2(0,T;H^1(\mathbb T^d))}^2\leq C_T\left(1+
 \mathbb E\|u_0\|_{H^1(\mathbb T^d)}^2\right)<\infty.
\end{align*}
Since $\kappa(\widetilde u)=\operatorname{div}\widetilde Y_n$ and
$|\operatorname{div}\widetilde Y_n|^2\leq d|\nabla \widetilde Y_n|^2$,
the stated curvature regularity follows.
\end{proof}

\begin{proposition}\label{prop-gradient-compactness-plan}
Under the hypotheses of Proposition~\ref{lem-strong-normal-limit}, along the
entire represented sequence,
\begin{equation}\label{target-strong-gradient-convergence}
 \nabla \widetilde u^{\nu_j}\longrightarrow\nabla \widetilde u
 \quad\text{strongly in }
 L^p(\widetilde\Omega\times(0,T)\times\mathbb T^d)
 \quad\text{for every }1\leq p<2.
\end{equation}
For each such $p$, the same convergence holds almost surely in
$L^p((0,T)\times\mathbb T^d)$.  In addition,
\begin{equation}\label{strong-square-root-mobility}
 H(\nabla \widetilde u^{\nu_j})^{1/2}\longrightarrow H(\nabla \widetilde u)^{1/2}
 \quad\text{strongly in }L^2((0,T)\times\mathbb T^d)
\end{equation}
almost surely and in $L^2(\widetilde\Omega\times(0,T)\times\mathbb T^d)$.
\end{proposition}

\begin{proof}
The map
\begin{equation*}
 \xi\longmapsto\frac{\xi}{H(\xi)}
\end{equation*}
is one-to-one from $\mathbb R^d$ onto the open unit ball, with continuous
inverse $\zeta\mapsto\zeta/\sqrt{1-|\zeta|^2}$.  From
\eqref{gradient-weak-limit}, $\nabla \widetilde u$ is finite almost everywhere and
therefore $|\nabla \widetilde u|/H(\nabla \widetilde u)<1$ almost everywhere.  Thus the inverse
is evaluated inside its domain.  Write $F(\zeta)=\zeta/\sqrt{1-|\zeta|^2}$ for $|\zeta|<1$.
For $0<\rho<1$, differentiation and the mean-value theorem show that 
\begin{align*}
 DF(\zeta)
 &=\frac{I}{\sqrt{1-|\zeta|^2}}
 +\frac{\zeta\otimes\zeta}{(1-|\zeta|^2)^{3/2}},\\
 |F(\zeta)-F(\eta)|
 &\leq L_\rho|\zeta-\eta|,
 \quad |\zeta|,|\eta|\leq1-\rho/2,\qquad
 L_\rho:=\bigl(1-(1-\rho/2)^2\bigr)^{-3/2}.
\end{align*}
On $\{|\widetilde Y_n|\leq1-\rho,\ |\widetilde Y_n^{\nu_j}-\widetilde Y_n|\leq\rho/2\}$, both arguments
of $F$ lie in this smaller ball and hence
\begin{equation*}
 |\nabla \widetilde u^{\nu_j}-\nabla \widetilde u|
 =|F(\widetilde Y_n^{\nu_j})-F(\widetilde Y_n)|
 \leq L_\rho|\widetilde Y_n^{\nu_j}-\widetilde Y_n|.
\end{equation*}
Fix a path in the common probability-one event.  For every $\delta>0$,
Chebyshev's inequality on $(0,T)\times\mathbb T^d$ gives
\begin{align*}
 &\int_0^T\!\int_{\mathbb T^d}
 \mathbf 1_{\{|\nabla \widetilde u^{\nu_j}-\nabla \widetilde u|>\delta\}}\,dx\,ds\\
\leq&\int_0^T\!\int_{\mathbb T^d}
 \mathbf 1_{\{|\widetilde Y_n|>1-\rho\}}\,dx\,ds
 +\left(\frac{4}{\rho^2}
 +\frac{L_\rho^2}{\delta^2}\right)
 \|\widetilde Y_n^{\nu_j}-\widetilde Y_n\|_{L^2((0,T)\times\mathbb T^d)}^2.
\end{align*}
By \eqref{strong-normal-convergence}, the last term tends to zero.
Since $|\widetilde Y_n|<1$ almost everywhere, it follows from the dominated convergence theorem that 
\begin{equation*}
 \limsup_{j\to\infty}\int_0^T\!\int_{\mathbb T^d}
 \mathbf 1_{\{|\nabla \widetilde u^{\nu_j}-\nabla \widetilde u|>\delta\}}\,dx\,ds
 \leq\lim_{\rho\downarrow0}\int_0^T\!\int_{\mathbb T^d}
 \mathbf 1_{\{|\widetilde Y_n|>1-\rho\}}\,dx\,ds=0.
\end{equation*}
The integral on the left is bounded by a deterministic constant for every
path.  The dominated convergence theorem on $\widetilde\Omega$ therefore also implies 
\begin{equation*}
 \widetilde{\mathbb E}\int_0^T\!\int_{\mathbb T^d}
 \mathbf 1_{\{|\nabla \widetilde u^{\nu_j}-\nabla \widetilde u|>\delta\}}\,dx\,ds
 \longrightarrow0.
\end{equation*}
Thus
\begin{equation}\label{gradient-convergence-in-measure}
 \nabla \widetilde u^{\nu_j}\longrightarrow\nabla \widetilde u
 \quad\text{in measure}
\end{equation}
both pathwise on $(0,T)\times\mathbb T^d$ and on
$\widetilde\Omega\times(0,T)\times\mathbb T^d$.

Let $1\leq p<2$.  For every measurable
$E\subset(0,T)\times\mathbb T^d$, H\"older's inequality and
\eqref{pathwise-gradient-curvature-bounds} give, almost surely,
\begin{equation*}
 \sup_{j\geq1}\int_E|\nabla \widetilde u^{\nu_j}-\nabla \widetilde u|^p\,dx\,dt
 \leq
 \left(\sup_{j\geq1}\|\nabla \widetilde u^{\nu_j}-\nabla \widetilde u\|_{L^2((0,T)\times\mathbb T^d)}^2
 \right)^{p/2}|E|^{1-p/2}.
\end{equation*}
Hence the $p$-th powers are uniformly integrable pathwise.  The analogous
estimate on the probability--space--time product follows from
\eqref{expected-gradient-bounds-represented}.  Vitali's theorem and
\eqref{gradient-convergence-in-measure} now prove all assertions in
\eqref{target-strong-gradient-convergence}.  The endpoint $p=2$ is not
obtained because the squared gradients are not known to be uniformly
integrable.

For every $\xi,\eta\in\mathbb R^d$,
\begin{equation*}
 |H(\xi)^{1/2}-H(\eta)^{1/2}|^2
 \leq|H(\xi)-H(\eta)|\leq|\xi-\eta|.
\end{equation*}
The strong convergence of the gradients, almost surely in
$L^1((0,T)\times\mathbb T^d)$ and in
$L^1(\widetilde\Omega\times(0,T)\times\mathbb T^d)$, therefore proves
\eqref{strong-square-root-mobility}, both almost surely and in
$L^2(\widetilde\Omega\times(0,T)\times\mathbb T^d)$.
\end{proof}

\begin{remark}
The $L^p((0,T)\times\mathbb T^d)$ gradient convergence in
\eqref{target-strong-gradient-convergence} has the same range
$1\leq p<2$ as \cite[Proposition~5.8]{HofmanovaRogerVonRenesse2017}
for its two-dimensional model.  In \cite{HofmanovaRogerVonRenesse2017}, the gradient limit is identified
by Young measures and a div--curl argument. 
\end{remark}

We identify the two drift terms separately.  For the correction, we
first identify the weighted Hessian.  For the curvature drift, we use
strong convergence of the square-root mobility.

\begin{proposition}\label{lem-weighted-hessian-correction-limit}
Assume the hypotheses of
Proposition~\ref{prop-gradient-compactness-plan}.  The weak limit $\widetilde Y_D$
is symmetric, belongs to
$L^2(\widetilde\Omega\times(0,T)\times\mathbb T^d;\mathbb R^{d\times d})$, and
satisfies, almost everywhere,
\begin{equation*}
 \widetilde Y_D=\frac{D^2\widetilde u}{H(\nabla \widetilde u)}
 :=\left(I-\frac{\nabla \widetilde u\otimes\nabla \widetilde u}{H(\nabla \widetilde u)^2}\right)^{-1}
 \nabla\left(\frac{\nabla \widetilde u}{H(\nabla \widetilde u)}\right).
\end{equation*}
Along the entire represented sequence, almost surely,
\begin{align}\label{geometric-correction-convergence}
 \frac{\nabla \widetilde u^{\nu_j}\cdot D^2\widetilde u^{\nu_j}\nabla \widetilde u^{\nu_j}}
 {H(\nabla \widetilde u^{\nu_j})^3}
 &\rightharpoonup
 \frac{\nabla \widetilde u\cdot D^2\widetilde u\,\nabla \widetilde u}{H(\nabla \widetilde u)^3}
 &&\text{weakly in }L^2((0,T)\times\mathbb T^d).
\end{align}
The limiting expressions have the meaning of
Definition~\ref{def-energy-martingale-solution}.
\end{proposition}

\begin{proof}
We work on the new probability space from
Proposition~\ref{prop-gradient-compactness-plan}.  By
\eqref{vanishing-viscosity-representation-convergence} and
\eqref{strong-normal-convergence}, $\widetilde Y_n^{\nu_j}\to \widetilde Y_n$ strongly in
$L^2((0,T)\times\mathbb T^d)$ and weakly in
$L^2(0,T;H^1(\mathbb T^d))$, while
$\widetilde Y_D^{\nu_j}\rightharpoonup \widetilde Y_D$ weakly in
$L^2((0,T)\times\mathbb T^d)$.
Using the chain rule of weak derivatives, we obtain
\begin{equation}\label{Rv-identity}
 \nabla\left(\frac{\nabla \widetilde u^{\nu_j}}{H(\nabla \widetilde u^{\nu_j})}\right)
 =\left(I-\frac{\nabla \widetilde u^{\nu_j}\otimes\nabla \widetilde u^{\nu_j}}
 {H(\nabla \widetilde u^{\nu_j})^2}\right)
 \frac{D^2\widetilde u^{\nu_j}}{H(\nabla \widetilde u^{\nu_j})}.
\end{equation}
The matrix in \eqref{Rv-identity} is a contraction on $\mathbb R^d$,
and $|\nabla\widetilde u^{\nu_j}|/H(\nabla\widetilde u^{\nu_j})\leq1$.
Thus \eqref{target-gradient-curvature-estimate}, transferred by
the identity of joint distributions, gives
\begin{align*}
 \sup_{j\geq1}\widetilde{\mathbb E}
 \left\|\frac{\nabla\widetilde u^{\nu_j}}
 {H(\nabla\widetilde u^{\nu_j})}\right\|_
 {L^2(0,T;H^1(\mathbb T^d;\mathbb R^d))}^2&\lesssim1
 +\sup_{j\geq1}\widetilde{\mathbb E}\int_0^T\!\int_{\mathbb T^d}
 \frac{|D^2\widetilde u^{\nu_j}|^2}
 {H(\nabla\widetilde u^{\nu_j})^2}\,dx\,dt\\
 &\leq C_T\left(1+\mathbb E
 \|u_0\|_{H^1(\mathbb T^d)}^2\right).
\end{align*}
The representation in
\eqref{vanishing-viscosity-representation-convergence} implies weak
convergence in this time--space Sobolev space, and
Proposition~\ref{lem-strong-normal-limit} identifies the limit as
$\nabla\widetilde u/H(\nabla\widetilde u)$.
In particular, this limit belongs almost surely to $L^2(0,T;H^1(\mathbb T^d;\mathbb R^d))$, and the continuity of spatial differentiation yields, almost surely,
\begin{equation*}
 \nabla\left(\frac{\nabla\widetilde u^{\nu_j}}
 {H(\nabla\widetilde u^{\nu_j})}\right)
 \rightharpoonup
 \nabla\left(\frac{\nabla\widetilde u}
 {H(\nabla\widetilde u)}\right)
 \quad\text{in }
 L^2((0,T)\times\mathbb T^d;\mathbb R^{d\times d}).
\end{equation*}
Moreover, we have 
\begin{equation}\label{normal-matrix-convergence}
\begin{aligned}
 &\left\|
 \left(I-\frac{\nabla \widetilde u^{\nu_j}\otimes\nabla \widetilde u^{\nu_j}}
 {H(\nabla \widetilde u^{\nu_j})^2}\right)
 -\left(I-\frac{\nabla \widetilde u\otimes\nabla \widetilde u}
 {H(\nabla \widetilde u)^2}\right)
 \right\|_{L^2((0,T)\times\mathbb T^d)}\\
 =&\left\|
 \frac{\nabla \widetilde u^{\nu_j}\otimes\nabla \widetilde u^{\nu_j}}{H(\nabla \widetilde u^{\nu_j})^2}
 -\frac{\nabla \widetilde u\otimes\nabla \widetilde u}{H(\nabla \widetilde u)^2}
 \right\|_{L^2((0,T)\times\mathbb T^d)}\\
\leq&2\left\|
 \frac{\nabla \widetilde u^{\nu_j}}{H(\nabla \widetilde u^{\nu_j})}
 -\frac{\nabla \widetilde u}{H(\nabla \widetilde u)}
 \right\|_{L^2((0,T)\times\mathbb T^d)}\longrightarrow0.
\end{aligned}
\end{equation}
Noting that the matrices $
 I-\frac{\nabla\widetilde u^{\nu_j}\otimes\nabla\widetilde u^{\nu_j}}
 {H(\nabla\widetilde u^{\nu_j})^2}$ and $I-\frac{\nabla\widetilde u\otimes\nabla\widetilde u}
 {H(\nabla\widetilde u)^2}$  have Frobenius norm at most $\sqrt d$.
Thus, by the dominated convergence theorem, we have that for every fixed
$\Phi\in L^2((0,T)\times\mathbb T^d;\mathbb R^{d\times d})$,
\begin{align*}
 \left(I-\frac{\nabla \widetilde u^{\nu_j}\otimes\nabla \widetilde u^{\nu_j}}
 {H(\nabla \widetilde u^{\nu_j})^2}\right)\Phi\longrightarrow
 \left(I-\frac{\nabla \widetilde u\otimes\nabla \widetilde u}
 {H(\nabla \widetilde u)^2}\right)\Phi
 \quad\text{strongly in }L^2((0,T)\times\mathbb T^d).
\end{align*}
Together with $\widetilde Y_D^{\nu_j}\rightharpoonup \widetilde Y_D$ in
$L^2((0,T)\times\mathbb T^d)$, it follows that 
\begin{align*}
 &\int_0^T\!\int_{\mathbb T^d}
 \nabla\left(\frac{\nabla \widetilde u}{H(\nabla \widetilde u)}\right):\Phi\,dx\,dt=\lim_{j\to\infty}\int_0^T\!\int_{\mathbb T^d}
 \nabla\left(\frac{\nabla \widetilde u^{\nu_j}}{H(\nabla \widetilde u^{\nu_j})}\right):\Phi\,dx\,dt\\
=&\lim_{j\to\infty}\int_0^T\!\int_{\mathbb T^d}
 \left[\left(I-\frac{\nabla \widetilde u^{\nu_j}\otimes\nabla \widetilde u^{\nu_j}}
 {H(\nabla \widetilde u^{\nu_j})^2}\right)
 \frac{D^2\widetilde u^{\nu_j}}{H(\nabla \widetilde u^{\nu_j})}\right]:\Phi\,dx\,dt\\
 =&\lim_{j\to\infty}\int_0^T\!\int_{\mathbb T^d}
 \frac{D^2\widetilde u^{\nu_j}}{H(\nabla \widetilde u^{\nu_j})}:
 \left(I-\frac{\nabla \widetilde u^{\nu_j}\otimes\nabla \widetilde u^{\nu_j}}
 {H(\nabla \widetilde u^{\nu_j})^2}\right)\Phi\,dx\,dt\\
 =&\int_0^T\!\int_{\mathbb T^d}\widetilde Y_D:
 \left(I-\frac{\nabla \widetilde u\otimes\nabla \widetilde u}
 {H(\nabla \widetilde u)^2}\right)\Phi\,dx\,dt\\
 =&\int_0^T\!\int_{\mathbb T^d}
 \left(I-\frac{\nabla \widetilde u\otimes\nabla \widetilde u}{H(\nabla \widetilde u)^2}\right)\widetilde Y_D:\Phi\,dx\,dt.
\end{align*}
Consequently, 
\begin{equation*}
 \nabla\left(\frac{\nabla \widetilde u}{H(\nabla \widetilde u)}\right)
 =\left(I-\frac{\nabla \widetilde u\otimes\nabla \widetilde u}{H(\nabla \widetilde u)^2}\right)\widetilde Y_D.
\end{equation*}
Since $\nabla \widetilde u\in L^2((0,T)\times\mathbb T^d)$,
$H(\nabla \widetilde u)$ is finite almost everywhere.
The matrix on the right has strictly positive eigenvalues $1$ and
$H(\nabla \widetilde u)^{-2}$, then its inverse exists almost everywhere.  Therefore, almost everywhere,
\begin{equation}\label{limit-weighted-hessian}
 \frac{D^2\widetilde u}{H(\nabla \widetilde u)}
 :=\left(I-\frac{\nabla \widetilde u\otimes\nabla \widetilde u}{H(\nabla \widetilde u)^2}\right)^{-1}
 \nabla\left(\frac{\nabla \widetilde u}{H(\nabla \widetilde u)}\right).
\end{equation}
This implies that $\widetilde Y_D=D^2\widetilde u/H(\nabla \widetilde u)$ with the explicit representation stated above.
In particular, the field belongs
to $L^2((0,T)\times\mathbb T^d)$ and is symmetric, since $\widetilde Y_D$ is the weak limit of symmetric
matrix fields. 

By \eqref{normal-matrix-convergence} and the uniform boundedness of the
tensor matrices, the same dominated convergence argument gives, for every
$\zeta\in L^2((0,T)\times\mathbb T^d)$,
\begin{align*}
 \zeta\frac{\nabla \widetilde u^{\nu_j}\otimes\nabla \widetilde u^{\nu_j}}
 {H(\nabla \widetilde u^{\nu_j})^2}\longrightarrow
 \zeta\frac{\nabla \widetilde u\otimes\nabla \widetilde u}
 {H(\nabla \widetilde u)^2}
 \quad\text{strongly in }L^2((0,T)\times\mathbb T^d).
\end{align*}
Pairing this strong convergence with the weak convergence of the
weighted Hessians, we obtain
\begin{align*}
 &\lim_{j\to\infty}\int_0^T\!\int_{\mathbb T^d}\zeta
 \frac{\nabla \widetilde u^{\nu_j}\cdot D^2\widetilde u^{\nu_j}\nabla \widetilde u^{\nu_j}}
 {H(\nabla \widetilde u^{\nu_j})^3}\,dx\,dt=\lim_{j\to\infty}\int_0^T\!\int_{\mathbb T^d}
 \frac{D^2\widetilde u^{\nu_j}}{H(\nabla \widetilde u^{\nu_j})}:
 \zeta\frac{\nabla \widetilde u^{\nu_j}\otimes\nabla \widetilde u^{\nu_j}}
 {H(\nabla \widetilde u^{\nu_j})^2}\,dx\,dt\\
 =&\int_0^T\!\int_{\mathbb T^d}\widetilde Y_D:
 \zeta\frac{\nabla \widetilde u\otimes\nabla \widetilde u}
 {H(\nabla \widetilde u)^2}\,dx\,dt=\int_0^T\!\int_{\mathbb T^d}
 \frac{D^2\widetilde u}{H(\nabla \widetilde u)}:
 \zeta\frac{\nabla \widetilde u\otimes\nabla \widetilde u}
 {H(\nabla \widetilde u)^2}\,dx\,dt\\
=&\int_0^T\!\int_{\mathbb T^d}\zeta
 \frac{\nabla \widetilde u\cdot D^2\widetilde u\,\nabla \widetilde u}
 {H(\nabla \widetilde u)^3}\,dx\,dt.
\end{align*}
Thus the correction converges weakly in
$L^2((0,T)\times\mathbb T^d)$ to
\begin{equation*}
 \frac{\nabla \widetilde u}{H(\nabla \widetilde u)}\cdot
 \frac{D^2\widetilde u}{H(\nabla \widetilde u)}
 \frac{\nabla \widetilde u}{H(\nabla \widetilde u)}.
\end{equation*}
This proves
\eqref{geometric-correction-convergence}.
Finally, thanks to weak lower semicontinuity, Fatou's lemma, equality of the
represented laws, and \eqref{target-gradient-curvature-estimate}, we have that 
\begin{align*}
 \widetilde{\mathbb E}\int_0^T\!\int_{\mathbb T^d}|\widetilde Y_D|^2\,dx\,dt\leq\liminf_{j\to\infty}\widetilde{\mathbb E}
 \int_0^T\!\int_{\mathbb T^d}|\widetilde Y_D^{\nu_j}|^2\,dx\,dt
 \leq C_T\left(1+\mathbb E\|u_0\|_{H^1(\mathbb T^d)}^2\right).
\end{align*}
\end{proof}

We next pass to the limit for the curvature drift. 

\begin{proposition}\label{lem-curvature-drift-limit}
Under the hypotheses of
Proposition~\ref{lem-weighted-hessian-correction-limit},
\begin{equation*}
 \widetilde Y_\kappa=H(\nabla \widetilde u)^{1/2}\kappa(\widetilde u)
 \in L^2(\widetilde\Omega\times(0,T)\times\mathbb T^d).
\end{equation*}
Along the entire represented sequence,
\begin{equation}\label{geometric-term-convergences}
 H(\nabla \widetilde u^{\nu_j})\kappa(\widetilde u^{\nu_j})
 \rightharpoonup H(\nabla \widetilde u)\kappa(\widetilde u)
 \quad\text{weakly in }L^1((0,T)\times\mathbb T^d)
 \quad\text{almost surely}.
\end{equation}
\end{proposition}

\begin{proof}
We work on the new probability space as before.
Using the weak $L^2(0,T;H^1(\mathbb T^d))$ convergence of the normalized
gradients to see that 
\begin{equation*}
 \kappa(\widetilde u^{\nu_j})=\operatorname{div}
 \left(\frac{\nabla \widetilde u^{\nu_j}}{H(\nabla \widetilde u^{\nu_j})}\right)
 \rightharpoonup\operatorname{div}
 \left(\frac{\nabla \widetilde u}{H(\nabla \widetilde u)}\right)=\kappa(\widetilde u)
 \quad\text{in }L^2((0,T)\times\mathbb T^d).
\end{equation*}
For every $\zeta\in L^\infty((0,T)\times\mathbb T^d)$,
\eqref{strong-square-root-mobility} and the preceding weak convergence provide that 
\begin{align*}
 &\int_0^T\!\int_{\mathbb T^d}\widetilde Y_\kappa\zeta\,dx\,dt=\lim_{j\to\infty}\int_0^T\!\int_{\mathbb T^d}
 \widetilde Y_\kappa^{\nu_j}\zeta\,dx\,dt=\lim_{j\to\infty}\int_0^T\!\int_{\mathbb T^d}
 \kappa(\widetilde u^{\nu_j})H(\nabla \widetilde u^{\nu_j})^{1/2}\zeta\,dx\,dt\\
 =&\lim_{j\to\infty}\int_0^T\!\int_{\mathbb T^d}
 \kappa(\widetilde u^{\nu_j})H(\nabla \widetilde u)^{1/2}\zeta\,dx\,dt=\int_0^T\!\int_{\mathbb T^d}
 \kappa(\widetilde u)H(\nabla \widetilde u)^{1/2}\zeta\,dx\,dt.
\end{align*}
Hence
\begin{equation*}
 \widetilde Y_\kappa=H(\nabla \widetilde u)^{1/2}\kappa(\widetilde u)
 \in L^2((0,T)\times\mathbb T^d),
\end{equation*}
which identifies the weighted-curvature coordinate and its weak
$L^2((0,T)\times\mathbb T^d)$
limit.

For every $\zeta\in L^\infty((0,T)\times\mathbb T^d)$,
\eqref{strong-square-root-mobility} implies
$\zeta H(\nabla \widetilde u^{\nu_j})^{1/2}\to\zeta H(\nabla \widetilde u)^{1/2}$
strongly in $L^2((0,T)\times\mathbb T^d)$.
Combining this with $\widetilde Y_\kappa^{\nu_j}\rightharpoonup \widetilde Y_\kappa$ weakly in
$L^2((0,T)\times\mathbb T^d)$ yields
\begin{align*}
 &\lim_{j\to\infty}\int_0^T\!\int_{\mathbb T^d}
 \zeta H(\nabla \widetilde u^{\nu_j})\kappa(\widetilde u^{\nu_j})\,dx\,dt=\lim_{j\to\infty}\int_0^T\!\int_{\mathbb T^d}
 \zeta H(\nabla \widetilde u^{\nu_j})^{1/2}\widetilde Y_\kappa^{\nu_j}\,dx\,dt\\
 &=\int_0^T\!\int_{\mathbb T^d}
 \zeta H(\nabla \widetilde u)^{1/2}\widetilde Y_\kappa\,dx\,dt=\int_0^T\!\int_{\mathbb T^d}
 \zeta H(\nabla \widetilde u)\kappa(\widetilde u)\,dx\,dt.
\end{align*}
This proves
\eqref{geometric-term-convergences} in the weak topology of
$L^1((0,T)\times\mathbb T^d)$.

Thanks to the equality of represented laws, Fatou's lemma and
\eqref{target-gradient-curvature-estimate}, we have 
\begin{align*}
 &\widetilde{\mathbb E}\int_0^T\!\int_{\mathbb T^d}|\widetilde Y_\kappa|^2\,dx\,dt\leq\liminf_{j\to\infty}\widetilde{\mathbb E}
 \int_0^T\!\int_{\mathbb T^d}|\widetilde Y_\kappa^{\nu_j}|^2\,dx\,dt
 \leq C_T\left(1+\mathbb E\|u_0\|_{H^1(\mathbb T^d)}^2\right).
\end{align*}
The weighted-Hessian integrability was established in
Proposition~\ref{lem-weighted-hessian-correction-limit}.
\end{proof}

\subsection{Existence by the vanishing-viscosity limit}
\label{subsec-limit-equation}

With the help of the above preparation, we are ready to construct energy martingale solutions. 

\begin{theorem}\label{thm-vanishing-viscosity-existence}
Fix an integer $d\geq2$, $T\in(0,\infty)$, $\varepsilon>0$, and
$K\in\mathbb N$, and suppose
Assumption~\ref{constant-noise-covariance} holds.
Let $\Lambda$ be a Borel probability measure on $H^1(\mathbb T^d)$
satisfying
\[
 \int_{H^1(\mathbb T^d)}\|z\|_{H^1(\mathbb T^d)}^2\,\Lambda(dz)<\infty.
\]
Then there exist a filtered probability space
$(\widetilde\Omega,\widetilde{\mathcal F},
(\widetilde{\mathcal G}_t)_{t\in[0,T]},\widetilde{\mathbb P})$
satisfying the usual conditions, a standard $K$-dimensional
$(\widetilde{\mathcal G}_t)$-Brownian motion
$\widetilde\beta=(\widetilde\beta^k)_{k=1,\ldots,K}$, and an
energy martingale solution $\widetilde u$ of \eqref{SPDE-Ito}
in the sense of Definition~\ref{def-energy-martingale-solution}, with
\[
 \mathcal L_{\widetilde{\mathbb P}}(\widetilde u(0))=\Lambda.
\]
\end{theorem}

\begin{proof}
We construct the solution by vanishing viscosity and prove
\eqref{theorem-area-bound} and \eqref{theorem-weighted-hessian-bound},
with $u$ and $\mathbb E$ replaced by $\widetilde u$ and
$\widetilde{\mathbb E}$. Recall the almost surely convergent representation
$\widetilde Z^j\to\widetilde Z$ from
Proposition~\ref{prop-vanishing-viscosity-representation}, constructed
from viscous solutions with initial law $\Lambda$.

\medskip
\noindent\textit{Step 1. Martingale residuals.}
We begin the martingale argument by defining the limiting residual. On the representation space, let $(\widetilde{\mathcal G}_t)$
be the usual augmentation of the raw canonical filtration
$$
 \widetilde{\mathcal G}_t^0
 :=\sigma\{\widetilde u(r),\widetilde\beta(r):0\leq r\leq t\}.
$$
For every $\varphi\in C^\infty(\mathbb T^d)$, we define
\begin{equation}\label{limit-residual-definition}
 \widetilde M_t^\varphi:=\langle \widetilde u(t)-\widetilde u(0),\varphi\rangle
 -\int_0^t\!\int_{\mathbb T^d}H(\nabla \widetilde u)\kappa(\widetilde u)\varphi\,dx\,dr-\frac{\varepsilon\Lambda_K}{8}\int_0^t\!\int_{\mathbb T^d}
 \frac{\nabla \widetilde u\cdot D^2\widetilde u\,\nabla \widetilde u}{H(\nabla \widetilde u)^3}
 \varphi\,dx\,dr.
\end{equation}
Since $\widetilde u$ is adapted and continuous in $H^{-m}(\mathbb T^d)$, it is
predictable relative to $(\widetilde{\mathcal G}_t)$.
It belongs to $H^1(\mathbb T^d)$ almost everywhere, so the Lusin--Souslin argument
used in Step 2 of Theorem~\ref{thm-fixed-nu-identification} gives an
$H^1(\mathbb T^d)$-valued predictable representative.  Thus its gradient, normalized
gradient, and noise coefficients have predictable versions.
The normalized gradient has an $H^1(\mathbb T^d;\mathbb R^d)$-valued
predictable representative, so its weak spatial derivatives and
$\kappa(\widetilde u)$ also have predictable versions. 
The explicit formula \eqref{limit-weighted-hessian}, using the Borel
matrix inversion map, then yields predictable versions of the weighted
Hessian and both geometric drift terms.

Their integrability has already been established in
Propositions~\ref{lem-strong-normal-limit},
\ref{lem-weighted-hessian-correction-limit}, and
\ref{lem-curvature-drift-limit}. 
For the noise, by Assumption~\ref{constant-noise-covariance}, we have
\[
 \widetilde{\mathbb E}\int_0^T\sum_{k=1}^K
 \|G_k(\widetilde u(r))\|_{L^2(\mathbb T^d)}^2\,dr
 =\Lambda_K\,\widetilde{\mathbb E}\int_0^T
 \mathcal A(\widetilde u(r))\,dr<\infty.
\]
Consequently, the drift integrals in \eqref{limit-residual-definition}
and the stochastic integrals used below are well defined, and
$\widetilde M^\varphi$ is continuous and adapted to
$(\widetilde{\mathcal G}_t)$.

\medskip
\noindent\textit{Step 2. Convergence of residuals and noise covariations.}
For each represented $\tilde{u}^{\nu_j}$, we define
\begin{equation}\label{viscous-residual-definition}
\begin{aligned}
 \widetilde M_t^{\nu_j,\varphi}:={}&\langle \widetilde u^{\nu_j}(t)-\widetilde u^{\nu_j}(0),\varphi\rangle
 -\nu_j\int_0^t\langle\Delta \widetilde u^{\nu_j}(r),\varphi\rangle\,dr-\int_0^t\!\int_{\mathbb T^d}
 H(\nabla \widetilde u^{\nu_j})\kappa(\widetilde u^{\nu_j})\varphi\,dx\,dr\\
 &-\frac{\varepsilon\Lambda_K}{8}\int_0^t\!\int_{\mathbb T^d}
 \frac{\nabla \widetilde u^{\nu_j}\cdot D^2\widetilde u^{\nu_j}\nabla \widetilde u^{\nu_j}}
 {H(\nabla \widetilde u^{\nu_j})^3}\varphi\,dx\,dr.
\end{aligned}
\end{equation}
We first pass to the limits for the drift terms. Propositions~\ref{lem-weighted-hessian-correction-limit} and
\ref{lem-curvature-drift-limit}, tested against
$\mathbf1_{[0,t]}\varphi$, give convergence of the two nonviscous drift
integrals in \eqref{viscous-residual-definition} for every $t\in[0,T]$.
We now show that these convergences are uniform in $t$.

For $0\leq s<t\leq T$, we obtain from H\"older's inequality that
\begin{align*}
 &\left|\int_s^t\!\int_{\mathbb T^d}
 H(\nabla\widetilde u^{\nu_j})\kappa(\widetilde u^{\nu_j})
 \varphi\,dx\,dr\right|\\
 \leq&\|\varphi\|_{L^\infty(\mathbb T^d)}
 \left(\int_s^t\!\int_{\mathbb T^d}
 H(\nabla\widetilde u^{\nu_j})\,dx\,dr\right)^{1/2}
 \left\|H(\nabla\widetilde u^{\nu_j})^{1/2}
 \kappa(\widetilde u^{\nu_j})\right\|_{L^2((0,T)\times\mathbb T^d)}.
\end{align*}
The last factor is bounded uniformly in $j$ by weak $L^2((0,T)\times\mathbb T^d)$
convergence.  Proposition~\ref{prop-gradient-compactness-plan} and
the Lipschitz continuity of $H$ give
$H(\nabla\widetilde u^{\nu_j})\to H(\nabla\widetilde u)$
strongly in $L^1((0,T)\times\mathbb T^d)$.
For all sufficiently large $j$, the $L^1((0,T)\times\mathbb T^d)$ difference
from the limit is uniformly small.  By absolute continuity of the
integral for the limit and the finitely many remaining terms, we obtain
$$
 \lim_{h\downarrow0}\ \sup_{j\geq1}\
 \sup_{\substack{0\leq s<t\leq T\\t-s\leq h}}
 \int_s^t\!\int_{\mathbb T^d}
 H(\nabla\widetilde u^{\nu_j})\,dx\,dr=0.
$$
Thus the curvature terms are equicontinuous.
For the correction term, using the Cauchy--Schwarz inequality and the uniform $L^2((0,T)\times\mathbb T^d)$
bound in Proposition~\ref{lem-weighted-hessian-correction-limit} to see that 
$$
 \left|\int_s^t\!\int_{\mathbb T^d}
 \frac{\nabla\widetilde u^{\nu_j}\cdot
 D^2\widetilde u^{\nu_j}\nabla\widetilde u^{\nu_j}}
 {H(\nabla\widetilde u^{\nu_j})^3}\varphi\,dx\,dr\right|
 \leq C_{\omega,\varphi}|t-s|^{1/2},
$$
where $C_{\omega,\varphi}$ is independent of $j$.

Let $\eta>0$ be a fixed arbitrary number.  The preceding equicontinuity and continuity
of the limiting time integral allow us to choose a finite partition
$0=t_0<\cdots<t_L=T$ such that both the approximate and limiting
curvature terms vary by at most $\eta$ on each subinterval,
uniformly in $j$.  Hence
\begin{align*}
 &\sup_{t\in[0,T]}\left|\int_0^t\!\int_{\mathbb T^d}
 \bigl[H(\nabla\widetilde u^{\nu_j})\kappa(\widetilde u^{\nu_j})
 -H(\nabla\widetilde u)\kappa(\widetilde u)\bigr]\varphi\,dx\,dr\right|\\
\leq&\max_{0\leq\ell\leq L}
 \left|\int_0^{t_\ell}\!\int_{\mathbb T^d}
 \bigl[H(\nabla\widetilde u^{\nu_j})\kappa(\widetilde u^{\nu_j})
 -H(\nabla\widetilde u)\kappa(\widetilde u)\bigr]\varphi\,dx\,dr\right|
 +2\eta.
\end{align*}
The finite maximum tends to zero by the convergence at each fixed time.
Letting $\eta\downarrow0$ gives
$$
 \sup_{t\in[0,T]}\left|\int_0^t\!\int_{\mathbb T^d}
 \bigl[H(\nabla\widetilde u^{\nu_j})\kappa(\widetilde u^{\nu_j})
 -H(\nabla\widetilde u)\kappa(\widetilde u)\bigr]\varphi\,dx\,dr\right|
 \longrightarrow0
 \quad\text{almost surely}.
$$
The same finite-partition argument applies to the correction terms:
\begin{align*}
 &\sup_{t\in[0,T]}\Bigg|
 \int_0^t\!\int_{\mathbb T^d}
 \frac{\nabla\widetilde u^{\nu_j}\cdot
 D^2\widetilde u^{\nu_j}\nabla\widetilde u^{\nu_j}}
 {H(\nabla\widetilde u^{\nu_j})^3}\varphi\,dx\,dr
 -\int_0^t\!\int_{\mathbb T^d}
 \frac{\nabla\widetilde u\cdot D^2\widetilde u\,\nabla\widetilde u}
 {H(\nabla\widetilde u)^3}\varphi\,dx\,dr
 \Bigg|\longrightarrow0
 \quad\text{almost surely}.
\end{align*}

For the viscosity term, integration by parts and
\eqref{pathwise-gradient-curvature-bounds} show that 
\begin{align*}
 &\sup_{0\leq t\leq T}\left|
 \nu_j\int_0^t\langle\Delta\widetilde u^{\nu_j},\varphi\rangle\,dr
 \right|\leq\nu_j\sqrt T\,\|\nabla\varphi\|_{L^2(\mathbb T^d)}
 \|\nabla\widetilde u^{\nu_j}\|_{L^2((0,T)\times\mathbb T^d)}
 \leq C_{\omega,\varphi,T}\nu_j\longrightarrow0
 \quad\text{almost surely}.
\end{align*}
Together with the path convergence in
\eqref{vanishing-viscosity-representation-convergence}, this yields
\begin{equation}\label{martingale-residual-pathwise-convergence}
 \widetilde M^{\nu_j,\varphi}\longrightarrow\widetilde M^\varphi
 \quad\text{almost surely in }C([0,T]).
\end{equation}

By \eqref{strong-square-root-mobility}, for $k=1,\ldots,K$,
\[
 \langle G_k(\widetilde u^{\nu_j}),\varphi\rangle
 \longrightarrow\langle G_k(\widetilde u),\varphi\rangle
\]
strongly in $L^2(\widetilde\Omega\times(0,T))$ and almost surely
in $L^2(0,T)$.  Hence, using H\"older's inequality, we obtain
\begin{align*}
 \sqrt\varepsilon\int_0^\cdot
 \langle G_k(\widetilde u^{\nu_j}(r)),\varphi\rangle\,dr&\longrightarrow\sqrt\varepsilon\int_0^\cdot
 \langle G_k(\widetilde u(r)),\varphi\rangle\,dr,
 \qquad k=1,\ldots,K,\\
 \varepsilon\sum_{k=1}^K\int_0^\cdot
 \langle G_k(\widetilde u^{\nu_j}(r)),\varphi\rangle
 \langle G_k(\widetilde u^{\nu_j}(r)),\psi\rangle\,dr&\longrightarrow\varepsilon\sum_{k=1}^K\int_0^\cdot
 \langle G_k(\widetilde u(r)),\varphi\rangle
 \langle G_k(\widetilde u(r)),\psi\rangle\,dr,
\end{align*}
both in $L^1(\widetilde\Omega;C([0,T]))$, for every smooth
$\varphi,\psi\in C^{\infty}(\mathbb{T}^d)$.

\medskip
\noindent\textit{Step 3. Identification of the  stochastic integral.}
Recall that $\mathcal X$ is defined by
\eqref{joint-tightness-space}. Let $v$ and $w$ denote the solution coordinate and Brownian coordinate of
$z\in\mathcal X$, respectively.  On the Borel set where
$v\in L^2(0,T;H^2(\mathbb T^d))$, define
\begin{align*}
 \mathscr R_j^\varphi(z)(t)
 &:=\langle v(t)-v(0),\varphi\rangle
 -\int_0^t\langle b^{\nu_j}(v(r)),\varphi\rangle\,dr,\\
 \mathscr A_k^\varphi(z)(t)
 &:=\sqrt\varepsilon\int_0^t
 \langle G_k(v(r)),\varphi\rangle\,dr,\qquad k=1,\ldots,K,\\
 \mathscr C^{\varphi,\psi}(z)(t)
 &:=\varepsilon\sum_{k=1}^K\int_0^t
 \langle G_k(v(r)),\varphi\rangle
 \langle G_k(v(r)),\psi\rangle\,dr.
\end{align*}
By the Borel reconstruction in
Lemma~\ref{lem-vanishing-viscosity-coordinate-identities} and the
continuity of the coefficient maps on $H^2(\mathbb T^d)$, these maps, extended by
zero outside that set, are Borel maps
\[
 \mathcal X\longrightarrow C([0,T]).
\]
Moreover, the past-coordinate map
\[
 \mathcal X\ni z\longmapsto
 (v|_{[0,s]},w|_{[0,s]})
 \in C([0,s];H^{-m}(\mathbb T^d))\times C([0,s];\mathbb R^K)
\]
is continuous.  These Borel maps allow us to transfer the moment
bounds and martingale identities to the representation space by
equality of joint laws.

We now upgrade the residual convergence to convergence in
expectation.  The Burkholder--Davis--Gundy inequality on the original
spaces, \eqref{second-moment-area-estimate}, and equality of joint laws
give
\[
 \sup_{j\geq1}\widetilde{\mathbb E}\sup_{0\leq t\leq T}
 |\widetilde M_t^{\nu_j,\varphi}|^4
 \leq C_{\varphi,T}\sup_{j\geq1}\widetilde{\mathbb E}\sup_{0\leq t\leq T}
 \mathcal A(\widetilde u^{\nu_j}(t))^2
 \leq C_{\varphi,T}\left(1+\mathbb E\|u_0\|_{H^1(\mathbb T^d)}^2\right).
\]
Fatou's lemma and Vitali's theorem, together with
\eqref{martingale-residual-pathwise-convergence}, therefore give
\begin{equation}\label{martingale-strong-convergence}
 \widetilde M^{\nu_j,\varphi}\longrightarrow\widetilde M^\varphi
 \quad\text{strongly in }L^2(\widetilde\Omega;C([0,T])).
\end{equation}
Similarly, the Brownian coordinates converge strongly in
$L^2(\widetilde\Omega;C([0,T];\mathbb R^K))$, by their almost sure
path convergence and uniform fourth moments.
Consequently, the products of residuals and Brownian coordinates
used below converge in
$L^1(\widetilde\Omega;C([0,T]))$, allowing passage to expectations
in the joint martingale identities.

Fix $0\leq s<t\leq T$ and a bounded continuous functional $\Psi$ of
$(\widetilde u,\widetilde{\beta})$ restricted to $[0,s]$.  Set
\begin{equation*}
 \widetilde\Psi^{\nu_j}=\Psi(\widetilde u^{\nu_j}|_{[0,s]},\widetilde{\beta}^{\nu_j}|_{[0,s]}),\qquad
 \widetilde\Psi^0=\Psi(\widetilde u|_{[0,s]},\widetilde{\beta}|_{[0,s]}).
\end{equation*}
The path convergence in
\eqref{vanishing-viscosity-representation-convergence} implies
$\widetilde\Psi^{\nu_j}\to\widetilde\Psi^0$ almost surely.  On each original stochastic basis the
corresponding past functional is measurable at time $s$.  After transfer by
equality of joint laws, the expected increment from $s$ to $t$ of each of
the following represented processes, multiplied by $\widetilde\Psi^{\nu_j}$, is zero:
\begin{align*}
 &\widetilde M^{\nu_j,\varphi},\quad \widetilde{\beta}^{\nu_j,k},\quad\widetilde{\beta}^{\nu_j,k}\widetilde{\beta}^{\nu_j,\ell}-\delta_{k\ell}(\cdot),\\
 &\widetilde M^{\nu_j,\varphi}\widetilde{\beta}^{\nu_j,k}
 -\sqrt\varepsilon\int_0^{\cdot}
 \langle f_kH(\nabla \widetilde u^{\nu_j}(r))^{1/2},\varphi\rangle\,dr,\\
 &\widetilde M^{\nu_j,\varphi}\widetilde M^{\nu_j,\psi}
 -\varepsilon\sum_{k=1}^K\int_0^{\cdot}
 \langle f_kH(\nabla \widetilde u^{\nu_j}(r))^{1/2},\varphi\rangle
 \langle f_kH(\nabla \widetilde u^{\nu_j}(r))^{1/2},\psi\rangle\,dr.
\end{align*}
The strong convergences proved above permit passage to expectations.
A monotone-class argument then replaces $\widetilde\Psi^0$ by any bounded
$\widetilde{\mathcal G}_s^0$-measurable random variable.  Thus the limiting processes
are first martingales in $(\widetilde{\mathcal G}_t^0)$.  The moment bounds and
compensator estimates above give them integrable time suprema, so their
continuous paths are also $L^1(\widetilde\Omega)$ continuous.  As in
Step 3 of the proof of Theorem~\ref{thm-fixed-nu-identification}, these identities also hold in the
right-continuous completed filtration $(\widetilde{\mathcal G}_t)$.

It follows from the convergence of the Brownian paths that $\widetilde{\beta}(0)=0$ almost surely.
The identities for
$\widetilde{\beta}^k\widetilde{\beta}^\ell-\delta_{k\ell}(\cdot)$ and the multidimensional
L\'evy characterization show that $\widetilde{\beta}$ is a standard
$K$-dimensional Brownian motion.  The limiting residuals are continuous
square-integrable martingales by the moment bounds above, Fatou's
lemma, and the martingale identities just established.  The remaining product
identities give
\begin{equation}\label{limiting-covariation}
 \langle \widetilde M^\varphi,\widetilde M^\psi\rangle_t
 =\varepsilon\sum_{k=1}^K\int_0^t
 \langle f_kH(\nabla \widetilde u(r))^{1/2},\varphi\rangle
 \langle f_kH(\nabla \widetilde u(r))^{1/2},\psi\rangle\,dr,
\end{equation}
and
\begin{equation*}
 \langle \widetilde M^\varphi,\widetilde{\beta}^k\rangle_t
 =\sqrt\varepsilon\int_0^t
 \langle f_kH(\nabla \widetilde u(r))^{1/2},\varphi\rangle\,dr.
\end{equation*}

We define
\begin{equation*}
 \widetilde N_t^\varphi=\widetilde M_t^\varphi-\sqrt\varepsilon\sum_{k=1}^K\int_0^t
 \langle f_kH(\nabla \widetilde u(r))^{1/2},\varphi\rangle\,d\widetilde{\beta}_r^k.
\end{equation*}
Then $\widetilde N^\varphi$ is a continuous square-integrable martingale starting at
zero.  The preceding covariation identities yield
\begin{align*}
 \langle \widetilde N^\varphi\rangle_t
 ={}&\langle \widetilde M^\varphi\rangle_t
 +\varepsilon\sum_{k=1}^K\int_0^t
 |\langle f_kH(\nabla \widetilde u(r))^{1/2},\varphi\rangle|^2\,dr\\
 &-2\sqrt\varepsilon\sum_{k=1}^K\int_0^t
 \langle f_kH(\nabla \widetilde u(r))^{1/2},\varphi\rangle
 \,d\langle \widetilde M^\varphi,\widetilde{\beta}^k\rangle_r=0.
\end{align*}
Hence $\widetilde N^\varphi=0$, which identifies $\widetilde M^\varphi$ as the stochastic
integral against the retained Brownian motion $\widetilde{\beta}$.

Substituting \eqref{limit-residual-definition}, we obtain
\eqref{weak-martingale-formulation} for every smooth test function,
almost surely for all $t\in[0,T]$.

\medskip
\noindent\textit{Step 4. Identification of the initial law.}
Since $\widetilde u^{\nu_j}(0)=\widetilde u_0^j$, every
$\widetilde u_0^j$ has law $\Lambda$, and
\eqref{vanishing-viscosity-representation-convergence} identifies
$\widetilde u(0)$ with the almost sure weak $H^1(\mathbb T^d)$ limit
$\widetilde u_0$ of these initial coordinates, we obtain
\begin{equation}\label{theorem-initial-law-identification}
 \widetilde{\mathbb P}\circ \widetilde u(0)^{-1}=\widetilde{\mathbb P}\circ \widetilde u_0^{-1}=\Lambda.
\end{equation}

\medskip
\noindent\textit{Step 5. Energy bounds and path regularity.}
All viscous initial coordinates and the limiting initial datum have
law $\Lambda$ by Step 4.  Their common
second moment is therefore
\begin{equation*}
 \widetilde{\mathbb E}\|\widetilde u_0\|_{H^1(\mathbb T^d)}^2
 =\int_{H^1(\mathbb T^d)}\|z\|_{H^1(\mathbb T^d)}^2\,\Lambda(dz)<\infty.
\end{equation*}
The Cauchy--Schwarz inequality in space also gives
\begin{align*}
 \widetilde{\mathbb E}\mathcal A(\widetilde u_0)^2
 &\lesssim1+\int_{H^1(\mathbb T^d)}
 \|z\|_{H^1(\mathbb T^d)}^2\,\Lambda(dz).
\end{align*}
Substituting the initial-moment identity into
\eqref{second-moment-area-estimate},
\eqref{target-gradient-curvature-estimate}, and \eqref{uniform-u-estimate},
and using equality of represented laws, yields their uniform bounds with
right-hand side
$C_T\bigl(1+\int_{H^1(\mathbb T^d)}\|z\|_{H^1(\mathbb T^d)}^2\,\Lambda(dz)\bigr)$.
The additive constant accounts for the finite volume of $\mathbb T^d$
and the data-independent terms in those estimates. $C_T$ is independent
of $\nu_j$ and $\Lambda$.

Since $H$ is Lipschitz,
\eqref{target-strong-gradient-convergence} with $p=1$ yields, almost surely,
\begin{align*}
 \int_0^T|\mathcal A(\widetilde u^{\nu_j}(t))-\mathcal A(\widetilde u(t))|\,dt
 &\leq\int_0^T\!\int_{\mathbb T^d}
 |\nabla \widetilde u^{\nu_j}-\nabla \widetilde u|\,dx\,dt\longrightarrow0.
\end{align*}
For the area supremum, take a subsequence realizing the lower limit of
the suprema and then an almost-everywhere-in-time convergent subsequence
of the displayed $L^1(0,T)$ convergence.  Lower semicontinuity of the
weak $L^2((0,T)\times\mathbb T^d)$ limits in
Propositions~\ref{lem-weighted-hessian-correction-limit} and
\ref{lem-curvature-drift-limit} then give, almost surely,
\begin{align*}
 \left(\operatorname*{ess\,sup}_{t<T}\mathcal A(\widetilde u(t))\right)^2
 &\leq\liminf_{j\to\infty}
 \left(\sup_{0\leq t\leq T}\mathcal A(\widetilde u^{\nu_j}(t))\right)^2,\\
 \left(\int_0^T\!\int_{\mathbb T^d}
 H(\nabla \widetilde u)|\kappa(\widetilde u)|^2\,dx\,dt\right)^2
 &\leq\liminf_{j\to\infty}
 \left(\int_0^T\!\int_{\mathbb T^d}
 H(\nabla \widetilde u^{\nu_j})|\kappa(\widetilde u^{\nu_j})|^2\,dx\,dt\right)^2,\\
 \int_0^T\!\int_{\mathbb T^d}
 \left|\frac{D^2\widetilde u}{H(\nabla \widetilde u)}\right|^2\,dx\,dt
 &\leq\liminf_{j\to\infty}
 \int_0^T\!\int_{\mathbb T^d}
 \left|\frac{D^2\widetilde u^{\nu_j}}{H(\nabla \widetilde u^{\nu_j})}\right|^2\,dx\,dt.
\end{align*}
Fatou's lemma, \eqref{second-moment-area-estimate} and
\eqref{target-gradient-curvature-estimate} therefore yield
\begin{align*}
 &\widetilde{\mathbb E}\left[\operatorname*{ess\,sup}_{t<T}
 \mathcal A(\widetilde u(t))^2\right]
 +\widetilde{\mathbb E}\left(\int_0^T\!\int_{\mathbb T^d}
 H(\nabla \widetilde u)|\kappa(\widetilde u)|^2\,dx\,dt\right)^2\\
 &\quad+\widetilde{\mathbb E}\int_0^T\!\int_{\mathbb T^d}
 \left|\frac{D^2\widetilde u}{H(\nabla \widetilde u)}\right|^2\,dx\,dt
 \leq C_T\bigl(1+\int_{H^1(\mathbb T^d)}\|z\|_{H^1(\mathbb T^d)}^2\,\Lambda(dz)\bigr).
\end{align*}
This proves \eqref{theorem-area-bound} and the weighted-Hessian part of
\eqref{theorem-weighted-hessian-bound}.

The uniform $L^2(\widetilde\Omega\times(0,T)\times\mathbb T^d)$ gradient bound in
\eqref{target-gradient-curvature-estimate}, together with the strong
$L^1(\widetilde\Omega\times(0,T)\times\mathbb T^d)$ convergence, identifies every weakly convergent subsequence with
$\nabla \widetilde u$.  Thus the entire represented sequence converges weakly in
this $L^2(\widetilde\Omega\times(0,T)\times\mathbb T^d)$ space.  For every measurable $E\subset(0,T)$,
\begin{equation*}
\begin{aligned}
 \int_E\widetilde{\mathbb E}\|\nabla \widetilde u(t)\|_{L^2(\mathbb T^d)}^2\,dt
 &=\|\mathbf1_E\nabla \widetilde u\|_{L^2(\widetilde\Omega\times(0,T)\times\mathbb T^d)}^2\\
 &\leq\liminf_{j\to\infty}
 \int_E\widetilde{\mathbb E}\|\nabla \widetilde u^{\nu_j}(t)\|_{L^2(\mathbb T^d)}^2\,dt\\
 &\leq C_T\bigl(1+\int_{H^1(\mathbb T^d)}\|z\|_{H^1(\mathbb T^d)}^2\,\Lambda(dz)\bigr)|E|.
\end{aligned}
\end{equation*}
It follows from the characterization of $L^\infty(0,T)$ that
\begin{equation*}
 \operatorname*{ess\,sup}_{t<T}
 \widetilde{\mathbb E}\|\nabla \widetilde u(t)\|_{L^2(\mathbb T^d)}^2
 \leq C_T\bigl(1+\int_{H^1(\mathbb T^d)}\|z\|_{H^1(\mathbb T^d)}^2\,\Lambda(dz)\bigr).
\end{equation*}
By \eqref{uniform-u-estimate}, equality of represented laws,
\eqref{vanishing-viscosity-representation-convergence}, and Fatou's lemma,
\begin{align*}
 \widetilde{\mathbb E}\|\widetilde u\|_{L^2((0,T)\times\mathbb T^d)}^2
 \leq\liminf_{j\to\infty}
 \widetilde{\mathbb E}\|\widetilde u^{\nu_j}\|_{L^2((0,T)\times\mathbb T^d)}^2\leq C_T\left(1+\int_{H^1(\mathbb T^d)}
 \|z\|_{H^1(\mathbb T^d)}^2\,\Lambda(dz)\right).
\end{align*}
Together with the gradient bound above, this proves
$\widetilde u\in L^2(\widetilde\Omega\times(0,T);H^1(\mathbb T^d))$. This completes the proof.

\end{proof}

\subsection*{Declaration of AI-Assisted Technologies}
The first author began studying this problem in 2023, identifying the Stratonovich-to-It\^o correction term and establishing the area estimates, but could not pass to the limit in the correction term. In 2026, the two authors resumed the project and established the $H^1$-estimate, while this difficulty remained unresolved. After extensive iterative discussions with OpenAI's ChatGPT 5.6 Pro, this obstacle was finally resolved. In particular, the central idea contributed by AI was subsequently rewritten by the authors and incorporated into Propositions \ref{lem-strong-normal-limit}, \ref{prop-gradient-compactness-plan}, and \ref{lem-weighted-hessian-correction-limit}. The authors also used AI assistance to iteratively check and revise the manuscript and improve its language. Figure 1 was generated by AI under the authors' instructions.

\bibliographystyle{abbrv}
\bibliography{references}
\end{document}